\documentclass{article}
\usepackage[utf8]{inputenc}
\usepackage[a4paper, total={6in, 9in}]{geometry}
\usepackage{hyperref}
\usepackage{setspace}
\usepackage[all,cmtip]{xy}
\usepackage{hyperref}
\usepackage{tikz}
\usepackage{tikz-cd}

\hypersetup{
    linktoc=all,     
    linkcolor=black,  
}

\usepackage{inputenc, mathtools, amssymb, graphicx}
\usepackage{xcolor}

\newcommand{\C}{\mathbb{C}}
\newcommand{\PP}{\mathbb{P}}
\newcommand{\Z}{\mathbb{Z}}
\newcommand{\R}{\mathbb{R}}

\newcommand{\Hb}{{\cal H}_b}
\newcommand{\eat}[1]{}
\newtheorem{theorem}{Theorem}[section]
\newtheorem{lemma}[theorem]{Lemma}
\newtheorem{prop}[theorem]{Proposition}

\newtheorem{remark}[theorem]{Remark}
\newtheorem{ex}[theorem]{Example}
\newtheorem{defn}[theorem]{Definition}

\newtheorem{assume}[theorem]{Assumption}
\newcommand{\ignore}[1]{}
\newcommand{\mf}[1]{\mathfrak{#1}}

\begin{document}
 
\author{
Bharat Adsul \thanks{Computer Science Department, Indian Institute of Technology, Mumbai, India. \tt{adsul@cse.iitb.ac.in}}
\and
Milind Sohoni \thanks{Computer Science Department, Indian Institute of Technology, Mumbai, and Indian Institute of Technology, Goa, India. \tt{sohoni@cse.iitb.ac.in}}
\and 
K. V. Subrahmanyam \thanks {Chennai Mathematical Institute, Chennai, India. \tt{kv@cmi.ac.in}}
}

\title{Geometric Complexity Theory - Lie Algebraic Methods for Projective Limits of Stable Points} 
\date{31st December, 2021} 

\maketitle
\section*{Abstract}

Let $G$ be a connected reductive algebraic group over ${\mathbb C}$, with Lie algebra ${\cal G}$, acting rationally on a complex vector space $V$, and the corresponding projective space $\PP V$. Let $x\in V$ and let $H\subseteq G$ be its stabilizer and
${\cal H}\subseteq {\cal G}$, its Lie algebra. Our primary objective is to understand the points $[y]$, and their stabilizers, which occur in the vicinity of $[x]$ in $\PP V$. 

Towards this, we construct an explicit Lie algebra action of ${\cal G}$ on a suitably parametrized neighbourhood of $x$. As a consequence, we show that the Lie algebras of the stabilizers of points in the neighbourhood of $x$ are parameterized by subspaces of ${\cal H}$. When ${\cal H}$ is reductive, our results imply that these are in fact, Lie subalgebras of ${\cal H}$. If the orbit of $x$ were closed this would also follow from a celebrated theorem of Luna~\cite{luna1973slices}. 

To construct our Lie algebra action we proceed as follows. We identify the tangent space to the orbit $O_x$ of $x$ with a complement of ${\cal H}$ in ${\cal G}$. Let $N$ be any linear subspace of $V$ that is complementary to $TO_x$, the
tangent space to the orbit $O_x$. We construct an explicit map $\theta : V \times N \rightarrow V$ which captures the action of the nilpotent part of ${\cal H}$ on $V$. We show that the map $\theta $ is intimately connected with the local curvature form of the orbit at that point. We call this data, of a neighbourhood of $x$ with the Lie algebra morphism from ${\cal G}$ to vector fields in this neighbourhood, a local model at $x$.

The action of $G$ on $V$ extends to an action on ${\mathbb P}(V)$. We illustrate the utility of the local model in understanding when $[x]\in {\mathbb P}V$ is in the projective orbit closure of $[y] \in {\mathbb P}V$ via two applications. 

The first application is when $V=Sym^d(X)$, the space of forms in the variables $X:=\{x_1,\ldots,x_k\}$. We consider a form $f$ (as the ``$y$" above) and elements in its orbit given by $f(t)=A(t)\cdot f$, where $A(t)$ is an invertible family in $GL(X)$, with $t \in {\mathbb C}$. 
We express $f(t)$ as $f(t)=t^a g +t^b f_b + \mbox{ higher terms}$, with $a<b$. We call $g\in Sym^d (X)$, the leading term as the limit point (the "$x$" above) and $f_b $ as the direction of approach. Let ${\cal K}$ be the Lie algebra of the stabilizer of $f$ and ${\cal H}$ the Lie algebra of the stabilizer of $g$. 
The local analysis gives us a flattening ${\cal K}_0$ of ${\cal K}$ as a subalgebra of ${\cal H}$, thereby connecting the two stabilizers. There is a natural action of ${\cal H}$ on $\overline{N} = V/TO_g $. We show that ${\cal K}_0$ also stabilizes $\overline{f_b} \in \overline{N}$. We show that there is an $\epsilon$-extension Lie algebra ${\cal K}(\epsilon)$ of ${\cal K}_0$ whose structure constants are incrementally closer to ${\cal K}$. The extendability of ${\cal K}_0$ to ${\cal K}(\epsilon)$ translates into certain Lie algebra cohomology conditions. 

We then specialize to the projective orbit closures $\overline{O(f)}$ of forms $f$ whose $SL(X)$-orbits are affine, and where the form $g$ is such that $\overline{O(g)}$ is of co-dimension $1$ in $\overline{O(f)}$. When $g$ is obtained as the limit of a 1-PS $\lambda (t)$, we examine the {\em tangent of exit}, i.e., the direction along which the path $\lambda (t)\cdot f$ leaves $f$. We show that (i) either ${\cal H}$ has a very simple structure, or (ii) conjugates of the elements of ${\cal K}$, the stabilizer of $f$ also stabilize $g$ and the tangent of exit. This is the triple stabilizer condition. 

The second application looks at the space $Mat_n (\C)$ of $n\times n$-matrices under conjugation by $GL_n$. We show that for a general diagonalizable matrix $X$, the signature of nilpotent matrices in the projective closure of its orbit are determined by the multiplicity data of the spectrum of $X$.

Finally, we formulate the {\em path problem} of finding paths with specific properties from $y$ to its limit points $x$. This connects the algebraic study of the neighborhood of $y$, with the local differential geometry of the orbit thereby allowing an optimization formulation. 

The questions we study are motivated from Geometric complexity theory, an approach to solve the fundamental lower bound problems in complexity theory. This theory was proposed by the second author and Ketan Mulmuley in \cite{mulmuley2001geometric}.

\section{Introduction}
\label{intro}
Let $G$ be a connected reductive {\em algebraic} group over ${\mathbb C}$ and let $\rho: G \rightarrow End(V)$ be a rational representation of $G$ on a complex vector space $V$. We say that $G$ acts on $V$ via the representation $\rho$.
Let $x,y \in V$. Denote by $O(y)$ or $G\cdot y$ the orbit $O(y)=\{\rho(g) \cdot y| g \in G\}$ of $y$ in $V$ under the given representation.  In general the orbit is not Zariski closed, it is only constructible.  
Clearly $G$ also acts on projective space ${\mathbb P}(V)$. For a nonzero $v \in V$, let $[v]$ denote the line corresponding to $v$ in ${\mathbb P}(V)$. 

A fundamental problem of invariant theory is that of obtaining a good description of the $G$ orbit closures of points $y \in V$ and deciding whether a point $x$ belongs to the orbit closure of $y$.
These problems go back to Hilbert and are of fundamental importance in the construction of moduli spaces, see \cite{mumford1994geometric}. The ubiquitous appearance of this problem in many areas of mathematics is also surveyed in the introduction of \cite{popov2009two}. 

 We recall some definitions from Geometric Invariant Theory which we will use throughout  the paper. Standard references for this are \cite{mumford1994geometric} and \cite{dolgachev2003lectures}. Unless stated otherwise, we work over $\C$, the field of complex numbers.

We say a point $v \in V$ is {\it unstable} (or equivalently, is in the {\it nullcone}) for the action of $G$ if $\overline{Gv}$ contains $0$. A point which is not unstable is said to be {\it semistable}. The semistable locus is open in $V$. A point is said to be {\it polystable} if its orbit is closed in the open set of semistable points. If, furthermore, the stabilizer in $G$ of a polystable point has dimension 0, we say it is {\it stable}. We use the same adjectives, viz.,  unstable, polystable and stable for the projective point $[v]$, based on the status of a representative point $v\in [v]$. 

The central motivation for our investigation arises from computational complexity, which is the study of algorithms to compute functions, typically polynomials and forms. 
In Appendix ~\ref{appendix-gct} we give a brief introduction to algebraic complexity theory and the fundamental lower bounds problems therein. We also discuss the {\em geometric complexity theory (GCT)} approach to these problems, as proposed by Mulmuley and Sohoni ~\cite{mulmuley2001geometric},\cite{mulmuley2008geometric}. 
In this formulation, lower bound problems from algebraic complexity translate to determining the membership of special points $x$ within the projective orbit closure of $y$, {\em where $y$ and $x$ have distinctive stabilizers}. In this setting one of the central problems is to determine the orbit closure in ${\mathbb P}(Sym^n(\C^{n^2})^*)$ of the determinant, as a polynomial of degree $n$ in $n^2$ variables and to determine the stabilizers which arise in this orbit closure. This problem, in turn, is closely connected to determining strata in orbit closures, and their classification. 

Determining projective closures is intimately connected to many classical problems. For example, when $y$ is an unstable point in $V$, determining the projective limits of $[y]$ requires an analysis of the stratification of the nullcone as a projective variety, \cite{hesselink1979desingularizations}. This subsumes, e.g., Hilbert's work on classifying null forms, and the combinatorial structure of the variety of nilpotent matrices. 

Another situation where projective closures arise is in the analysis of leading terms. Consider, for example, a finite dimensional $SL(n)$-module $V$, and let $\lambda:  \C^* \rightarrow SL(n) $ be a one parameter subgroup (1-ps) of $SL(n)$. Let the action of $\lambda$ on a point $y\in V$ be given by:
\[ \lambda (t) \cdot y =\sum_j t^j y_j. \]  

Let $x$ be the least degree term i.e., $x=y_k$ where $y_k \neq 0$, $y_i =0$ for all $i<k$. We call $x$ the leading term of $y$ under $\lambda$. Then it is easily seen that 
$[x]$ is in the orbit closure of $[y]$ in projective space ${\mathbb P}(V)$. Whence, projective orbit closures contain all such limit points which arise by the action of 1-ps on $[y]$. Thus, while the orbit of $y$ may be closed in $V$, the orbit $[y]\in \PP V$ is almost always not closed, and picks up points from the nullcone which are intimately connected to it. 

We give a brief review of what is known (to the best of our knowledge) about the affine $G$-orbit closures of points in $V$ and the projective $G$-orbit closure of points in ${\mathbb P}(V)$.
  
The set of unstable points in $V$ is described by the Hilbert-Mumford theory, see \cite{mumford1994geometric, kempf1978instability}. Every unstable point is driven to $0$ by a one parameter subgroup. In \cite{kempf1978instability}, Kempf showed that for every unstable point $y$, there is a canonical parabolic subgroup and an optimal one parameter subgroup which drives $y$ to zero. The corresponding computational question remains far from solved, except in a few special cases, see \cite{burgisser2018alternating}, \cite{garg2019operator}, \cite{derksen2017polynomial}, \cite{ivanyos2017non}. 
    
Hesselink \cite{hesselink1979desingularizations} gave a description of the the strata of the nullcone in terms of Kempf's optimal 1-PS. These strata are locally closed, irreducible and $G$-invariant. Popov \cite{popov2003cone} gave a combinatorial algorithm to determine the strata. 

When $x$ is polystable it follows from Matsushima \cite{matsushima1960espaces} that the stabilizer $H$ of $x$ is reductive.   
In this case Luna's \'{e}tale slice theorem gives information about stabilizers in a Zariski neighbourhoods of $x$. Luna's theorem asserts the existence of an $H$-stable affine variety $N_x \subset V$ passing through $x$ and a $G$-invariant open neighbourhood $U $ of the orbit $Gx$ in $V$, having the following properties. Let $G \times^H N_x$ be the fibre bundle over $G/H$ with fibre $N_x$(see Section 3.2 for definitions). Luna's theorem states that the map 
$\psi: G \times^H N_x \rightarrow U$, sending $[g,v]$ to $gv$ is strongly \'{e}tale. \eat{Strong \'{e}tale-ness means: 
\begin{enumerate}
        \item The map of $G$-quotients, $\psi/G : G \times_H N_x//G \simeq N_x//H \rightarrow U//G$ is \'{e}tale.
        \item Let $G \times_H N_x \rightarrow N_x//H$ be the quotient map. Then $G \times_H N_x \rightarrow  U \times_{U/G} N_x//H$ is an isomorphism.
    \end{enumerate}
    }
It follows from Luna's theorem that the every point in a Zariski neighbourhood of $x$ in $V$ has a stabilizer which is contained in a conjugate of $H$. So, if $H$ does not contain a conjugate of the stabilizer of $y$, then $x$ cannot belong to the {\em affine} orbit closure of $y$.

Of particular interest in GCT is the study of forms $f\in Sym^n (W^*)$ with reductive stabilizers and their orbit closures. In this connection, we have Matsushima's result $\cite{matsushima1960espaces}$ that  
when $f$ has a reductive stabilizer, the $SL(W)$-orbit of $f\in V$ is an affine variety. Whence, the closure in projective space $\overline{O([y])} \backslash O([y])$ is either zero or is pure of codimension one. Thus, there are forms $g_1 ,\ldots, g_k$ whose orbits constitute the closure. 

The action of $GL(W)$ on $Sym^n(W^*)$ can be extended to a natural action of $GL(W) \otimes \C((t))$ on $Sym^n(W^*)$, where $\C((t))$ denotes the Laurent series in $t$. Hilbert \cite{hilbert1893vollen} showed that for any $[g]$ in the orbit closure of $[f]$ there is a $\sigma \in GL(W) \otimes \C((t))$ such that $P:=[f] \circ \sigma $ satisfies
$(P)_{t=0}=[g]$. In \cite{lehmkuhl1989order} (see also \cite[Section 9.4] {burgisser2011overview}, \cite[Section 7.1]{huttenhain2016boundary}), this was reproved. They show the existence of a $\tilde{\sigma} \in End(W) \otimes \C[t]$ such that 
$f \circ \tilde{\sigma} = t^K g \: mod \: (t^{K+1})$. Thus, leading terms under 1-parameter substitutions do lead to all points in the closure, and therefore also the $g_i$ above. 

Although the above results tell us how a form $g$ arises in the closure of a form $f$, how is one to determine these forms? In GCT, however one starts with forms $f$ which are characterized by their stabilizers and one expects the problem to be tractable. Even so very little is known. For example, determining the codimension one components of the boundary of the determinant form is an important open problem.
     
The orbit closure in projective space of a form $f$ as above can in principle be determined by resolving the indeterminacy locus of a certain rational map. This is outlined in  \cite{burgisser2011overview} and makes clear the connections to Geometric invariant theory.  In \cite{huttenhain2017geometric}, the author  works this out in the case of the determinant form of a 3 x 3 matrix and shows that there are precisely two codimension one components, thereby proving a conjecture of Landsberg\cite{landsberg2015geometric}.

In Popov \cite{popov2009two}
the author gives a constructive algorithm to determine if $x$ is in the orbit closure of $y$. He does this in both, the affine as well as the projective setting. This establishes the decidability of the orbit closure membership problem in the sense of computability theory.

\subsection{Our results and the organization of the paper.}
In the spirit of Luna \cite{luna1973slices}, in  this paper we develop a model for a neighbourhood of $x$, when $x$ has a non trivial stabilizer. 
We aim to study the set of such points $y$ for which $x$ may appear as a limit point. Moreover, in the sense of \cite{mulmuley2001geometric}, we consider points $x$ and $y$ which have distinctive stabilizers. While the stabilizer $K \subseteq G$ of $y$ is reductive, the stabilizer $H$ of $x$ is typically not reductive.

To describe our approach it will be useful to  think of $G$ as a complex Lie group with Lie algebra ${\cal G}$, see \cite{lee2001structure}. Let $H=\{g | \rho(g) \cdot x = x\}$ be the stabilizer of $x$ in $G$ and ${\cal H}$ its Lie algebra. Our main technical construction, Theorem~\ref{thm:main}, is an explicit Lie algebra homomorphism from ${\cal G}$ to vector fields at points in a neighbourhood of $x$. As a consequence we show that the Lie algebras of the stabilizers of points in this neighbourhood are parameterized by subspaces of ${\cal H}$. When $x$ has a reductive stabilizer, our results imply that the Lie algebras of stabilizers of points in the neighbourhood are Lie subalgebras of ${\cal H}$. Such a conclusion would also follow from Luna's theorem for polystable $x$.

To describe the ${\cal G}$-action, we start with $N$, a linear subspace of $V$ that is complementary to $TO_x$, the tangent space to the orbit $O_x$. 
Let $M$ be any locally closed submanifold of $G$ passing through $e$ that is complementary to $H$. Then $M \times N$ is a locally closed submanifold of $G \times V$ containing $(e,0)$. Next, the map $ G \times V \to G\times^H V$ is ${\cal G}$-equivariant, and the locally closed embedding of $M \times N \to G\times V$ gives us a quotient ${\cal G}$ action on $M\times N$. 
The map from $G \times^H V$ to $V$ sending $[g,v] \rightarrow g(x+v)$
when restricted to $M \times N$ gives us a ${\cal G}$-equivariant diffeomorphim of a neighbourhood $W$ of $(e,0) \in M\times N$ to a neighbourhood of $x$ in $V$.  An important step
in this paper is an explicit description
of this action of ${\cal G}$ on $W$ in terms of the projections $W \to M$ and
$W \to N$ coming from the product $M\times N$.
 An essential component of this construction, described in Section~\ref{thetamap}, is a map $\theta : V \times N \rightarrow V$ which captures the action of the nilpotent part of ${\cal H}$. In Remark~\ref{rem:SFTtheta} we show that this map is intimately related to the second fundamental form at $x$.  
\eat{$M \times N$ is a locally closed submanifold of $G \times V$. The natural map from $G \times V$ to $G \times^H V$ embeds $M \times N$ as a locally closed submanifold $\tilde{W}$ of 
$G \times^H V$ containing $[e,0]$. The map from $M \times N$ to $V$ given by $(m,n) \to m(x+n)$ factors through $\tilde{W}$ to give an open embedding of a neighbourhood $W$ of $[e,0]$ in $\tilde{W}$ to $V$. There is an action of $G$ on $V$. Therefore the open embedding of $W$ into $V$
induces a local action of $G$ on $W$,
and hence an action of ${\cal G}$ on $W$. An important contribution 
in this paper is an explicit description
of this action of ${\cal G}$ on $W$ in terms of the projections $W \to M$ and
$W \to N$ coming from the product $M\times N$. To do this we define a map $\theta : V \times N \rightarrow V$ which captures the action of the nilpotent part of ${\cal H}$.} 
We call this data, of a neighbourhood of $x$ with an explicitly defined ${\cal G}$-action, a local model at $x$. We develop the local model in Section~\ref{sec:prelim}.

We illustrate the utility of the local model in two applications. 
In Section~\ref{sec:forms} we apply this to the situation $V=Sym^d(X)$, the space of degree $d$ forms in the variables $X$. Let $f$ and $g$ be two forms with stabilizer (Lie algebras) ${\cal K}$ and ${\cal H}$ respectively. 
We consider a form $f$ (as the ``$y$" above) and elements in its orbit, of form $f(t)=A(t)\cdot f$, where $A(t)$ is an invertible family in $GL(X)$. We express $f(t)$ as $f(t)=t^a g +t^b f_b + \mbox{ higher terms}$, with $a<b$. We call $g\in Sym^d (X)$, the leading term as the limit point (and ``$x$" above) and $f_b $ as the tangent of approach. We construct a basis $\{ \mf{k}_i (t) = \mf{h}_i (t) +\mf{s}_i (t)\}_{i=1}^{dim({\cal K})}$ of $f(t)$ for generic $t$, with the coefficients of $\mf{h}_i (t)$ belonging to ${\cal H}$. In Theorem~\ref{theorem:mainform} we relate the stabilizer ${\cal H}$ of $g$ with ${\cal K}$, the stabilizer of $f$ as follows. The above data gives us an action of ${\cal H}$ on $\overline{N} = V/TO_g $. We show that there is a Lie subalegbra ${\cal K}_0 \subset {\cal H}$ of the same dimension as ${\cal K}$ such that ${\cal K}_0$ stabilizes $\overline{f_b} \in \overline{N}$. Furthermore ${\cal K}_0$ is the Lie subalgebra spanned by $\{ {\mf h}_i(0)\}$.  

It turns out that 
there is an $\epsilon$-extension Lie algebra ${\cal K}(\epsilon)$ of ${\cal K}_0$ whose structure constants are incrementally closer to ${\cal K}$. The extendability of ${\cal K}_0$ to ${\cal K}(\epsilon)$ translates into certain cohomology conditions. We discuss this connection to Lie algebra cohomology in Section~\ref{sec:cohomo}. 

In Section~\ref{lambdat} we continue our discussion of forms and consider the case when $g$ is the limit of a 2 block 1-parameter subgroup acting on $f$ (see Section~\label{sec:1pssimple} for definition). The grading we have on the underlying space by assuming so allows us to prove stronger results. In Proposition~\ref{prop:stabgeneric} we show that either there is an element which stabilizes $f, g$ and the tangent of approach (i.e, a triple stabilizer), or the limiting Lie algebra ${\cal K}_0$ is nilpotent. In Proposition~\ref{prop:generalstab} we extend this to the case when $g$ is the limit of a general $1$-ps acting on $f$. 

In Section~\ref{sec:codim1} we illustrate the results of Section~\ref{sec:1pssimple} when $f$ has a reductive stabilizer and the orbit of $g$ is of codimension $1$ in the projective closure of $f$. In particular we give details when $f$ is the determinant of a $3 \times 3$ matrix.

In Section~\ref{sec:conj} we apply the local model to $Mat_n (\C)$, the vector space of $n\times n$-matrices under conjugation by $GL_n$. We show in Theorem~\ref{thm:conjugation:main} that for a general diagonalizable matrix $X$, the signature of nilpotent matrices in the projective closure of its orbit are determined by the multiplicity data of the spectrum of $X$. 

In Section~\ref{sec:diffgeom}
we begin a study of differential geometry at the point $y$ in the local model. We believe that differential geometric techniques may be useful to further relate properties of the point $y$, the tangent vector of exit from $y$ and the ultimate limit $x$. This is illustrated by an example of the cyclic shift matrix in Section~\ref{sec:cyclic}, which requires us to compute the Reimannian curvature tensor. We also make connections to Kempf's theorem in Proposition ~\ref{prop:diff-kempf}. 

\noindent {\bf Acknowledgements}:
The second author would like to thank Prof Ketan Mulmuley, Computer Science Department, University of Chicago, for many discussions in the period 1999-2002 during which the $\theta$ map was discovered. The third author was supported by a MATRICS grant from the Department of Science and Technology, India and by a grant from the Infosys foundation. 
\section{The local model}
\label{sec:prelim}
Let $G$ be a complex algebraic, connected, reductive group, with Lie algebra ${\cal G}$, acting on a complex vector space $V$, by a rational representation $\rho:G \rightarrow 
GL(V)$. Our interest is in a local model for the action of ${\cal G}$, the Lie algebra of $G$, on a neighbourhood of a point $x \in V$, with stabilizer $H$, and Lie algbra ${\cal H}$. Note that $\rho :G\rightarrow GL(V)$ induces a Lie algebra map $\rho :{\cal G} \rightarrow End(V)$. 

An elementary  description of a neighborhood of $x$ is as follows. 
\begin{itemize}
\item The orbit $O(x)=G\cdot x$ is a smooth manifold of dimension $dim(G)-dim(H)$.
\item Let ${\cal S}\subseteq {\cal G}$ be a complement to ${\cal H}$ within ${\cal G}$, i.e., ${\cal G}={\cal S}\oplus {\cal H}$. Then 
\[ {\cal S} \cdot x= \{ \rho (\mf{s})\cdot x | \mf{s} \in {\cal S}  \} \] is the tangent space $T_x O(x)$ to the orbit $O(x)$.
\item Let $N\subseteq T_x V$ be a complement to $T_x O \subseteq T_x V$, such that $T_x V = T_x O \oplus N = {\cal S}\cdot x \oplus N$. In any chosen neighborhood $U \subseteq V$ of $x$, we identify $T_u V$ with $V$, for all $u\in U$. We may choose a neighborhood $U$ of $x$ such that for all $u\in U$, $T_u V = {\cal S}\cdot u \oplus N$.
\end{itemize}

We would like to describe the action of any $\mf{g}\in {\cal G}$ in this neighborhood $U$ in terms of the above division of $T_u V = {\cal S}\cdot u \oplus N$. 
In particular, we will choose $u=x+n$, i.e.,  the section $x+N$, as the domain of interest. The specification of $\rho (\mf{g})$ on such a section may be suitably extended to sections at other points $x'\in O(x)$ in the vicinity of $x$ and therefore to a neighborhood of $x$. 
We show that the action of $\mf{g}$ on this section is better understood by (i) using the decomposition ${\cal G}={\cal S}\oplus {\cal H}$ to express $\mf{g}\in {\cal G}$ as $\mf{g}=\mf{s}+\mf{h}$ with $\mf{s}\in {\cal S}$ and $\mf{h} \in {\cal H}$, and (ii) obtaining the descriptions of $\mf{s}$ and $\mf{h}$ separately, on the section $x+N$. 
This description is the local model and it relies on the following basic definition and theorem on the formation of quotients.

\subsection{The quotient model}
\begin{defn}
Let ${\cal X}$ be a smooth manifold and let $Vec({\cal X})$ be the Lie algebra of smooth vector fields on ${\cal X}$. Let ${\cal G}$ be a Lie algebra. We say that ${\cal X}$ is a ${\cal G}$-manifold via $\gamma$ if there is a Lie algebra homomorphism $\gamma : {\cal G}\rightarrow Vec({\cal X})$. We denote this by the data $({\cal G},\gamma, {\cal X})$.

Next, let $({\cal G},\gamma ,{\cal X})$ and $({\cal G},\rho ,{\cal Y})$ be two ${\cal G}$-manifolds. A map $\mu : {\cal X}\rightarrow {\cal Y}$ is called ${\cal G}$-equivariant iff for all $x\in {\cal X}$ and $\mf{g}\in {\cal G}$, we have $\mu^*_x (\gamma (\mf{g})(x))=\rho (\mf{g}) (\mu (x))$, 
where $\mu^*_x :T_x {\cal X} \rightarrow T_{\mu(x)} {\cal Y}$ is the tangent map corresponding to the map $\mu $. 

We define ${\cal D}_{\mu} (x) \subseteq T_x {\cal X}$ as those elements $\Delta x \in T_x {\cal X}$ such that $\mu^*_x (\Delta x)=0$. Thus ${\cal D}_{\mu}$ records the kernel of the map $\mu^*$ at every $x \in {\cal X}$.  
\end{defn}
\begin{theorem} \label{theorem:Gequiv}
Let $({\cal G},\gamma ,{\cal X})$ and $({\cal G},\rho ,{\cal Y})$ be two ${\cal G}$-manifolds of dimensions $m$ and $n$, respectively.  
Let $\mu :{\cal X}\rightarrow {\cal Y}$
be ${\cal G}$-equivariant. Suppose that at the point $x\in {\cal X}$, we have $\mu^*_x (T_x {\cal X}) =T_{\mu(x)} {\cal Y}$. 
Moreover, suppose that ${\cal M}$ is a locally closed submanifold of ${\cal X}$ of dimension $m-n$ such that $x\in {\cal M}$ and the tangent space $T_x (M)$ is transverse to ${\cal D}_{\mu}$ at $x$. Then:
\begin{enumerate}
    \item There is an open submanifold ${\cal M}'\subseteq {\cal M}$ containing $x$ such that for all $m\in {\cal M}'$, the tangent space $T_x {\cal M}'$ is complementary to ${\cal D}_{\mu} (m)$. 
    \item There is a corresponding ${\cal G}$-action $\gamma_{{\cal M}'}$ on ${\cal M}'$ such that the map $\mu : {\cal M}' \rightarrow {\cal Y}$ is a ${\cal G}$-equivariant diffeomorphism in a neighborhood of $x$.
    \item Moreover, for $m\in {\cal M}'$, since we have $T_m {\cal X} = T_m {\cal M}' \oplus {\cal D}_{\mu} (m)$,  let $\pi_{\cal M'}: T_m {\cal X} \rightarrow T_m {\cal M}'$ be the corresponding projection. Then we have:
    \[ \gamma_{{\cal M}'} (\mf{g})(m) = \pi_{{\cal M}'} (\gamma (\mf{g}) (m)) \mbox{   for all $\mf{g}\in {\cal G}$ and $m\in {\cal M}'$}\] 
\end{enumerate}
\end{theorem}

For the definition of a ${\cal G}$-manifold see ~\cite{alekseevsky1995differential}, ~\cite{michor2008topics}. The proof of the theorem is straightforward. It allows the construction of a local quotient model ${\cal M}'$ within ${\cal X}$ of the space ${\cal Y}$. 

\begin{defn}
Let $G$ act on $V$ via $\rho$ as before. The $G$-space  $G\times V$ is the collection of tuples $(g,v)$ with $g\in G$ and $v\in V$.  The action $\gamma $ of $G$ on $G\times V$ is given by $\gamma (g')(g,v)=(g'g,v)$. Let $x\in V$ be a point with stabilizer $H\subseteq G$, as above. The map $\mu :G\times V \rightarrow V$ is given by $\mu (g,v)=\rho(g)(x+v)$ (abbreviated as $g \cdot (x+v)$). \end{defn}

\begin{lemma} \label{lemma:main}
The following assertions hold:
\begin{enumerate} 
\item[(i)] The map $\mu$ is $G$-equivariant. 
\item[(ii)] Let $v\in V$, then the tangent space $T_{(e,v)} (G\times V)$ is $(\mf{g},\Delta v)$, where $\mf{g}\in {\cal G}$ and $\Delta v \in T_v V \cong V$.  
\item[(iii)] The tangent map is given as follows: $\mu^*_{(e,v)} (\mf{g},\Delta v)= \mf{g}\cdot (x+v)+\Delta v$. Whence, the distribution ${\cal D}_{\mu } (e,v)=\{ (\mf{g}, \Delta v) |\mf{g}\cdot (x+v)+\Delta v =0 \}$. 
\item[(iv)] $\mu^*_{(e,0)} = \mf{g} \cdot x +\Delta v$. Thus $\mu^* $ is surjective at $(e,0)$. 
\end{enumerate}
\end{lemma}

The proof is straightforward.

\begin{prop}
\label{prop:main}
Let $N$ be a vector space complement to $TO(x)_x$ within $V$. Let ${\cal S}$ be a complement to ${\cal H}$ within ${\cal G}$. 
Let $M$ be a locally closed submanifold of $G$ of dimension $dim(G)-dim(H)$ such that $e\in M$ and $T_e M ={\cal S}$. Then there is an open submanifold $M' \subseteq M$ containing $e$, such that ${\cal G}$ acts on $M'\times N \subseteq G\times V$ via $\gamma '$. 
This action is derived from the action $\gamma $ of ${\cal G}$ on $G\times V$ modulo the distribution ${\cal D}_{\mu}$. 
Moreover, for this neighborhood of $M'\times N$ containing $(e,0)$, we have: 
\begin{itemize} 
\item[(i)] The map $\mu :M' \times N \rightarrow V$ is a ${\cal G}$-equivariant diffeomorphism from the space $({\cal G},\gamma',M'\times N)$ to $({\cal G},\rho ,V)$. It maps a neighborhood of $[e,0]\in M'\times N$ to a neighborhood of $x\in V$. 
\item[(ii)] For any element $\mf{g}=\mf{s}+\mf{h}$, with $\mf{s}\in {\cal S}$ and $\mf{h}\in {\cal H}$, the action of $\gamma '(\mf{g})$ at the point $(e,n)\in M'\times N$ is given by
\[ \gamma '(\mf{g}) (e,n)=(\mf{s}+\mf{s}', n') \in T_{(e,n)} (M'\times N)\]
where $\mf{s}'\in {\cal S}$ and $n'\in N$ are unique and satisfy the equation 
\[ \mf{s}'\cdot (x+n) +n'=\mf{h}\cdot n.\]
\end{itemize}
\end{prop}

\noindent 
{\bf Proof}: Let us verify that the conditions required by Theorem \ref{theorem:Gequiv} hold for the map $\mu$. From lemma \ref{lemma:main}, (i) and (iv), $\mu$ is $G$-equivariant and $\mu^*$ is surjective at $(e,0)$. 
By the choice of $M$ and $N$, we see that $\mu^* (e,0) (M\times N)={\cal S}\cdot x+N=V$. Thus $T_{(e,0)}M\times N$ is transverse to ${\cal D}_{\mu}$. 
Thus Theorem \ref{theorem:Gequiv} applies at the point $(e,0)$ and a suitable action $\gamma'$ on $M'\times N$ is available for which $\mu$ is an equivariant diffeomorphism. Thus (i) is proved. 
We also have  that $M'\times N$ is transverse to the distribution $D_{\mu}$. 
To prove (ii), we must make an explicit computation of $\gamma'$ at the points  $(e,n) \in M'\times N \subseteq G\times V$ modulo the distribution $D_{\mu }$. For this, we must evaluate $\pi_{M'} : T_{(e,n)} (G\times V) \rightarrow T_{(e,n)} (M\times N)$. We do this in two steps:
\[ T_{(e,n)} G\times V \rightarrow T_{(e,n)} M\times V \rightarrow T_{(e,n)} M\times N \]

From Lemma \ref{lemma:main}(ii), for a typical element $\mf{g}=\mf{s}+\mf{h}$, we have $\gamma (\mf{g})(e,n)=
(\mf{s}+\mf{h}, 0)\in T_{(e,n)}(G\times V)$. We observe that for any $\mf{h} \in {\cal H}$, we have 
\[ \mu^*_{(e,n)} (\mf{h},0)= \mf{h}\cdot (x+n)=\mf{h}\cdot n =\mu^*_{(e,n)} (0,\mf{h} \cdot n)  \] 
 Thus, we have:
\[ \gamma_M (\mf{s}+\mf{h})(e,n)=(\mf{s}+\mf{h},0)=(\mf{s}, \mf{h} \cdot n) \in T_{(e,n)} M\times V \mbox{   (modulo ${\cal D}_{\mu}$)}\]
Thus, we have reduced the first component from ${\cal G}$ to ${\cal S}=T_e M$. 

Using Lemma \ref{lemma:main}(iii) it follows that 
\[ \mu^*_{(e,n)}(\mf{s}, \mf{h} \cdot n) = \mf{s}\cdot (x+n)+ \mf{h}\cdot n \in T_{(x+n)}V \]

Consider the vector $\mf{h}\cdot n \in T_{(x+n)}V$. Since $\mu $ is a local diffeomorphism between $ (e,0)\in M\times N$ and $x\in V$, we may assume that $\mu$ is regular at $x+n$. 
Whence, there must be a unique $\mf{s'}\in {\cal S}$ and $n' \in N$ such that for $(\mf{s'},n') \in T_{(e,n)}M\times N$, we have $\mu^*_{(e,n)} (\mf{s}',n')=\mf{h}\cdot n$, or in other words:
\[ \mf{s}'\cdot (x+n) +n'=\mf{h}\cdot n. \]
Hence
\[ \mu_{* (e,n)} (\mf{g}\cdot [e,n])= \mu^*_{(e,n)} (\mf{s}, \mf{h}\cdot n)=\mf{s}\cdot (x+n)+ \mf{h}\cdot n=(\mf{s}+\mf{s}')\cdot (x+n)+n'=\mu^*_{ (e,n)} (\mf{s}+\mf{s'},n')\] 
Since $(\mf{s}+\mf{s'},n')\in T(M\times N)_{(e,n)}$, it must be the projection of $(\mf{s}, \mf{h} \cdot n)$ to $T_{(e,n)} (M\times N)$, and thus from Theorem~\ref{theorem:Gequiv}(3) it describes the action of $\mf{g}$ on $M\times N$ at the point $(e,n) \in M\times N$ . This proves (ii). $\Box $

\subsection{The $\theta$ map.}
\label{thetamap}It is now clear that the computation modulo the distribution $
D_{\mu }$ at $(e,n)$ essentially revolves around the
solution to the following equation in $V$, where ${\Delta v}\in V\cong T_{(x+n)}V$ is a given
and ${\mf s}\in {\cal S}$, ${\mf n} \in N$ are the unknowns:
\[ {\mf s} \cdot (x+n) +{\mf n} =\Delta v. \]

Since $V ={\cal S}\cdot (x+n)\oplus N$, we may define projection $\lambda_{\cal S} (n):V \rightarrow {\cal S}$ and $\lambda_N (n):V \rightarrow N$ so that the above equation may be written as:
\[ \Delta v = (\lambda_{\cal S}(n)(\Delta v) )(x+n)+(\lambda_N (n) (\Delta v)) \]
We now solve the above equation, with $n$ as a parameter. We note that for $n=0$, the equation is easily solved by our choice of $N$ and the isomorphism $V=T_x V =T_x O(x) \oplus N$ with $T_x O(x)={\cal S}\cdot x$. We may write the projections $\lambda_{\cal S} (0): V\rightarrow {\cal S}$ and $\lambda_N (0): V\rightarrow N$, or simply $\lambda_{\cal S}, \lambda_N$, such that
\[ \lambda_{\cal S} (\Delta v) \cdot x + \lambda_N (\Delta v) =\Delta v.\]
By adding an ``error term'' $\lambda_{\cal S}(\Delta v) \cdot n$ on both sides, we get
\[ \lambda_{\cal S} (\Delta v)\cdot (x+n)+\lambda_N (\Delta
v)=\Delta v +\lambda_{\cal S}(\Delta v) \cdot n. \]

For a fixed $n \in N$, we may now define $\theta(n) :V\rightarrow V$,
as follows: for $\Delta v\in V$, 
we define $\theta(n) (\Delta v)=\lambda_{\cal S} (\Delta v)\cdot n$. 
We see that $\theta(n) :V \rightarrow V$ is a linear operator $V$ which depends linearly on $n$. 
Thus, the above equation may be rewritten as:
\[ \lambda_{\cal S} (\Delta v)\cdot (x+n)+\lambda_N (\Delta
v)=\Delta v +\theta (n)(\Delta v)\]
Substituting $\theta^i (\Delta v)$ for $\Delta v$ in the above equaltion gives us the $i$-th approximation term as:
\[ \lambda_{\cal S} ( \theta^i(n) (\Delta v))\cdot (x+n) + 
 \lambda_N ( \theta^i(n) (\Delta v)) =\theta^i(n) (\Delta
 v)+\theta^{i+1} (\Delta v). \]

 Since
$\theta(n) $ depends linearly on $n$, we may assume that $n$ is small enough and 
hence all the eigen-values of $\theta (n)$ may be assumed to be
less than $1$ in magnitude. Thus the system, 
\[ {\mf s}\cdot (x+n) +{\mf n}=\Delta v, \]
has the solution:
\[ [\lambda_{\cal S} \circ (1-\theta(n) +\theta^2(n) -\cdots )](\Delta v)\cdot
(x+n)+ 
[\lambda_N \circ (1-\theta(n) +\theta^2(n) -\cdots )](\Delta v) =\Delta v
\label{main:eqn}.\]
Thus, we have proved the following basic lemma:
\begin{lemma} \label{lemma:thetasolve}
For any element $\Delta v \in T_{(x+n)}V$, there are unique $\mf{s} \in{\cal S}$ and $n'\in N$ such that:
\[ \mf{s} \cdot (x+n)+n' =\Delta v \]
Moreover $\mf{s}=\lambda_{\cal S} \circ ((1+\theta (n))^{-1} (\Delta v))$ and $n'=\lambda_N \circ (1+\theta (n))^{-1} (\Delta v)$. Thus, if ${\cal I}(n):V \rightarrow V$ is the map $(1+\theta (n))^{-1}$, then $\lambda_{\cal S} (n)= \lambda_{\cal S} \circ {\cal I}(n)$ and $\lambda_N (n) =\lambda_N \circ {\cal I}(n)$. 
\end{lemma}

We will use the above calculations and Proposition \ref{prop:main} to prove our main theorem.

Recall that $H$ stabilizes the point $x$. We may choose a Levi decomposition $H=RQ$, where $Q$ is the unipotent radical of $H$ and $R$ is a reductive complement. Let ${\cal R}$ denote the Lie algebra of $R$ and ${\cal Q}$ denote the Lie algebra of $Q$. 

We may choose $N$ in the above proposition to be an $R$-invariant complement to the space $O_x $ inside $V$. We may also choose ${\cal S}$ to be an $R$-invariant complement to ${\cal H}$ with ${\cal G}$.  Thus ${\cal S}$ now acquires an $R$-module
structure. 

 \begin{theorem}
\label{thm:main}
Let $\rho$ be a representation of $G$ on a vector space $V$. Let $x \in V$ be a point with stabilizer $H$. Let $H=RQ$ be a Levi decomposition and ${\cal H}, {\cal Q}$ and ${\cal R}$ be the Lie algebras of $H,Q$ and $R$, respectively. Then there is a submanifold  $M \times N \subset G \times_H V$ in the neighbourhood of $(e,0)$ with an action of ${\cal G}$. The map $\mu :M \times N \rightarrow V$ given by $\mu(g,n) = g\cdot (x+n)$ is a ${\cal G}$-equivariant diffeomorphism which maps a neighborhood of $(e,0)\in M\times N$ to a neighborhood of $x\in V$. We then have:
\begin{enumerate}
\item Let $\mf{g}$ be a typical element of ${\cal G}$. Let $\mf{g}=\mf{s} + \mf{r} + \mf{q}$, with $\mf{s}\in {\cal S}, \mf{r}\in {\cal R}$ and $\mf{q}\in {\cal Q}$. Then the action $\gamma'$ of $\mf{g}$ on $M\times N$ is given as: 
\[ \mf{g}\circ (e,n)=(\mf{s}+ \lambda_{\cal S} \circ (1+\theta(n))^{-1}(\mf{q} \cdot n), \: \mf{r}\cdot n + \lambda_N \circ (1+\theta(n))^{-1}(\mf{q} \cdot n ) ) \]

\item If ${\cal G}_n \subseteq {\cal G}$ is the stabilizer of $x+n\in V$, then there is a subspace ${\cal H}_n \subseteq {\cal H}$ of the same dimension as ${\cal G}_n$ such that for every $\mf{g}\in {\cal G}_n$, there is a unique $\mf{h}\in {\cal H}_n$ with $\mf{h}=\mf{q}+\mf{r}$ such that: 
\[ \mf{r}\cdot n + \lambda_N \circ (1+\theta(n))^{-1}(\mf{q} \cdot n ) = 0\]
Moreover,  $\mf{g}=\mf{h}+\mf{s}$, where $\mf{s}$ is as given below:
\[ \mf{s} = -\lambda_{\cal S} \circ (1+\theta(n))^{-1}(\mf{q} \cdot n)\]
The element $\mf{g}$ is called the ${\cal S}$-completion of $\mf{h}$ at the point $x+n$. 
\end{enumerate}
\end{theorem}

\noindent
{\bf Proof}:  Let us express $\mf{g}=\mf{s}+\mf{h}$ and apply Prop. \ref{prop:main}, part (ii). We then have that $\gamma'(\mf{g})(e,n)=(\mf{s}+\mf{s'}, n')$, where $\mf{s}'$ and $n'$ solve the equation:
\[ \mf{h} \cdot n =\mf{s}'(x+n)+n' \]
By lemma \ref{lemma:thetasolve}, these are precisely $\lambda_{\cal S} (1+\theta (n)^{-1} (\mf{h}\cdot n)$ and $n'=\lambda_N (1+\theta (n))^{-1} (\mf{h}\cdot n)$. 
Now, breaking $\mf{h}$ as $\mf{h}=\mf{q}+\mf{r}$, and recalling that $\mf{r} \cdot N \subseteq N$, we see $\lambda_N (\mf{r}\cdot n)=\mf{r}\cdot n$, and $\lambda_{\cal S} (\mf{r}\cdot n)=0$. 
Thus $\theta  (n) (\mf{r} \cdot n)=0$ and thus $\lambda_N (1+\theta (n))^{-1} (\mf{r} \cdot n +\mf{q} \cdot n)=\lambda_N (\mf{r} \cdot n) +\lambda_N (1+\theta (n))^{-1} (\mf{q} \cdot n)=\mf{r} \cdot n +\lambda_N (1+\theta (n))^{-1} (\mf{q} \cdot n)$. 
Similarly, $\lambda_{\cal S} (1+\theta (n))^{-1} (\mf{r} \cdot n +\mf{q} \cdot n)=\lambda_{\cal S} (\mf{r} \cdot n) +\lambda_{\cal S} (1+\theta (n))^{-1} (\mf{q} \cdot n)=\lambda_{\cal S} (1+\theta (n))^{-1} (\mf{q} \cdot n)$. 

This proves (1). 

The hypothesis that $\mf{g} \in {\cal G}_n$ is equivalent to the claim that $\gamma' (\mf{g}) (e,n)=0$. This requires both the components of the action of $\mf{g} \cdot (e,n)$ to vanish. This implies that $\mf{r}\cdot n + \lambda_N \circ (1+\theta(n))^{-1}(\mf{q} \cdot n ) = 0$. This is the first assertion. 

Now, given that $\mf{h}=\mf{q}+\mf{r}$ satisfies the above, $\mf{s}$ is uniquely determined by the requirement that the first component vanishes, i.e.,
$\mf{s} + \lambda_{\cal S} \circ (1+\theta(n))^{-1}(\mf{q} \cdot n)$ must vanish. This determines $\mf{s}$ as required.

Finally, if there were $\mf{g}'=\mf{h}+\mf{s}$ and $\mf{g}'=\mf{h}+\mf{s}'$ both elements of ${\cal G}_n$, then we would have $\mf{s}-\mf{s'}\in {\cal G}_n$. Whence, $(\mf{s}-\mf{s}')\cdot (e,n)=(\mf{s}-\mf{s}', 0) \neq 0$ which violates the requirement that the first component vanishes. 
$\Box $

\begin{remark}
Consider the case when the stabilizer $H$ is reductive, whence $R=H$ and $Q=0$. We may select an $N$ which is $H$-invariant as above. Since ${\cal Q}=0$, we see that $\mf{q}\cdot n=0$, and thus, the ${\cal S}$-completion of any element $\mf{h}\in H$ is $\mf{h}$ itself. Thus the stabilizer ${\cal G}_n$ of $x+n$ must be a subalgebra of ${\cal H}$. This is the Lie algebraic version of Luna's slice theorem. 
\end{remark}

\noindent
We now prove certain equivariant properties of maps closely related to $\theta $. We continue with the Levi decomposition of $H=RQ$ and the choice of $N$ and ${\cal S}$. 

The Lie algebra action defined as a map $\rho : {\cal G}\rightarrow End(V)$ may well be defined as $\hat{\rho}: {\cal G}\otimes V \rightarrow V$. 
Note that ${\cal G}$ is an $R$-module under the adjoint action, and thus so is ${\cal G}\otimes V$. Then we have the following lemma:
\begin{lemma}
The map $\hat{\rho}$ is $R$-equivariant. Next, let $T_x O$ be the tangent space of the orbit of $x$ under the action of $G$ and ${\cal S}$ be an $R$-invariant complement to ${\cal H}\subseteq {\cal G}$. Then the map $\alpha :{\cal S}\rightarrow T_x O$ given by $\alpha (\mf{s})=\rho(\mf{s}) \cdot x$ is an isomorphism of $R$-modules. 
\end{lemma}

\noindent 
{\bf Proof}. For the first assertion, note that for $r\in R, v\in V$ and $\mf{g}\in {\cal G}$, we have:
\[ r \cdot \hat{\rho} (\mf{g}\otimes v) = 
 r\cdot \rho (\mf{g}) \cdot r^{-1} \otimes r\cdot v= r\cdot (\rho (\mf{g})\cdot v)\]
Thus $\hat{\rho}$ is indeed $R$-equivariant. Next, the map $\alpha :{\cal S}\rightarrow  T_x O$ is merely the restriction to $\hat{\rho}$ to the tensor product of two $R$-submodules, viz., ${\cal S}\subset {\cal G}$ and $\C \cdot x \subseteq V$. This proves the second assertion. $\Box$

Recall that $\theta (n):V\rightarrow V$ is defined as the composite:
\[ \theta : V\stackrel{\lambda_{\cal S}}\rightarrow {\cal S} \stackrel{\cdot n}\rightarrow V \]
We define two associated maps:
$\Theta :V \otimes N \rightarrow V$ and $\Phi : {\cal S}\times N \rightarrow {\cal S}$ as:
\[ \Theta (v\otimes n)=\theta(n)(v)= \lambda_{\cal S} (v) \cdot n \]
and:
\[ \Phi (\mf{s} \otimes n )=\lambda_{\cal S} (\mf{s}\cdot n) \]
Then, we have the following:

\begin{prop}
The maps $\Theta $ and $\Phi$ are $R$-equivariant.
\end{prop}

{\bf Proof}: $\Theta$ is merely the composite:
\[ V\otimes N \stackrel{\lambda_{\cal S}\otimes{\tt id}}{\longrightarrow} {\cal S}\otimes N \stackrel{\hat{\rho}}{\longrightarrow} V \] 
As for $\Phi $, it is:
\[ {\cal S}\otimes N \stackrel{\hat{\rho}}{\longrightarrow} V \stackrel{\lambda_{\cal S}}{\longrightarrow} {\cal S} \]
This proves the proposition. $\Box $

\begin{remark}
\label{rem:SFTtheta}
The map $\Phi$ above is the algebraic analogue of the Second Fundamental Tensor at the point $x$ for the orbit $O(x)$ regarded as a submanifold of $V$. 
Indeed, suppose that $G$ is the orthogonal group and $V$ is equipped with a $G$-invariant inner product. Let $N$ be chosen orthogonal to $T_x O$, the tangent space of the orbit $O(x)$ at $x$. Let $v_1 ,\ldots ,v_K$ be an orthonormal basis for $N$. Let $\mf{s}_1 ,\ldots ,\mf{s}_L$ be chosen such that the elements $\{ v_{K+i} =(\mf{s}_i \cdot x)| i=1,\ldots ,L\}$ is an orthonormal basis for $T_x O$. 
The data for $\Phi $ is the tuple $\{(\alpha_{ij}^r)|i,j=1,\ldots,L, r=1,\ldots,K\}$, where we have:
\[ \lambda_{\cal S} (\mf{s}_i \cdot v_r )=\sum_j \alpha_{ij}^r (\mf{s}_j \cdot x) \]
Now, one must recall that $\rho (\mf{s}_i)$ may be given as a linear vector field as follows:
\[ \rho (\mf{s}_i)=  \sum_{j=1}^{K+L} \sum_{k=1}^{K+L} S_i (j,k)v_k \frac{\partial}{\partial v_j } \] where $S_i \in \C^{ (K+L)\times (K+L)}$. We abbreviate this as $\rho(\mf{s}_i)=\nabla S_i \overline{v}$, where $\nabla $ is the row vector $[(\partial /\partial v_1 ),\ldots ,(\partial /\partial v_{K+L})]$ and  $\overline{v}$ is the column vector $[v_1 ,\ldots ,v_{K+L}]^T$.

Note that $\rho$ preserves inner products implies that $S_i^T =-S_i$. The data $\alpha_{ij}^r$ may then be obtained as:
\[ \alpha_{ij}^r = (S_j x)^T (S_i v_r )= -x^T S_j S_i v_r \] 
We also have the Second Fundamental Form of the orbit $O(x)$ as a map $\Pi : T_x O \times T_x O \rightarrow N$ defined as follows. For $v_i ,v_j \in T_x O$, let $X_i,X_j$ be vector fields on $O(x)$ such $X_i (x)=v_i$ and $X_j (x)=v_j$. Then $\Pi (v_i ,v_j )=\lambda_N (D_{X_i} (X_j))=\sum_{r=1}^K \beta_{ij}^r v_r$.  
Now, we have $\rho(\mf{s}_i)$ available to us as the $X_i$ above. Then, it is easily checked that $D_{X_i}(X_j)$ is given  $\nabla S_j S_i \overline{v}$. 
Whence the evaluation of $\Pi  (v_i ,v_j)$ at the point $x$ and in the direction $v_r$ is given by $\beta_{ij}^r=v_r^T S_j S_i x$.
However, this is equal to $v_r^T S_i S_j x +v_r^T [S_j ,S_i ]x$. 
Since there is a $\mf{g}=[\mf{s}_i ,\mf{s}_j]\in {\cal G}$ such that $\rho (\mf{g})=\nabla [S_j ,S_i] \overline{v}$, this evaluated at $x$ must be in the tangent space $T_x O$, and thus $v_r^T [S_j ,S_i ]x=0$. Whence 
\[ \beta_{ij}^r =v_r^T S_i S_j x=x^T S_j^T S_i^T v_r =x^T S_j S_i v_r =-\alpha_{ij}^r\]
This proves the claim. $\Box$ 
\end{remark}

We explore further uses of the Reimannian curvature tensor 
to study projective orbit closures in Section~\ref{sec:diff}.

\subsection{Applying the local model}

We close this section with a general application of the local model to projective limits of a particular type, and a simple concrete example. 
\begin{prop} \label{prop:family}
Let $V$ be a $G$-module, where $G\subseteq GL(X)$, for a suitable $X$ and let $p,q \in V$. Let $I$ be the identity matrix and $tI$ act on $p$ as $tI\cdot p=t^d p$.  Let $A(t)\subseteq G$ be an algebraic family of elements in $G$ such that:
\[ A(t)\cdot q =t^k \cdot p + \mbox{ higher degree terms } \]
Suppose that the stabilizer of $p$ is $H$. Then the local model at $p$ contains a scalar multiple of a $G$-conjugate of $q$ which is of the form $p+n$, i.e., appears on the normal slice $N$ at $p$. 
\end{prop}

\noindent 
{\bf Proof}. Consider the family $A'(t)=t^{-k/d} \cdot A(t)$, i.e., whose elements are scalar multiples of $A(t)$ and thus are elements of $G$. Applying this family to $q$, we see that:
\[ A'(t) \cdot  q= p + \mbox{ positive powers of $t$} \]
Thus, $\lim_{t\rightarrow 0} A'(t)q=p$, and for every $\epsilon >0$, there is a $q_0 =A'(t_0 ) \cdot q$, i.e., a scalar multiple of a $G$-conjugate of $q$ such $\| q_0 -p\|<\epsilon$. 
Thus, in the local model $\mu : M\times N \rightarrow  V$, we have  
a $(g,n)$ such that $g\cdot (p+n)=q_0$. Applying $g^{-1}$ on both sides, we see that $g^{-1} q_0 =p+n$. This proves the claim. $\Box $

Let us consider a small example to illustrate the computational ease of the local model. 

\begin{ex}
Consider $sl(2)$ given as the matrix below acting of $\C \cdot X$, where $X=\{ x,y\} $. Let $V=Sym^2 (X)$ with the ordered basis $B=\{ x^2 ,xy, y^2 \}$.  The action of the generic group element is also shown below:
\[ sl(2) =\left\{ g_{a,b,c}=\left[ \begin{array}{cr} 
a & b \\
c & -a \end{array} \right] \mbox{ with } a,b,c\in \C \right\}  \]
On $Sym^2 (X)$, it may be useful to interpret this as the operation:
\[ a(x\frac{\partial}{\partial x}  -y\frac{\partial}{\partial y}) +b y \frac{\partial}{\partial x} +c x \frac{\partial}{\partial y} \]
Whence:
\[
g_{a,b,c} \cdot \left[ \begin{array}{c}
x^2 \\
xy \\
y^2 \end{array} \right]
=\left[ \begin{array}{ccc}
2a & 2b & 0 \\
c & 0 & b \\
0 & 2c & -2a \end{array} \right] 
\left[ \begin{array}{c}
x^2 \\
xy \\
y^2 \end{array} \right]
\]
Let us consider three forms $p_0 =0, p_1 =x^2$ and $p_2 =x^2 +y^2$ and apply the local model to each.
\begin{enumerate}
    \item For $p_0 =0$, we see that ${\cal H}_0$, the stabilizer of $p_0$ is the complete group $sl(2)$. Thus ${\cal S}=0$, the orbit is $0$-dimensional, and the maps $\theta$ is $0$. The only choice for $N$ is $V$, and for any $n$, the stabilizer ${\cal G}_n \subseteq {\cal H}_0=sl(2)$.
    \item Consider next $p_2 =x^2 +y^2$. We see that the stabilizer ${\cal H}_2$ is the group $O\subset sl(2)$, given below. A complementary subspace is also shown. 
    \[ O =\left\{ \left[ \begin{array}{cr} 
0 & b \\
-b & 0 \end{array} \right] | b\in \C \right\} \: \: \: 
{\cal S} =\left\{ \left[ \begin{array}{cr} 
a & c \\
c & -a \end{array} \right] | a,c\in \C \right\}
\]
Since ${\cal S}$ is two dimensional, so is the orbit.  Applying ${\cal S}$ to $p_2$ gives $T_{p_2} O= \C \cdot (x^2 -y^2)+\C\cdot xy$. Thus a normal $n$ to the orbit is $x^2+y^2$, i.e., $p_2$ itself. Next, we compute $\theta (n) : V\rightarrow V$. For this, we must first compute $\lambda_S$, the ${\cal S}$-parametrization of the projection of $V\rightarrow T_{p_2} O$, and then apply this to $n$. But since $n=p_2$, this makes $\theta (n)$ a projection onto $T_{p_2} O$. On the other hand, since ${\cal H}_2$ is reductive, ${\cal Q}=0$, and the stabilizer of $p_2+n$, is $O_2$. 
\item Let us consider next $p_1 =x^2$. The stabilizer ${\cal H}_1$ and ${\cal S}$ are given below:
\[ {\cal H}_1 =\left\{ \left[ \begin{array}{cr} 
0 & 0 \\
c & 0 \end{array} \right] | c\in \C \right\} \: \: \: 
{\cal S} =\left\{ \left[ \begin{array}{cr} 
a & b \\
0 & -a \end{array} \right] | a,c\in \C \right\}
\]
 We see that ${\cal S}\cdot p_1 =2ax^2+2bxy$, and thus, we may choose $n=y^2$. Note that ${\cal R}=0$ and ${\cal Q}={\cal H}$. Let us choose $\mf{q}$ as the element $g_{0,0,1}$, and note that $\mf{q}\cdot n =2xy$. Let us compute $\theta (n)$ next. Since $\lambda_{\cal S} (y^2)=0$, we see:
\[
\lambda_{\cal S} \cdot \left[ \begin{array}{c}
x^2 \\
xy \\
y^2 \end{array} \right]
=\left[ \begin{array}{c}
g_{1/2,0,0} \\
g_{0,1/2,0} \\
0
\end{array} \right] \: \:\mbox{and thus  }
\theta (n) \cdot 
\left[ \begin{array}{c}
x^2 \\
xy \\
y^2 \end{array} \right]
=\left[ \begin{array}{ccc}
0 & 0 & -1 \\
0 & 0 & 0 \\
0 & 0 & 0 \end{array} \right] 
\left[ \begin{array}{c}
x^2 \\
xy \\
y^2 \end{array} \right]
\]
This gives us:
\[ (1+\theta (n))^{-1} \cdot 
\left[ \begin{array}{c}
x^2 \\
xy \\
y^2 \end{array} \right]
=\left[ \begin{array}{ccc}
1 & 0 & 1 \\
0 & 1 & 0 \\
0 & 0 & 1 \end{array} \right] 
\left[ \begin{array}{c}
x^2 \\
xy \\
y^2 \end{array} \right]
\]
Whence $\lambda_N (1+\theta (n))^{-1} (\mf{q}\cdot n)=0$ and thus ${\cal H}_n ={\cal H}_1$. The ${\cal S}$-completion is given by the condition that $\mf{s}+\lambda_{\cal S}(1+\theta (n))^{-1}(\mf{q} \cdot n)=0$. Since 
$\mf{q} \cdot n = 2xy$, we get $\lambda_{\cal S}(1+\theta (n))^{-1}(\mf{q} \cdot n) = g_{0,1,0}$. So $\mf{s}=-g_{0,-1,0}$, and thus, the ${\cal S}$-completion of the element $g_{0,0,1}$ is $g_{0,-1,1}$, the generator of $O$. Thus, the stabilizer of $p_2=p_1+y^2$ matches the one computed earlier. 
\end{enumerate}
\end{ex}

\section{Forms and their limits}
\label{sec:forms}
We will now specialize the above local model to the space of forms, and examine projective limits when they arise as leading terms under a 1-parameter family. This question is central to the framing of several questions in computational complexity theory. In Section~\ref{sec:A(t)} we start the investigation with a general 1-parameter family. Let $g$ be the limit of a form $f$ under the action of such a 1-parameter family $A(t)$. So $g$ is the leading term of $A(t) \cdot f$. We show a connection (Theorem~\ref{theorem:mainform}) between the stabilizer of $f$ at a generic value of the parameter $t$, the stabilizer of $g$ and the stabilizer of second leading term of $A(t) \cdot f$, what we call the tangent of approach.   
In Section ~\ref{lambdat} we consider the case when the 1-parameter family has more structure, and is in fact a 1-parameter subgroup (1-PS). 
In this case one can exploit the fact that the stabilizers come with a natural grading to get a tighter relation between the stabizers of $f, g$ and the stabilizers of the tangent of approach and the tangent of exit. 

\subsection{Stabilizers of Limits of Forms under a 1-parameter family}
\label{sec:A(t)}

Let $X$ be a set of indeterminates and let $V$ be the $GL(X)$-module $Sym^d (X)$ and $f,g\in V$ be non-zero forms. 
Let ${\cal G}$ be the  algebra $gl(X)$, ${\cal K}\subseteq {\cal G}$ be the stabilizer of $f$, and its dimension be $k$. Let ${\cal H}$ be the Lie algebra stabilizing $g$ and its dimension be $r$. 

Suppose that we have an algebraic family $A(t)\subseteq GL(X)$, parametrized by $t$ such that $A(1)=e$, the identity element. Let  
 \[ f(t)=A(t).f = t^{a} g+t^{b} f_{b} + t^{b+1} f_{b+1}+\ldots + t^{D} f_{D}, \]
By normalizing the family $A(t)$, as in the proof of Prop. \ref{prop:family}, and using suitable powers of $t$, we may assume that:
\[ f(t)=A(t).f =  g+t^b f_b + t^{b+1} f_{b+1}+\ldots + t^D f_D, \]
with $f_b \neq 0$. 
We will write $f^+(t)$ for the sum $\sum_{i=b}^{D} t^{i} f_i$. Thus $\lim_{t\rightarrow 0} f(t)=g$ and $\lim_{t\rightarrow 0} t^{-b} f^+ (t)=f_b$. We call $f_b$ as the tangent of approach. 

\begin{assume}
\label{assume:traneverse}
{\bf Transversality Assumption}. We assume that the vector space spanned by $f_b ,\ldots,f_D$ intersects $T_g O(g)$ trivially, where $T_g O(g)$ is the tangent space of the $GL(X)$-orbit $O(g)$ at point $g$. 
\end{assume}

 Let ${\cal S}$ be a complement to ${\cal H}$ in ${\cal G}$. Note that  
 the tangent space $T_g O(g)$ of the $GL(X)$-orbit of $g$ is precisely ${\cal S}\cdot g$ and is of the same dimension as ${\cal S}$. Let $N$ be a complement to $T_g O(g)$ in $TV_g \cong V$. By the transversality assumption, we may  assume that $f^+ (t)\in N$ for all $t$. This completes the data required to compute the local model. 
 
For a Lie algebra element $\mf{g}=\mf{s}+\mf{h}$, with $\mf{s}\in {\cal S}$ and $ \mf{h} \in {\cal H}$, let us compute the action of $\mf{g}$ on $f(t)=g+f^+(t)$. Identifying the tangent space $TV_{g}$ with ${\cal S}\oplus N$ as in Section \ref{thetamap}, we have:
\begin{equation} \label{eqn:both}
(\mf{s}+\mf{h})\cdot f(t)=(\mf{s}+\lambda_{\cal S} ((1+\theta (f^+ (t)))^{-1} (\mf{h}\cdot  f^+ (t)),\lambda_N ((1+\theta (f^+ (t)))^{-1} (\mf{h}\cdot  f^+ (t)).
\end{equation}
Thus, for a generic $t_0 \in \C$,  every element $\mf{k}\in {\cal K}(t_0 )$, the Lie algebra of the stabilizer of $f(t_0 )$, may be expressed as $\mf{k}=\mf{h}+\mf{s}$ with $\mf{h}$ satisfying the condition below: 
 \begin{equation}
\label{sec:stab:N}     
 \lambda_N ((1+\theta (f^+(t_0 )))^{-1} ((\mf{h}\cdot f^+(t_0)) =0. 
\end{equation}
and $\mf{s}$, the ${\cal S}$-completion of $\mf{h}$ at the point $f^+ (t_0 )$. We then have the following:
\begin{lemma} \label{lemma:stabdim}
For a generic $t_0 \in \C$ the stabilizer ${\cal K}(t_0)$ of $f(t_0)$ is conjugate to the stabilizer ${\cal K}$ of $f$. The dimension of ${\cal K}(t_0)$ is also $k$. 
\end{lemma}

\noindent 
{\bf Proof}: If $t_0$ is so chosen so that $A(t_0 )$ is invertible, then $f(t_0 )$ is in the $GL(X)$-orbit of $f$ and the assertion follows. $\Box $

\hspace*{0.2cm}

{\em Our objective is to construct a {\em uniform} basis $\{ \mf{k}_1 (t) ,\ldots ,\mf{k}_k (t) \}$ so that for a generic $t_0 $, the elements $\mf{k}_1 (t_0 ),\ldots, \mf{k}_k (t_0 )$ is  a $\C$-basis for ${\cal K}(t_0 )$.}

\hspace*{0.2cm}

Let $u_1 ,\ldots, u_m $ be a basis of $N$ and let us extend it by the elements $u_{m+1}=\mf{s}_1 \cdot g,\ldots,u_{m+p}=\mf{s}_p \cdot g$ to a basis of $TV_g =V$, where $\mf{s}_1,\ldots, \mf{s}_p $ is a basis of ${\cal S}$. This basis implements the isomorphism $TV_g =N \oplus TO_g$. 
 Let us denote these bases by $\mf{u}_N$ and $\mf{u}_{\cal S}$.
 
 Note that $(1+\theta (f^+ (t)))^{-1} $ is a map  $V \rightarrow V$ and we may represent it as the $(m+p)\times (m+p)$ matrix ${\cal I}(\theta , t)$  with entries in $\C(t)$ such that $(1+\theta (f^+(t)))^{-1}[ \mf{u}_N, \mf{u}_{\cal S}]=[ \mf{u}_N, \mf{u}_{\cal S}]\cdot {\cal I}(\theta , t) $. \eat{For simplicity the dependence of this matrix upon $f^+(t)$ is not reflected in our notation.}
 
 \begin{lemma} \label{lemma:Itheta}
 Let $\Delta (t)$ be the determinant of the matrix $I+\theta (f^+ (t))$ . Then $\Delta =1+c.t^{2b}+\mbox{ higher degree terms}$ (where $c\in \C$). 
 The matrix $\Delta \cdot {\cal I}(\theta,t)$ has entries in $\C[t]$ and is of the following form:
 \[ \Delta \cdot {\cal I}(\theta, t)= I - t^b \theta (f_b ) + \text{  higher degree terms}. \]
 \end{lemma}
 
 \noindent 
 {\bf Proof}: By linearity of $\theta$, we have:
 \[ \theta (f^+ (t))= t^b \theta (f_b )+ t^{b+1} \theta (f_{b+1})+ \text{ higher degree terms}. \]
 Let $J=I+\theta (f^+ (t))$ in the above basis, then $J_{ij}\in \C[t]$ with the following properties: (i) for $i\neq j$, we have $t^b$ divides $J_{ij}$, and (ii) $t^b$ divides $J_{ii}-1$. The first assertion on $\Delta (t)$, the determinant of $J$, follows. For the second assertion, again for $i\neq j$ each minor  $\tilde{J}_{ji} (t)$ of J is $-t^b\theta(f_b)_{ij}+\mbox{ higher terms}$, while for $i=j$, we have $\tilde{J}_{ii}=1-t^b\theta(f_b)_{ii}$ +\mbox{higher terms}. This proves the second part. $\Box $
 
Thus, if we express ${\cal I}(\theta ,t)$ as a matrix with entries in $\C [[t]]$, then: 
\[ {\cal I}(\theta,t)= I- t^b \theta (f_b )+ \mbox{higher terms}. \]

Let ${\cal M}_N (t): {\cal H}\rightarrow N$ be the map $\mf{h}\rightarrow \lambda_N \circ (1+\theta (f^+(t)))^{-1} (\mf{h}\cdot f^+(t))$ and ${\cal M}_{\cal S}(t):{\cal H}\rightarrow {\cal S}$ be the map $\mf{h} \rightarrow \lambda_{\cal S} \circ (1+\theta (f^+(t)))^{-1} (\mf{h}\cdot f^+(t))$. 
Note that ${\cal M}_N$ and ${\cal M}_{\cal S}$ are linear maps which depend on the parameter $t$. The kernel $ker({\cal M}_N (t))$ determines the ${\cal H}$-part of every element of ${\cal K}(t)$, while ${\cal M}_{\cal S}$ is the ${\cal S}$-completion at the point $f^+ (t)$.  

Let us now compute matrices for $M_N$ and $M_{\cal S}$ for ${\cal M}_N$ and ${\cal M}_{\cal S}$, respectively.  

 Let $\mf{h}_1,\ldots, \mf{h}_r$ be a basis of ${\cal H}$ and let us denote this sequence by $\overline{\mf{b}}$.
 Let us form the matrix $M_N$, an $m\times r$ matrix with entries in $\C(t)$ as follows. Let ${\cal M}_N (\mf{h}_i )=\sum_{k=1}^m (M_N)_{ki} u_k$.  Similarly, let the $p\times r$-matrix $M_{\cal S}$ be defined as follows: Let $ {\cal M}_{\cal S} (\mf{h}_i )=\sum_{k=1}^p (M_{\cal S})_{ki} \mf{s}_{k}$. We collect this information in the lemma below:
 \begin{lemma} \label{lemma:Mstruct}
In the above notation, we have:
\begin{enumerate}
    \item The matrices $M_N (t)$ and $M_{\cal S} (t)$ satisfy the following equations:
    \[ \mf{u}_N M_N (t) = \lambda_N \cdot {\cal I}(\theta ,t) \cdot (\overline{\mf{b}} \cdot f^+ (t))=[{\cal M}_N (\overline{\mf{b}})] \: \:\mbox{ and } \: \mf{u}_{\cal S} M_{\cal S} (t) = \lambda_{\cal S} \cdot {\cal I}(\theta ,t) \cdot (\overline{\mf{b}}\cdot f^+ (t)) \]
    \item Let ${\cal K}(t_0)$ be the stabilizer of $f(t_0 )$ for some generic $t_0 \in \C$ and let $\mf{k}\in {\cal K} (t_0 )$. Express $\mf{k}$ as $\mf{k}=\mf{h}+\mf{s}$ and let $\mf{h}=\overline{\mf{b}}\cdot \overline{\alpha} $, i.e., a linear combination $\overline{\alpha}$ of the basis elements $\overline{\mf{b}}$. Then
    \[ \mf{u}_N M_N (t_0)  \overline{\alpha}=0 \mbox{  and  } \mf{s}\cdot g=-\mf{u}_{\cal S} M_{\cal S} (t_0) \overline{\alpha}  .\]
    \item Conversely, any element $\mf{k}=\mf{h}+\mf{s}$ with $\mf{h}=\overline{\mf{b}}\cdot \overline{\alpha}$, satisfying the above conditions, is an element of ${\cal K}(t_0)$. 
\end{enumerate}
\end{lemma}

\noindent 
{\bf Proof}: All the assertions follow from the construction of $M_N$ and $M_{\cal S}$ and the requirements of Eqn. \ref{eqn:both}. $\Box $
 
\begin{lemma} \label{lemma:Mstruct2}
\begin{enumerate}
\item[(i)] Let $\mf{R}\subseteq \C (t)$ be the localization of $\C[t]$ at $\Delta$. In other words, $\mf{R}= \{ p/\Delta^d | p\in \C[t] \mbox{ and } d\geq 0\}$. Then, the entries of $M_N (t)$ and $M_{\cal S} (t)$ are elements of $\mf{R}$. 
\item[(ii)] $t^b$ divides the matrix $M_{\cal S} (t)$ and $M_N (t)$ in $\mf{R}$. In particular, $M_{\cal S} (0)=0$. 
\item[(iii)] For any generic $t_0 \in \C$, the rank of $M_N (t_0)$ is a constant and equals the $\C(t)$-rank of $M_N (t)$. This equals $dim({\cal H})-dim({\cal K})$, i.e., $r-k$. 
\item[(iv)] The $\C (t)$-rank of $M_{\cal S} (t)$ is zero iff ${\cal H} \cdot f(t) \subseteq N$.
\end{enumerate}
\end{lemma} 

\noindent 
{\bf Proof}: The matrix $M_N (t)$ is the product of operations expressed in appropriate bases, viz., (a) $\lambda_N :V \rightarrow N$, (b), $(1+\theta (f^+ (t))^{-1}:V \rightarrow V$ and (c) $ (\overline{\mf{b}}\cdot f^+ (t)):{\cal H}\rightarrow V$. Let us look at each one in turn. The first is a projection and is a matrix over $\C$ of rank $m+p$. 
For (b), this is precisely the matrix ${\cal I}(\theta ,t)$ which we have shown has entries in $\mf{R}$. Moreover, we know that ${\cal I}(\theta ,0)$ is the identity matrix. 
Finally, since ${\cal H}$ acts linearly on $V$ and $t^b$ divides $f^+ (t)$, the matrix $(\overline{\mf{b}}\cdot f^+ (t))$ has entries in $t^b \cdot \C[t]$.  Thus the product, viz., $M_N (t)$ has entries in $\mf{R}$ and $t^b$ divides it (in $\mf{R}$). The same reasoning holds for $M_{\cal S}$ as well. This proves (i) and (ii).  

For (iii), by lemma \ref{lemma:Mstruct}, the dimension of ${\cal K}(t_0)$ for a generic $t_0 $ equals $k$ and this must be the nullity of $M_N (t_0 )$. But this is the same as nullity of the matrix $M_N (t)$ as a matrix over $\C(t)$. Thus $k=r-rank_{\C(t)} (M_N (t))=r-rank_{\C} (M_N (t_0 ))$. This proves (iii).  

Finally, for (iv), the condition that $M_{\cal S}$ is identically zero is equivalent to ${\cal M}_{\cal S}$ being zero, which is equivalent to ${\cal H}\cdot f(t) \subset N$. $\Box $
 
Lemma \ref{lemma:Mstruct2} now allows us to construct ${\cal K}(t)$ as the column annihilator space of $M_N (t)$. Since $rank(M_N (t))=r-k$, we may choose column vectors $\overline{\alpha}_1 (t) ,\ldots, \overline{\alpha}_k (t)\in \C (t)^r$ such that $M_N \cdot \overline{\alpha}_i =0$. Since the entries of $M_N $ are $\mf{R}$, we may assume that these column annihilators may be chosen over $\C [t]$ and in a standard form as given by the following lemma: 

\begin{lemma} \label{lemma:normalize}
Let $A=[\overline{\alpha}_1 ,\ldots ,\overline{\alpha}_k ]$ be the $r \times k$-matrix with entries in $\C(t)$ and the columns $\overline{\alpha}_1 ,\ldots , \overline{\alpha}_k$ as above. 
Then we may assume that there is another basis (over $\C (t)$) $\overline{\alpha}'_1, \ldots ,\overline{\alpha}'_k $, for the column space of $A$ and a matrix $A' =[\overline{\alpha}'_1 ,\ldots ,\overline{\alpha}'_k ]$ such that all entries of $A'$ are in $\C[t]$ and $A'(0)$ is of rank $k$. 
Moreover, if $c$ is any column in the column space of $A$, such that $c(0)$ is defined, then $c(0)$ is in the column space of $A'(0)$. 
\end{lemma}

\noindent
{\bf Proof}: We first note that elementary column operations on $A$ do not change the corresponding ${\cal K}(t)$. Next, by multiplying by a suitable polynomial $p(t)$ and pulling out an appropriate power of $t$ from each column, we may assume that $A$ has entries in $\C [t]$ and that $A(0)$ exists and each of its column is non-zero. 
Finally, if a column, say, $\overline{\alpha}_k$ of $A$ is such that $\overline{\alpha}_k (0)=\sum_{i=1}^{k-1} \mu_i \overline{\alpha}_i (0)$, with $\mu_i \in \C$, then we may replace $\overline{\alpha}_k$ by $t^{-a} (\overline{\alpha}_k - \sum_{i=1}^{k-1} \mu_i \overline{\alpha}_i)$ for some suitable $a>0$, to get a new column of a lower degree. This process may be repeated till the assertions in the lemma are satisfied. Let us come to the second part. Let $R=\{ r_1 ,\ldots ,r_k \}$ be the row indices such that $ det(A'(0)[R])\neq 0$ and let $\Delta (t)=det(A'(t)[R])$. Since $c$ is in the column space of $A'$, there are polynomials $p_1 ,\ldots ,p_k \in \C[t]$ such that $\Delta (t) c=\sum_{i=1}^k p_i \overline{\alpha}'_i$. The result follows. $\Box $

Henceforth, we assume that the chosen $\overline{\alpha_1},\ldots, \overline{\alpha_k}$ have the properties assured to us by lemma \ref{lemma:normalize}. This has an immediate corollary.

\begin{prop} \label{prop:K0}\begin{enumerate}
\item There is a $\C (t)$-basis $\{ \mf{k}_i (t)\}_{i=1}^k$ of ${\cal K}(t)$, the stabilizer Lie algebra of $f(t)$ and a large number $D$ such that
\[ \mf{k}_i (t) = \sum_{j=0}^D (\mf{s}_{ij}+ \mf{h}_{ij})t^j \]
for suitable elements $\mf{s}_{ij} \in {\cal S}, \mf{h}_{ij}\in {\cal H}$.
\item The elements $\mf{s}_{ij}$ are zero for all $i$, and for all $j=0,\ldots,b-1$. As a result, we also have:
\[ \mf{k}_i (t) = \mf{h}_i (t)+ t^b\mf{s}_i (t)\]
where $\mf{h}_i (t) \in \C[t]\otimes {\cal H}$ and $\mf{s}_i (t) \in \C[t]\otimes {\cal S}$. 
\item The space ${\cal K}_0 =\C \cdot \{ \mf{k}_1 (0), \ldots, \mf{k}_k (0)\}$, the subspace formed by $\C$-linear combinations of the leading terms $\mf{k}_i (0) =\mf{k}_{i0}=\mf{h}_{i0}$ is a Lie subalgebra of ${\cal H}$ and of dimension $k$. Moreover, if $\mf{k}(t)\in {\cal K}(t)$ is any element such that $\mf{k}(0)$ is defined, then $\mf{k}(0)\in {\cal K}_0$. \item For any element $ \mf{h}\in {\cal K}_0$, we have $\lambda_N (\mf{h} \cdot f_b) =0$ and thus, there is an $\mf{s}\in {\cal S}$ such that $\mf{s}\cdot g +\lambda_{\cal S}(\mf{h}\cdot f_b)=0$
\end{enumerate}
\end{prop}

\noindent
{\bf Proof}: For (1), let us assume that $\overline{\alpha}_1 (t),\ldots ,\overline{\alpha}_k (t)$ are as in Lemma \ref{lemma:normalize}. We may then compute the element $\mf{h}_i (t) =\mf{b}\cdot \overline{\alpha}_i (t)$. Next, we define $\mf{s}_i (t)$ such that $\mf{s}_i (t)\cdot g =-\mf{u}_{\cal S} \cdot M_{\cal S} \overline{\alpha}_i (t)$.
Then $\mf{k}_i (t)=\mf{s}_i (t) +\mf{h}_i (t)$ stabilizes $f(t)$. 
Note that these are elements of ${\cal G}\otimes \C[t]$, i.e., linear combinations of Lie algebra elements with coefficients in $\C[t]$. Assertion (1) is an expansion of these as a polynomial in $t$ and $D$ is a bound on the degree of these polynomials. Assertion (2) follows from the fact that the element $\mf{s}_i (t)=\sum_{j=0}^D \mf{s}_{ij}t^j $
arises from the matrix equation $\mf{s}_i (t)\cdot g=-\mf{u}_{\cal S}  M_{\cal S} \cdot \overline{\alpha_i} (t)$, and that $t^b$ divides $M_{\cal S}$ as in lemma \ref{lemma:Mstruct2} (ii).   

Let us now come to (3). Since ${\cal K}(t)$ is a Lie algebra, for any two elements $\mf{k}_i (t) ,\mf{k}_j (t)$, we have:
\[ \begin{array}{rcl}[ \mf{k}_i (t),\mf{k}_j  (t)]&=& [\mf{k}_{i0},\mf{k}_{j0}]+ (
[\mf{k}_{i0},\mf{k}_{j1}]
[\mf{k}_{i1},\mf{k}_{j0}]) t^1 + \ldots \\
\end{array}
\]
This implies that $\lim_{t\rightarrow 0} [\mf{k}_i (t),\mf{k}_j  (t)]$ exists and equals $[\mf{k}_{i0},\mf{k}_{j0}]$. 
Next, note that lemma \ref{lemma:normalize} implies that ${\cal K}(t)$ is a $k$-dimensional subspace of $gl(X)\otimes \C (t)$ with the basis $\{ \mf{k}_1 (t), \ldots, \mf{k}_k (t)\}$ such that $\{ \mf{k}_1 (0), \ldots, \mf{k}_k (0)\}$ continue to be linearly independent over $\C$ in $gl(X)$. Thus we may express $[\mf{k}_i (t),\mf{k}_j (t)]$ in this basis as 
\[ \begin{array}{rcl}[ \mf{k}_i (t),\mf{k}_j  (t)]&=& \sum_{u=1}^k \beta^{ij}_u (t) \mf{k}_u (t), \\
\end{array}
\]
where $\beta^{ij}_u (t)\in \C(t)$. The condition on $\Delta (t)$ in the proof of lemma \ref{lemma:normalize} give us that $\beta_{ij}(0)$ exist. Taking limits on both sides gives us:
\[ \begin{array}{rcl}[ \mf{k}_{i}(0),\mf{k}_{j}(0)  ]&=& \sum_{u=1}^k \beta^{ij}_u (0) \mf{k}_{u}(0). \\
\end{array}
\] 
This proves that ${\cal K}_0$ is a Lie algebra. Note that since $\mf{s}_{i0}=0$ for all $i$, the elements $\mf{k}_{i0}$ are elements of ${\cal H}$ and thus ${\cal K}_0 \subseteq {\cal H}$. The assertion that $\{ h_{i0} \}$ are linearly independent, and the second assertion, follow from Lemma \ref{lemma:normalize}. 

Finally, let us come to the (4). Recall that $\mf{b}$ is the ordered basis $\{ \mf{h}_i ,\ldots ,\mf{h}_r \}$. We see that the condition that
$\mf{h}_i (t)=\mf{b} \cdot \overline{\alpha}_i (t)$ where $M_N \cdot \overline{\alpha}_i (t) =0$ is tantamount to asserting that
\[ \mf{u}_N M_N \overline{\alpha}_i (t)= \lambda_N \cdot (1+\theta (f^+))^{-1} \cdot (\mf{h}_i (t)\cdot f^+ (t))=0,\]
and, setting $\mf{s}_i \cdot g =-\mf{u}_{\cal S} M_{\cal S} \overline{\alpha}_i$ that,
\[ \mf{s}_i (t) \cdot g + \lambda_{\cal S} \cdot (1+\theta (f^+ ))^{-1} \cdot (\mf{h}_i (t) \cdot f^+ (t)) =0.\] 
Both together imply that $\mf{k}_i (t)$ stabilizes $f(t)$. Now comparing coefficient of the lowest degree, i.e., $t^{b}$, on both sides of the two equations, gives us:
\[ \lambda_N \cdot \mf{h}_{i0}\cdot f_b =0\: \mbox{ and }\: \mf{s}_{i,b}\cdot g+\lambda_{\cal S} \cdot \mf{h}_{i0}\cdot f_b =0.\]
This finishes the proof of the proposition. 
$\Box$

\begin{remark}
Surprisingly, the fact that the next lowest degree term in $(1+\theta (f^+))^{-1}$ is $-t^b \theta (f_b)$ is neither used nor required. 
\end{remark}
 We now have the following proposition.

\begin{prop} \label{prop:struct}
Let ${\cal K}_0$ be as above. Then:
\begin{enumerate}
\item Let $N\subseteq V$ be as above. Define the $\star$ action of ${\cal H}$ on $N$ as below:
\[ \mf{h} \star  n = \lambda_N  (\mf{h} \cdot n) \]
Then $\star : {\cal H} \times N \rightarrow N$ is a Lie algebra action on $N$ and matches the action of the quotient module $\overline{N}=V/TO_{g}$. The isomorphism from $N$ to $\overline{N}$ is given by $n\rightarrow \overline{n}=n+TO_g$. 
\item The subalgebra ${\cal K}_0 \subseteq {\cal H}$ stabilizes the element $\overline{f_b} \in \overline{N}$. In other words, ${\cal K}_0 \subseteq {\cal H}_b \subseteq {\cal H}$, the stabilizer of $\overline{f_b}\in \overline{N}$.   
\end{enumerate}
\end{prop}

\noindent
{\bf Proof}: Part (1) is straightforward.

For part (2), note that Prop. \ref{prop:K0} (4) tells us that $\lambda_N (\mf{h}\cdot f_b )=0$, for any $\mf{h} \in {\cal K}_0$. Thus, $\mf{h}$ stabilizes $\overline{f_b}$. 
$\Box $

\begin{defn} \label{defn:Ht}
Let ${\cal K}(t)$ be the stabilizer of $f(t)$ and let $\mf{k}_i (t) =\mf{h}_i (t) +\mf{s}_i (t)$, for $i=1,\ldots, k$, be a basis for ${\cal K}(t)$ as above. We define the space ${\cal H}(t)$ as the $\C[t]$-space spanned by $\{ \mf{h}_i (t)|i=1,\ldots,k\}$. For a $t_0 \in \C$, let ${\cal H} (t_0 )$ be the $\C$-space spanned by the vectors $\{ \mf{h}_i (t_0) |i=1,\ldots,k \}$. Note that, by definition, ${\cal H}(0)={\cal K}_0$. 
\end{defn}

\begin{lemma}
\label{lemma:Ht0}
Let $t_0 \in \C$ be generic and $\mf{k}=\mf{s}+\mf{h}$ be an element of ${\cal K}(t_0)$, then $\mf{h}\in {\cal H}(t_0)$. 
\end{lemma}

\noindent  
{\bf Proof}: The assertion is straightforward and follows from the fact that for generic $t_0$, the rank of the matrix $A(t_0 )$ of lemma \ref{lemma:normalize} is $k$. $\Box $ 

We now come to the main theorem of this section which summarizes the above results. We will recall the definition of ${\cal S}$-completion. 

\begin{defn}
For any element $\mf{h}\in {\cal H}$, and $n\in N$, we define its ${\cal S}$ completion as the element $\mf{h}+\mf{s}$, with $\mf{s}\in {\cal S}$  such that $\lambda_{\cal S} (\mf{h}\cdot n)=\mf{s}\cdot g$. For a set ${\cal J}\subseteq {\cal H}$, we define the ${\cal S}$-completion of ${\cal J}$ as the collection of all ${\cal S}$-completions of every element $\mf{j}\in {\cal J}$. 
\end{defn}

Note that since the map ${\cal S}\rightarrow {\cal S}\cdot g$ is a bijection, the ${\cal S}$-completion of an element $\mf{h}\in {\cal H}$ always exists and is unique. 

\begin{theorem} \label{theorem:mainform}
Let $A(t)\cdot f =f(t)=g +t^b f_b + \mbox{  higher terms}, $ as above, with ${\cal K}(t)$ of generic dimension $k$ as its stabilizer, and $f_b$, the tangent of approach. Then 
\begin{enumerate}
    \item[(i)] There is a subalgebra ${\cal K}_0 $ of  ${\cal H}_b \subseteq {\cal H}$ the stabilizer of $\overline{f_b}$ for the $\star$-action of ${\cal H}$ whose dimension is the same as the dimension of ${\cal K}(t_0)$ for a generic $t_0\in \C$.
    \item[(ii)] There is a basis $\{\mf{k}_i (t) | i=1,\ldots ,k \}$ of ${\cal K}(t)$ such that $\{ \mf{k}_i (0)|i=1,\ldots,k\}$ is a basis for ${\cal K}_0$. Moreover, ${\cal K}_0 ={\cal H}(0)$. 
    \item[(iii)] Let ${\cal H}(t)$ be the $\C [t]$-span of the vectors $\mf{h}_1 (t),\ldots, \mf{h}_k (t)$. Then, for any generic $t_0 \in \C$, the subspace ${\cal H}(t_0 )\subseteq {\cal H}$ has dimension $k$ and ${\cal K}(t_0)$ is the ${\cal S}$-completion of ${\cal H}(t_0 )$ for the point $f^+ (t_0 )$.   
\end{enumerate}
\end{theorem}

\noindent 
{\bf Proof}: Part (i) and (ii) are shown in Prop. \ref{prop:K0}. Part (iii) follows from lemma \ref{lemma:Ht0}.  This proves the theorem. $\Box $

\section{Projective limits under the action of 1-ps}
\label{lambdat}
We specialize the results of the previous section to the case when the 1-parameter family is a 1-parameter subgroup $\lambda (t)$ (1-PS). Each $\lambda $ effectively partitions the variable set $X$ into several parts $X= X_1 \cup \ldots \cup X_{k}$ with $\lambda $ operating on $X_i$ simply as $\lambda (t)(X_i)=t^{d_i} X_i$.
In Section~\ref{sec:1pssimple} we tackle the two block case, when the variable set $X$ is partitioned into two parts $Y$ and $Z$, with the 1-ps acting trivially on $Y$ and scaling by $t$ the variables in $Z$. Proposition~\ref{prop:stabgeneric} shows that in this case either there are elements which stabilizes $f, g$ and the tangent of approach and exit, i.e., a triple stabilizer, or the limiting Lie algebra ${\cal K}_0$ is nilpotent. Next, we study the general 1-ps case in Section~\ref{sec:1psgeneral}. We prove Proposition~\ref{prop:stabgeneral}, which is similar in spirit to that proved in the two block 1-ps.

\subsection{Two block $\lambda $}
\label{sec:1pssimple}
Let us partition $X$ as $X=Y\cup Z$. Let $\lambda : \C^* \rightarrow GL(X)$ be a 1-PS such that $\lambda (t) y=y$ for all $y\in Y$ but $\lambda (t) z= t z$. Let us denote by $V_i $ the space $Sym^{d-i} (Y)\otimes Sym^i (Z)$ and note that $Sym^d (X)=\oplus_{i=0}^d V_i$. Moreover, for any $v_i \in V_i$, we have $\lambda (t) v_i =t^i v_i $. Let $f$ and $g$ be forms so that:
\[ \lambda (t) f =f(t)= t^a g +t^b f_b +\ldots +t^d f_d \]
with $f_b \neq 0$. In other words, $g\in V_a$ is the leading term for the action of $\lambda $ on $f$ and $f_b$ is the tangent of approach. The family $A(t)=t^{-a/d}\lambda(t)$ is a suitable family and satisfies the conditions specified in the previous section. Whence, if the transversality assumption holds, then all the results, including Theorem \ref{theorem:mainform} are available. The special structure of $A(t)$ as above allows some more observations. We continue to use the same notation ${\cal K}(t)$, ${\cal H}(t)$ from the previous section.

First note that $\lambda (t)$ acts on ${\cal G}$ by conjugation whence, we have the weight space decomposition ${\cal G}_{-1}\oplus {\cal G}_0 \oplus {\cal G}_1$, where ${\cal G}_{-1}=Hom(Z,Y)$, ${\cal G}_0 =Hom(Y,Y)\oplus Hom(Z,Z)$ and ${\cal G}_1 = Hom(Y,Z)$. 
Note that $\lambda(t) \mf{g}_i \lambda(t)^{-1} =t^{i} \mf{g}_i$ for any $\mf{g}_i \in {\cal G}_i$, and that ${\cal G}_i \cdot V_{j} =V_{i+j}$. 

Since $g\in V_a$ and is homogeneous in the degree in $Y$, we have $\lambda (t) g =t^a g$ is the leading term of $\lambda (t) f$. Since the Lie algebra stabilizers of $t^a g$ and $g$ are the same, we have:
\begin{lemma}\label{lemma:basic}
In the above notation, the stabilizer ${\cal H}$ of $g$ is graded, i.e., ${\cal H}=\oplus_{i=-1,0,1} {\cal H}_i$ such that $\lambda (t) {\cal H}_i \lambda (t)^{-1} =t^i {\cal H}_i$. Moreover, we may choose ${\cal S}_i$ as the complement of ${\cal H}_i \subseteq {\cal G}_i$, then ${\cal S}=\oplus_{i=-1,0,1} {\cal S}_i$ is a complement of ${\cal H}$ in ${\cal G}$, and $\lambda (t) {\cal S}_i \lambda (t)^{-1} \subseteq t^i {\cal S}_i$.  Moreover, let $\pi_{\cal S}$ and $\pi_{\cal H}$ be the projections from ${\cal G}={\cal S}\oplus {\cal H}$, to the first and second components, then $\lambda (t)$ intertwines with these projections. In other words, for any $\mf{g}\in {\cal G}$, we have $\pi_{\cal S}( (\lambda (t) \mf{g} \lambda (t)^{-1}) = \lambda (t) \pi_{\cal S} (\mf{g}) \lambda (t)^{-1}$ and $\pi_{\cal H}( (\lambda (t) \mf{g} \lambda (t)^{-1}) = \lambda (t) \pi_{\cal H} (\mf{g}) \lambda (t)^{-1}$ 
\end{lemma}

Thus, $T_g O(g) ={\cal S}\cdot g \subseteq V_{a-1}\oplus V_a \oplus V_{a+1}$. The transversality condition translates into the requirement that $f_b \not \in T_g O(g)$. 

In view of Lemma~\ref{lemma:basic} and the
resulting discussion, the stabilization conditions of Theorem \ref{theorem:mainform} can be further studied. We present this
analysis in the Appendix~\ref{sec:stabcond}.

\begin{lemma} \label{lemma:uvprod}
Let $u,v\in \C$ be generic. Then $\lambda (u) {\cal H} (v) \lambda (u)^{-1} \subseteq {\cal H}(uv)$. 
\end{lemma}

\noindent 
{\bf Proof}: Let ${\cal K}$ be the stabilizer of $f=f(1)$. Now ${\cal K}(u)=\lambda(u) {\cal K} \lambda(u)^{-1}$ and that $\lambda (uv)=\lambda(u)\lambda(v)$, tells us that ${\cal K}(uv)=\lambda(u) {\cal K}(v) \lambda(u)^{-1}$. Moreover, note that $\pi_{\cal H} ({\cal K} (t_0))={\cal H}(t_0)$, for any generic $t_0 \in \C$, and thus may be applied to ${\cal K}(uv)$ above. This proves the lemma. $\Box $

We also have the important:

\begin{prop}
\label{prop:gradedK0}
The algebra ${\cal K}_0$ has a weight space decomposition ${\cal K}_0 ={\cal K}_{0,-1} \oplus {\cal K}_{0,0} \oplus {\cal K}_{0,1}$, where ${\cal K}_{0,i} \subseteq {\cal H}_i$. 
\end{prop}

\noindent 
{\bf Proof}: Let $\{ \mf{h}_i | i=1,\ldots,k \}$ be a basis of ${\cal H}(1)$ and let $\mf{h}_i =\oplus_j \mf{h}_{ij}$ be its the weight space decomposition. We may then write:
\[ \mf{b}_i (t) =\sum_j t^j \mf{h}_{ij} \]
The previous lemma implies that $\{ \mf{b}_i (t) | i=1,\ldots,k \}$ is another $\C (t)$-basis for ${\cal H}(t)$. The leading terms of this basis would decompose by weights and generate ${\cal H}(0)={\cal K}_0$.  
$\Box $

\begin{lemma}
\label{lemma:ellcondition}
 \begin{enumerate}
\item  Let $\ell=\log_t (\lambda (t)) \in gl(X)$, i.e., $\ell=diag([0,\ldots,0,1,\ldots, 1])$ in a suitable basis.Then for any $\mf{g}\in {\cal G}_i$, $[\ell , \mf{g}]=i \mf{g}$. For any $v\in V_i$, $\ell \cdot v= i\cdot v$. 
\item The algebra ${\cal H}_b$ is graded. Let $\ell' =\log_t (t^{-a/d} \lambda (t))=\ell -a/d \cdot I$. Then $\ell' \in {\cal H}$, $\ell'$ normalizes both ${\cal H}_b$ and ${\cal K}_0$, but $\ell' \not \in {\cal H}_b$. 
\end{enumerate}
\end{lemma}

\noindent 
{\bf Proof}: The first part of (1) is straightforward. For the action of $\ell$ on $V$, note that $\ell$ is also the differential operator $\sum_{z\in Z} z \frac{\partial}{\partial z}$. Thus $\ell v =i\cdot v$ for $v\in V_i$. Coming to (2), note that $\overline{f_b}$ is homogeneous and ${\cal H}$ is graded, hence ${\cal H}_b$ is graded. $[\ell' , {\cal H}_{b,i} ]\subseteq {\cal H}_{b,i}$. By the same token $[\ell' ,{\cal K}_{0,i} ]\subseteq {\cal K}_{0,i}$. By construction, $\ell' \cdot \overline{f_b} =(b-a) \overline{f_b } \neq 0$. $\Box $

The graded case also allows us to analyse the {\em tangent of exit}. 
\begin{defn}
Let $\lambda(t)$ be a 1-PS acting on the form $f$ and let:
\[ f(t)=\lambda (t) \cdot f = t^a g +t^b f_b +\ldots +t^D f_D \]
Then the tangent of exit is the form $\lim_{t\rightarrow 1} \frac{f(t)-f(1)}{t-1}$. 
\end{defn}
Since
\[ \ell f =af_a +bf_b +\ldots + Df_D \]
the tangent of exit is $\ell f -f$. Note that both $f$ and $\ell f$ are elements of $TO_f$, the tangent space of the orbit $O(f)$ at the point $f$. Thus, the tangent of exit is given (upto addition of a suitable multiple of the form $f$) by the action $\ell f $, where $\ell \not \in {\cal K}$. 

An important question is the stabilizer of $\ell f$ within $TO_f$ under the action of ${\cal K}$, the stabilizer of $f$. This is answered by the following lemma:

\begin{lemma} \label{lemma:pure}
$\Pi_i :{\cal G} \rightarrow {\cal G}_i$ be the weight projections as per $\lambda $. Then the stabilizer ${\cal K}_{\ell f} \subseteq {\cal K}$ within the stabilizer of $f$  consists of precisely those elements $\mf{k} \in {\cal K}$ such that $[\mf{k}, \ell] \in {\cal K}$. In particular, it contains all elements of ${\cal K}$ which are pure in weight, i.e., elements $\mf{k}$ such that $\Pi_i (\mf{k})=i\cdot \mf{k}$ for some $i$. Thus, these elements $\mf{k} \in {\cal K}$ of pure weight are {\bf triple stabilizers}, i.e., members of ${\cal K}\cap {\cal H} \cap {\cal K}_{\ell f}$. 
\end{lemma}

\noindent
{\bf Proof}: Let $\mf{k}\in {\cal K}$ stabilize $\ell f$. Then, since $\mf{k}f=0$, we have $\mf{k} \ell f = [ \mf{k} , \ell ]f+\ell \mf{k} f=[\mf{k} ,\ell] f $. This, in turn, is equivalent to $[\mf{k}, \ell]\in {\cal K}$. This proves the first assertion. For the second, note that if $\mf{k}\in {\cal K}$ is of pure weight, then $[\mf{k}, \ell]= i \cdot \mf{k}$ is certainly an element of ${\cal K}$. The last assertion comes from the fact that pure weight elements are also elements of ${\cal K}_0 \subseteq {\cal H}$. Note $\Box $

\begin{lemma} 
Pure elements $\mf{k}\in {\cal K}$ stabilize all components in the decomposition:
\[ f= f_a +f_b +\ldots + f_D \]
\end{lemma}

\noindent 
{\bf Proof}: Let $[\mf{k},\ell]=r\mf{k}$. We claim that $\mf{k} \ell^k f=0$ for all $k\geq 0$, which we prove by induction on $k$. The assertion for $k=0$ is obvious. Then, for general $k$, we have:
\[ \begin{array}{rcl}
0=\mf{k} \ell^k f &=& \ell \mf{k} \ell^{k-1} f + r \mf{k}\ell^{k-1} f \\
&=& 0 +0
\end{array} 
\]
This proves the assertion. 
Next, we apply $\mf{k}\ell^k$ to the above equation to get:

\[ \begin{array}{rcl}
0=\mf{k} (\ell^k f)&=& \mf{k} (\ell^k f_a)+ \mf{k} (\ell^k f_b)+\ldots + \mf{k} (\ell^k f_D ) \\
&=& \mf{k}(a^k f_a +b^k f_b +\ldots + D^k f_D)\\
&=& a^k (\mf{k} f_a) + b^k (\mf{k} f_b) +\ldots + D^k (\mf{k} f_D) 
\end{array}
\]
Since all of these are zero, we must have $\mf{k} f_c =0$ for all $c$. $\Box$

Thus, pure elements stabilize limits, tangents of approach and exit and all components in between. That there is an element $\mf{k} \in {\cal K}$ which is pure with respect $\ell$ is a significant alignment between ${\cal K}$ and $\ell$. 

\begin{remark} 
\label{remark:exit}
The above condition also connects with the weight structure of ${\cal K}_0$. Let $\mf{k}\in {\cal K}$ be such that its weight space decomposition is $\mf{k}=\mf{k}^{-1}+ \mf{k}^0 +\mf{k}^1$. We define  
\[ {\cal K}^{\geq 0} =\{ \mf{k}\in {\cal K} | \mf{k}^{-1}=0 \} \mbox{   and    } 
{\cal K}^{\geq 1} =\{ \mf{k}\in {\cal K} | \mf{k}^0 =0, \mf{k}^{-1}=0 \} \]
Note that ${\cal K}^{\geq 1} \subseteq {\cal K}^{\geq 0} \subseteq {\cal K}^{\geq (-1)}={\cal K}$ is a nested sequence of subalgebras such that ${\cal K}^{\geq i} /{\cal K}^{\geq (i+1)} \cong ({\cal K}_{0})_i$. Since $({\cal K}_0)_i \subseteq {\cal H}_i$, we have $dim(({\cal K}_0)_i ) \leq dim({\cal H}_i)$. 
\end{remark}

In Section \ref{subsection:det3} we work out this decomposition for various 1-PS acting on the form $det_3$, the $3\times 3$-determinant. 

\begin{prop} \label{prop:unistable}
Let $f,g, \lambda $ and $\ell$ be as above. Let $P(\lambda )$ be defined as below:
\[ P(\lambda )=\{ g \in G | \lim_{t\rightarrow 0} \lambda(t) g \lambda (t)^{-1} \mbox{ exists} \} \]
and let $U(\lambda)$ be its unipotent radical. Suppose that $\mf{k} \in {\cal P}(\lambda) \cap {\cal K}$ is a semi-simple element, then there is a unipotent element $u\in U(\lambda)$ such that:
\begin{enumerate}
\item If $f^u = u\cdot f$, then:
    \begin{equation} \label{eqn:apply}        \lambda (t) f^u = t^{a} g + \mbox{ higher terms }
    \end{equation}
    Thus, $g$ is the limit of $f^u$ under $\lambda$. 
    \item The element $\mf{k}^u =u\mf{k}u^{-1}$ stabilizes $f^u$, $g$ and $\ell f^u$. Thus, $\mf{k}^u$ is a triple stabilizer of the limit $g$, $f^u$ and the tangent of exit $\ell f^u$. 
    \end{enumerate}
\end{prop}

\noindent 
{\bf Proof}: Let $L(\lambda )$ be the elements of $GL(X)$ which commute with $\lambda $. Then $P(\lambda )=L(\lambda ) U(\lambda )=U(\lambda ) L(\lambda )$ is a Levi factorization, with $L(\lambda )$ as a reductive complement. 
Applying this to the corresponding Lie algebras, for any semisimple $\mf{k}$ as above, there is a $u\in U(\lambda)$ so that $\mf{k}^u=u\mf{k} u^{-1} \in {\cal L}(\lambda)$ and thus is an element of pure weight. Moreover $\mf{k}^u $ stabilizes $f^u =u\cdot f$.
Since $u\in U(\lambda)$ is unipotent, for any $h_c\in V_c$, the element $u\cdot h$ has the following weight decomposition with respect to $\lambda (t)$:
\[ uh_c =h_c +h'_{c+1} +\mbox{higher terms} \]
Now since $f$ has the weight decomposition:
\begin{equation} 
\label{eqn:apply2}       
f = g_a + f_b + \mbox{ higher terms }
    \end{equation}
on applying $u$ to Eq. \ref{eqn:apply2} we have:
\[ \begin{array}{rcl}
u f &=&  u\cdot g+u\cdot f_b +\mbox{higher order terms} \\
&=&  g + f'_{a+1} +  \mbox{ higher terms} \\
\end{array} \]
Thus, $g$ continues to be the leading term of $f^u$. Whence, we have:
\[ \lambda (t) (uf)= t^{a} g+ t^{a+1} f'_{a+1} + \mbox{higher terms} \]
Thus for the data $f^u , \lambda , \ell $ and $g$, we have that $g$ is the limit of $f^u$ under the action of $\lambda $. Now note that $\mf{k}^u \in {\cal K}_{f^u}$ and is of pure weight. Thus, by lemma \ref{lemma:pure}, $\mf{k}^u$ stabilizes the tangent of exit $\ell f^u$ and $g$. $\Box$ 

\begin{prop} \label{prop:stabgeneric}
Let $f,g, \lambda $ and $\ell$ be as above. Then at least one of the following hold: 
\begin{enumerate}
    \item[{\bf (A)}] ${\cal K}_0 $ is a nilpotent algebra, or
    \item[{\bf (B)}] there is a unipotent element $u\in U(\lambda)$ and an element $\mf{k} \in {\cal K}$ such that $g$ is a limit of $f^u$ under $\lambda $ and $\mf{k}^u$ is a triple stabilizer for the data $(f^u,\ell f^u ,g)$. 
\end{enumerate}

\end{prop}

\noindent 
{\bf Proof}: Note that ${\cal G}_0 \oplus {\cal G}_{1}={\cal P}(\lambda )$ and ${\cal G}_{1} ={\cal U}(\lambda)$, the Lie algebras of $P(\lambda )$ and $U(\lambda )$, and that the above decomposition is a Levi decomposition of ${\cal P}(\lambda)$ with ${\cal L}(\lambda )={\cal G}_0$ as a semisimple complement. 

\noindent 
{\bf Case 1}: Suppose that $dim(\Pi_{-1} ({\cal K})) =dim({\cal K})$, i.e., for every $\mf{k}\in {\cal K}$, if $\mf{k}\neq 0$ then so is $\Pi_{-1} (\mf{k})$. If this happens, then ${\cal K}_0 \subseteq {\cal G}_{-1}$ and ${\cal K}_0$ is nilpotent. Thus (A) holds. 

\noindent 
{\bf Case 2}: On the other hand, if 
$dim(\Pi_{-1} ({\cal K})) <dim({\cal K})$, then either:
\begin{enumerate}
\item[{\bf 2a}] There is an element $\mf{k} \in {\cal U}(\lambda) \cap {\cal K}$. In this case, $\mf{k}$ is a pure element and by Lemma \ref{lemma:pure}, $\mf{k}$ stabilizes $f,g$ and $\ell f$. Thus (B) holds with $u=1$.  \item[{\bf 2b}] All elements $\mf{k}\in {\cal P}(\lambda)$ may be written as $\mf{k}=\mf{k}_0 +\mf{k}_1$ with $\mf{k}_0 \neq 0$. Let us use Jordan's decomposition to the element $\mf{k}$ to express $\mf{k}=\mf{k}_{ss}+\mf{k}_n$, where $\mf{k}_{ss}$ is semisimple and $\mf{k}_n$ is nilpotent. Both are elements of ${\cal P}(\lambda ) \cap {\cal K}$. If there is a $\mf{k}_{ss} \neq 0$ as above, then by Prop. \ref{prop:unistable}, (B) holds. 
\item[{\bf 2c}] We come to the final case where for all $\mf{k}\in {\cal P}(\lambda)\cap {\cal K}$, we have $\mf{k}=\mf{k}_n=\mf{k}_0 +\mf{k}_1$ with $\mf{k}_0 \neq 0$. Since $\mf{k}_n^k=0$ for some $k$, we have $\mf{k}_0^k=0$ as well and thus $\mf{k}_0$ is nilpotent. Now, note that $\mf{k}_0$ is the leading term of $\mf{k}$ and thus $\mf{k}_0\in {\cal K}_0$. Whence we have ${\cal K}_0 =({\cal K}_0)_0 +({\cal K}_0)_{-1}$ where the general element $\mf{k}' \in {\cal K}_0$ may be written as $\mf{k}'_0 +\mf{k}'_{-1}$, where both are individually nilpotent. The grading ensures that $\mf{k}'$ is nilpotent too. Thus ${\cal K}_0$ is nilpotent, and (A) holds.
\end{enumerate}
$\Box$

\begin{remark}
The condition that ${\cal K}_0 \not \subseteq {\cal G}_1$ is tantamount to saying that $dim(\Pi_{-1} ({\cal K}))< dim({\cal K})$, i.e., ${\cal K}$ is not {\em generically placed} with respect to $\ell$. If indeed that is so, then all leading terms of ${\cal K}$ have weight $-1$ and ${\cal K}_0$ is nilpotent.  
\end{remark}

\subsection{General 1-PS}
\label{sec:1psgeneral}

We extend the results of the previous section to a general 1-PS $\lambda$, beyond the two-component case considered earlier. Thus, let $X=X_1 \cup \ldots \cup X_k$, with $|X_i |=n_i$, and let $\lambda (t) \in GL(X)=G$ be defined so that $\lambda (t) (x_i )=t^{d_i} x_i$. We call $\overline{d}=(d_1^{n_1},\ldots ,d_k^{n_k})$ as the type of $\lambda $. Thus, in the chosen basis of $(x_i )=X$, $\lambda (t)$ is the diagonal matrix $diag(t^{\overline{d}})$. As before, $V$ is a $G$-module via the map $\rho : GL(X)\rightarrow GL(V)$. We also have the Lie algebra map (also) $\rho : gl(X)\rightarrow gl(V)$. We now have $V=\oplus_{\chi} V^{\chi}$, the weight space decomposition of $V$. Note that $\chi \in \Z$. Let $\ell =\log_t (\lambda (t))$ as before. Then for $v^{\chi} \in V^{\chi}$, we have $\rho(\ell) \cdot v^{\chi} =\chi \cdot v^{\chi}$.  We also have ${\cal G}= \oplus_{\chi} {\cal G}^{\chi }$. 

\begin{defn}
For any ${\cal G}$-module $M$, we say that $m\in M$ is pure if $\ell \cdot m =\chi \cdot m$ for some $\chi \in \Z$. We say that $M$ is graded if $M=\oplus_{\chi} M^{\chi}$, the sum of its pure elements.
\end{defn}

\begin{prop} \label{prop:generalstab}
For the above situation, we have the following:
\begin{enumerate}
    \item The stabilizer ${\cal H} \subseteq {\cal G}$ of $x$ is graded. 
    \item Let ${\cal K}$ be the stabilizer of $y$, and ${\cal K}_0$ be the limit of ${\cal K}$ under the action of $\lambda$, then ${\cal K}_0$ is graded. 
    \item The tangent space $T_y O(y)$ of $y$ as an element of its $G$-orbit $O(y)$ is a ${\cal K}$-module. Let ${\cal K}_{\ell y} \subseteq {\cal K}$ be the stabilizer of the tangent vector $\ell y$. Then $\mf{k} \in {\cal K}$ stabilizes $\ell y$ if and only if $[\ell , \mf{k}] \in {\cal K}$. Consequently, if $\mf{k} \in {\cal K}$ is of pure weight, then $\mf{k} \in {\cal K}_{\ell y}$ as well.
    \end{enumerate}
\end{prop} 

\noindent 
{\bf Proof}: The arguments of Section \ref{sec:1pssimple} are easily modified for the general case. $\Box $

\noindent
Prop. \ref{prop:unistable} holds in the general case as well and Prop. \ref{prop:unistablegeneral} below is a restatement in the notation of this section. However, Prop. \ref{prop:stabgeneric}, on the structure of ${\cal K}_0$ does not; its modification is Prop. \ref{prop:stabgeneral}. 

\begin{prop} \label{prop:unistablegeneral}
Let $y,x, \lambda $ and $\ell$ be as above. Let $P(\lambda )$ be defined as below:
\[ P(\lambda )=\{ g \in G | \lim_{t\rightarrow 0} \lambda(t) g \lambda (t)^{-1} \mbox{ exists} \} \]
and let $U(\lambda)$ be its unipotent radical. Suppose that $\mf{k} \in {\cal P}(\lambda) \cap {\cal K}$ is a semi-simple element, then there is a unipotent element $u\in U(\lambda)$ such that:
\begin{enumerate}
\item If $y^u = u\cdot y$, then:
    \begin{equation} \label{eqn:apply}        \lambda (t) y^u = t^{a} x + \mbox{ higher terms }
    \end{equation}
    Thus, there is a $y^u$ in the orbit of $y$ such that $x$ is the limit of $y^u$ under $\lambda$. 
    \item The element $\mf{k}^u =u\mf{k}u^{-1}$ stabilizes $y^u$, $x$ and $\ell y^u$. Thus, $\mf{k}^u$ is a triple stabilizer of the limit $x$, $y^u$ and the tangent of exit $\ell y^u$. 
    \end{enumerate}
\end{prop}

The proof of this is a straightforward adaptation of the proof of Prop. \ref{prop:unistable}. 

\begin{prop} \label{prop:stabgeneral}
Let $y,x, \lambda $ and $\ell$ be as above. Then at least one of the following holds: 
\begin{enumerate}
    \item[{\bf (A)}] Let ${\cal K}'={\cal K}_0 \oplus \C \ell$, then ${\cal K}'$ is a Lie algebra of rank 1, i.e., the dimension of any maximal torus in ${\cal K}'$ is $1$, $[\ell , {\cal K}_0 ] \subseteq {\cal K}_0$, and every graded component $({\cal K}_0)_i$ is composed of nilpotent elements or
    \item[{\bf (B)}] there is a unipotent element $u\in U(\lambda)$ and an element $\mf{k} \in {\cal K}$ such that $x$ is a limit of $y^u$ under $\lambda $ and $\mf{k}^u$ is a triple stabilizer for the data $(y^u,\ell y^u ,x)$. 
\end{enumerate}

\end{prop}

\noindent 
{\bf Proof}: Note that ${\cal P} (\lambda )=\oplus_{\chi \geq 0} \: {\cal G}_{\chi}$ and ${\cal U}(\lambda)=\oplus_{\chi >0} \: {\cal G}_\chi$, the Lie algebras of $P(\lambda )$ and $U(\lambda )$, and that the above decomposition is a Levi decomposition of ${\cal P}(\lambda)$ with ${\cal G}_0$ as a semisimple complement. Next note that ${\cal K}'={\cal K}_0 \oplus \C \ell$ is indeed a Lie algebra since $[\ell ,{\cal K}_0] \subseteq {\cal K}_0$. Now we do a case analysis. 

\noindent 
{\bf Case 1}: Let $\Pi_{-}:{\cal G} \rightarrow \oplus_{\chi <0} \: {\cal G}_{\chi}$ be the projection. Suppose that $dim(\Pi_{-} ({\cal K})) =dim({\cal K})$. If this happens, then ${\cal K}_0 \subseteq \oplus_{\chi<0} {\cal G}_{\chi}$ and ${\cal K}_0$ is nilpotent. Even more, ${\cal K}' ={\cal K}_0 \oplus \C \ell$ is a Levi decomposition of ${\cal K}'$. Thus (A) holds. 

\noindent 
{\bf Case 2}: On the other hand, if 
$dim(\Pi_{-1} ({\cal K})) <dim({\cal K})$, then ${\cal K}\cap {\cal P}(\lambda) \neq 0$.

\begin{enumerate}
\item[{\bf 2a}]
If there is an $\ell' \in {\cal P}(\lambda )\cap {\cal K}$ such that $[\ell ' ,\ell]=0$, then it must be that $\ell' \in {\cal G}_0$, i.e., it must be of pure weight and therefore a triple stabilizer, and therefore (B) holds. 

\item[{\bf 2b}] If there is a semi-simple element $\mf{k} \in {\cal P}(\lambda) \cap {\cal K}$ then, by Prop. \ref{prop:stabgeneral}, there is a $u$ and a triple stabilizer as required. Thus (B) holds. 

\item[{\bf 2d}] We come to the final case where all $\mf{k}\in {\cal P}(\lambda)\cap {\cal K}$ are nilpotent. Let $\mf{k}_{\chi_0}$ be the leading term of $\mf{k}=\sum_{\chi \leq 0} \mf{k}_{\chi}$, i.e., the smallest degree $d$ such that $\mf{k}_d \neq 0$. If $\chi_0 >0$  then $\mf{k}_{\chi_0} \in {\cal K}_0$ is nilpotent. Otherwise, we may write $\mf{k}=\mf{k}_{0} + \sum_{\chi >0} \mf{k}_{\chi}$. Since $\mf{k}^k=0$ for some $k$, we have $\mf{k}_0^k=0$ as well and thus $\mf{k}_0$ is nilpotent. 
Thus, for every $\chi$, ${\cal K}_0 \cap {\cal G}_{\chi}$ consists of nilpotent elements. From Case 2a above, we already know that the rank of ${\cal K}'={\cal K}\oplus \C \ell$ is $1$ 
\end{enumerate}
This proves the proposition. $\Box$

\section{Closures of affine orbits}
\label{sec:codim1}
We apply the results of the last two sections to the closure of forms whose $SL(X)$-orbit is affine. As mentioned in the introduction (see also Appendix~\ref{appendix-gct}), the motivation for studying comes from geometric complexity theory. Examples of such forms are $det(X)$, the determinant of the matrix $X$, and $perm(X)$, its permanent. In the next Section, we cover some general results. In Section~\ref{subsection:det3} we apply these and earlier results to the $3 \times 3$ determinant form.

\subsection{General Results}

Let $f(X) \in Sym^d (X)$ be such that the stabilizer $K\subseteq GL(X)$ is reductive (of dimension $k$). This implies that the $SL(X)$-orbit of $f$ is an affine variety. Let $O(f)$ be the $GL(X)$-orbit of $f$ and $\overline{O(f)}$ be its projective closure. It is known that the closed set $\overline{O(f)} -O(f)$ is of co-dimension one. Moreover, there are forms $g_1 ,\ldots ,g_m$ such that it is contained in the union of the projective closures $\{ \overline{O(g_i )}\}_{i=1}^m$.

Separately, it is known that for any form $f$, forms $g$ which appear in the closure of the orbit $O(f)$ appear as leading terms under a 1-parameter substitutions $A(t)\cdot X$, where $A(t) \in GL(X)$ has entries which are polynomials in $t$. Thus, for the above $g_i$, there are linear transformations $A_i (t)$, with polynomial entries, such that $f(A_i (t)\cdot X)=t^d g_i + $ higher degree terms. It is also clear from dimension counting that for all $i$, the stabilizer $H_i$ of $g_i $ is of dimension $k+1$.

Our focus will be on one such those $g_i =g$. We begin with a lemma:
\begin{lemma} \label{lemma:invertible}
Given a transformation $A(t)\in gl(X)$ such that $g$ is the leading term of $f(t)=f(A(t) \cdot X)$, we may assume that $A(t)\in GL(X)$, i.e., $A(t)$ is invertible for most $t\in \C$.
\end{lemma}

\noindent 
{\bf Proof}: Let us consider $A'(t)=A(t)+t^D I$ for a large power $D>0$. Then $g$ continues to be the leading term of $f'(t)=f(A'(t)\cdot X)$. Moreover, $A'(t)$ is invertible for most $t\in \C$. $\Box $

By the above lemma, we may assume that $A(t)$ is generically invertible and that:
\[ f(t)=f(A(t)\cdot X)= t^a g + t^b f_b +\ldots + t^d f_d \]
 
Let ${\cal H}$ be the Lie algebra stabilizer of $g$ as before and ${\cal K}(t)=A(t){\cal K}A(t)^{-1}$ that of $f(t)$, with $dim({\cal K}(t))=k$. Let ${\cal K}_0$ be the limit of ${\cal K}(t)$ and note that $dim({\cal K}_0)=k$ and that ${\cal K}_0 \subseteq {\cal H}_b$, the stabilizer for $\star$-action of ${\cal H}$ on $\overline{N}$. We have the following lemma:
\begin{lemma}
The algebra ${\cal K}_0$ is precisely ${\cal H}_b$, the stabilizer of $\overline{f_b} \in \overline{N}$ under the $\star$-action of ${\cal H}$. 
\end{lemma}

\noindent 
{\bf Proof}: We have the containment ${\cal K}_0 \subseteq {\cal H}_b \subset {\cal H}$. Next, the orbit of $g$ is of co-dimension 1 within the orbit of $f(t)$. This, gives us that $dim({\cal H})=dim({\cal K}(t))-1=k+1$. Since ${\cal H}_b \neq {\cal H}$, we must have ${\cal H}_b ={\cal K}_0$. 
$\Box $  

\begin{lemma} \label{lemma:1PScodim1}
Suppose further that $A(t)$ is a one-parameter subgroup $\lambda (t)$. So we have
\[ \lambda(t) \cdot f=f(t) = t^a \cdot g + t^b \cdot f_b + \mbox{higher terms} \]
Let $\ell =\log (t^{-a/d} \lambda(t)) \in gl(X)$ be the Lie algebra element generating the 1-PS $t^{-a} \lambda(t)$. Suppose that $\overline{f_b}\neq 0$. Then (i) $\ell \in {\cal H}$ but does not stabilize $\overline{f_b}$, (ii) ${\cal K}_0={\cal H}_b$, and (iii) ${\cal H}={\cal K}_0 \oplus \C\cdot \ell$. 
\end{lemma}

\noindent
{\bf Proof}: Since $t^{-a} \lambda (t) \cdot g=g$ but $t^{-a} \lambda (t) f_b =t^{b-a} f_b$, (i) is clear. Thus $\ell \in {\cal H}-{\cal K}_0$. Since the dimension of ${\cal K}_0$ is one less than the dimension of ${\cal H}$, we must have ${\cal H}={\cal K}_0 \oplus \ell$. This proves the lemma. $\Box $

When $g$ is obtained as the limit of $\lambda (t)$ where $\ell $ has only two distinct entries (i.e., the 2-block case), we have the following result:
\begin{prop} \label{prop:codim1longchain}
Let $g$ be a projective limit of codimension $1$ of $f$ under $\lambda $ as above. Then at least one of the following conditions hold:
\begin{enumerate}
    \item ${\cal K}_0$ is nilpotent and ${\cal H}={\cal K}_0 \oplus \C \ell$ is a Levi factorization of ${\cal H}$. 
    \item There is a unipotent element $u \in U(\lambda)$ (see \ref{prop:stabgeneric}) and an element $\mf{k}\in {\cal K}$ such that $u\mf{k}u^{-1}$ is a stabilizer of $g,f^u$ and $\ell f^u $, and $g$ is the limit of $f^u$ under $\lambda$. 
\end{enumerate}
\end{prop}

\noindent 
{\bf Proof}: This is a consequence of Prop \ref{prop:stabgeneric} and the fact that ${\cal H}={\cal K}_0 \oplus \C \ell$. $\Box $

\subsection{The analysis of $det_3$, the $3\times 3$ - determinant.} \label{subsection:det3}

We consider here the set $X=\{ x_1 ,\ldots , x_9 \}$ and the $165$-dimensional space $Sym^3 (X)$ of homogeneous forms of degree $3$, acted upon by $GL(X)$. The $3\times 3$-determinant, $det_3 (X)$ is a special element and its definition is given below:
\[ det_3 (X)=det\left( \left[
\begin{array}{ccc}
x_1 & x_2 & x_3 \\
x_4 & x_5 & x_6 \\
x_7 & x_8 & x_9 \\ \end{array} \right] \right) 
\]
The obvious question are the codimension one forms in the orbit closure of the orbit of $det_3 (X)$. 

These have been identified by H\"uttenhain\cite{huttenhain2017geometric} and are given by the two forms $Q_1$ and $Q_2$ below:
\[ \begin{array}{rcl}
Q_1 (X)&=&det\left( \left[
\begin{array}{ccc}
x_1 & x_2 & x_3 \\
x_4 & x_5 & x_6 \\
x_7 & x_8 & -x_5 -x_1  \\ \end{array} \right] \right) \\ 
Q_2 (X) &=& x_4 x_1^2 +x_5 x_2^2 +x_6 x_3^2 + x_7 x_1 x_2 +x_8 x_2 x_3 +x_9 x_1 x_3 \end{array}
\]

We will study each in turn. Before we begin, let us look at the stabilizer ${\cal K}$ of $det_3 (X)$. This will serve as our form $f$, in the notation of this section.
Within $Sym^3 (X)$ the dimension of the $GL(X)$-orbit $O(det_3 (X))$ of $det_3 (X)$ is given by $dim(GL(X))-dim({\cal K}))$. It is known that this is composed of the transformations $X\rightarrow AXB$ with $det(AB)=1$ and the $X\rightarrow X^T$. The condition $det(AB)=1$ and that $(cA)X(c^{-1}B)$ lead to the same linear transformation on $X$, give us two independent constraints on $GL_3 \times GL_3 \rightarrow {\cal K}$. This gives us that
$dim({\cal K}(1))=18-2=16$ and $dim(O(det_3 (X))=65$.

\noindent
Coming to $Q_1$, we see that:
\[ det_3 (X)=Q_1 (X)+(x_1+x_5+x_9)(x_1 x_5 -x_2 x_4 ) \]
We denote by $Q_1'$ the form 
$(x_1+x_5+x_9)(x_1 x_5 -x_2 x_4)$.  
This motivates us to define $Y=\{ x_1 ,\ldots ,x_8 \}$ and $Z=\{ x_1 +x_5 +x_9\}$. We define $\lambda_1 (t) \in GL(X)$ as $\lambda_1 (t) x_i =x_i $ for $i=1,\ldots, 8$ and $\lambda (t)(z)=tz$, where $z=(x_1+x_5+x_9)$. We thus see that:
\[ \lambda_1 (t)\cdot det_3 (X)=Q_1 +t\cdot Q_1 ' \]
Thus in the notation of this section, $f=det_3 (X)$, $g=Q_1$ and $f_b = Q_1'$, with $a=0$ and $b=1$.

\begin{lemma}
The stabilizer of ${\cal H}_1$ of $Q_1$ within $gl(X)$ has dimension 17. Moreover, the stabilizer ${\cal H}_1$ has a direct sum decomposition into degree 0 and degree 1 parts, as in Section \ref{lambdat}. The tangent of exit is the form $Q'_1$ and the dimension of its stabilizer within ${\cal K}$, the stabilizer of $det^3(X)$, is $4$. 
\end{lemma}

\noindent
{\bf Proof}: Let us identify $Y$ as the space of all trace zero matrices $y$ labelled with their 9 entries minus the entry $y(3,3)$. Then, $Q_1$ is precisely the determinant of $y\in Y$. 
For any element $A\in GL_3$, the action of $A$ on $y$ given by $y\rightarrow AyA^{-1}$ preserves the space $Y$ of trace zero matrices as well as $Q_1$. Whence, we get a map $GL_3 \rightarrow {\cal H}_1$, the stabilizer of $Q_1$. 
Since this map has a  1-dimensional kernel (viz., the matrices $cI$, with $c\neq 0$), we have an $8$-dimensional image ${\cal R}_Y \subseteq {\cal H}_1$.  Next, let $Hom(Z,Y\oplus Z)$ be the collection of all linear substitutions for $z$ in terms of the variables of $Y$and $z$. Since $Q_1$ does not involve $z$, these substitutions leave $Q_1$ invariant. 
This gives us a map $Hom(Z,Y)\oplus Hom(Z,Z)\rightarrow {\cal H}_1$. We call this image as the spaces ${\cal Q}_1 \oplus {\cal R}_Z$. It can be shown that ${\cal Q}_1 , {\cal R}_Y$ and ${\cal R}_Z$ intersect trivially and thus ${\cal R}_Y \oplus {\cal R}_Z \oplus {\cal Q}_1 \subseteq {\cal H}_1$ to obtain a $17$-dimensional subalgebra. 
The result of \cite{huttenhain2017geometric} implies that the dimension can not be more, since we know that the dimension of the orbit of $Q_1$ is only one less than that of $det_3 (X)$, and hence its stabilizer dimension must be $17$. 

We then see that ${\cal H}_1 = {\cal R}_Y \oplus {\cal R}_Z \oplus {\cal Q}_1= {\cal R}\oplus {\cal Q}_1$ and splits by degree, i.e., $({\cal H}_1)_0 = {\cal R}$ and $({\cal H}_1)_{-1}={\cal Q}_1$. 

The stabilizer dimension of $Q'_1$ is through a direct computation. $\Box $ 

\begin{lemma}
We have $H_b ={\cal R}_Y \oplus {\cal Q}_1$ and ${\cal K}_0 ={\cal H}_b$. Let $\ell_1 =\log (\lambda (t))$, then 
$\ell_1 \in {\cal H}_1$ so that ${\cal H}_1 ={\cal H}_b \oplus \ell_1 $ and $[\ell_1 ,{\cal H}_b ]\subseteq {\cal H}_b$. Moreover, $\ell_1 \cdot det^3 (X)=Q'_1$, the tangent of exit, and ${\cal H}_b$ is an ideal of ${\cal H}_1$. Finally, as in Case (i), Section~\ref{sec:stabcond} with $b-a=1$, (i) for every element $\mf{r}\in {\cal R}_Y$, we have an element $\mf{s}$ in ${\cal S}_1$ such that $\mf{r}\cdot f_b =\mf{s} \cdot Q_1$, and (ii) for any $\mf{q}\in {\cal Q}_1$, there is an $\mf{s} \in {\cal S}_0$ such that $\mf{q}\cdot f_b =\mf{s}\cdot Q_1$. 
\end{lemma}

\noindent 
{\bf Proof}: The equality of ${\cal K}_0$ and $\Hb$ comes because $Q_1$ is one of the co-dimension 1 forms in the closure of $det_3 (X)$. The second assertion follows from Lemma \ref{lemma:1PScodim1}. The last set of assertions follow from the fact that $\Hb \cdot f_b \subseteq {\cal S}\cdot Q_1$. $\Box $
\begin{remark} 
The statement that in this case ${\cal H}_b$ is an ideal in ${\cal H}$ follows from a characterization of 
codimension 1 Lie subalgebras of a Lie algebra due to Hoffman~\cite{hofmann1965lie}. In Section~\ref{sec:cohomo} when we discuss connections to Lie algebra cohomology, Hoffman's theorem is stated as Theorem \ref{theorem:hoffman}.
\end{remark}

\begin{ex}
We do an example of an $\mf{r}$ and a corresponding $\mf{s}$ such that $\mf{r}\cdot f_b =\mf{s} \cdot Q_1$. However, constructing elements $\mf{r}\in {\cal R}_Y$ needs some care.  
Consider the matrix $E_{ij}$ such that $E_{ij} (i,j)=1$ and $E_{ij}(k,l)=0$ for all other tuples $(k,l)$. Let $Y$ be the matrix such that $det(Y)=Q_1$, as given below: 
\[ Y=\left[ \begin{array}{ccc}
x_1 & x_2 & x_3 \\
x_4 & x_5 & x_6 \\
x_7 & x_8 & -x_1 -x_5 \end{array} \right]
\]
We construct $[E_{23},Y]=E_{23}Y-YE_{23}$ as below:
\[ [E_{23} ,Y]= R_{23}=
\left[ \begin{array}{ccc}
0 & x_3 & 0 \\
-x_7 & x_6 -x_8  & x_1+x_5 \\
0 & -x_1 -x_5 & 0 \end{array} \right]
\]
The corresponding element $\mf{r}$ is to be constructed as the linear operator below:
\[ \mf{r}= x_3 \frac{\partial}{\partial x_2} -x_7 \frac{\partial}{\partial x_4} +(x_6 -x_8)\frac{\partial}{\partial x_5} +(x_1+x_5)\frac{\partial}{\partial x_6 } -(x_1+x_5)\frac{\partial}{\partial x_8} \]
We may check that $\mf{r}\cdot Q_1=0$. On the other hand, we have: 
\[ \mf{r}\cdot f_b = z((x_2 x_7 -x_1 x_8)- (x_3 x_4 -x_1 x_6))= (-z \frac{\partial}{\partial x_6}+z \frac{\partial}{\partial x_8} )\cdot Q_1 =\mf{s}\cdot Q_1\]
The second operator $\mf{s}$ is clearly an element of ${\cal S}_1$. 
\end{ex}

\noindent
Let us now come to $Q_2$.

That $Q_2 $ is in the orbit closure of $det_3 $ is shown by the following lemma.
 
 \begin{lemma} \cite{huttenhain2017geometric}
 Let $Y,Z$ be the generic matrices below and let $X=Y\oplus Z$. 
 \[ Y=\left[ \begin{array}{ccc}
 0 & x_1 & -x_2 \\
 -x_1 & 0 & x_3 \\
 x_2 & -x_3 & 0 \\
 \end{array} \right] \: \: \: \: 
 Z=\left[ \begin{array}{ccc}
 2x_6 & x_8 & x_9 \\
 x_8 & 2x_5 & x_7 \\
 x_9 & x_7 & 2x_4 \\
 \end{array} \right] 
 \]
Let $\lambda_2 (t)$ be such that $\lambda_2 (t)\cdot Y=Y$ and $\lambda_2 (t)\cdot Z=tZ$. Let us define $det^3 (X)$ as the determinant of the matrix $Y+Z$. Then:
\[ det^3(\lambda_2 (t)\cdot X))= det(Y+tZ)= tQ_2 +t^3 Q_3 \]
where:
\[ \begin{array}{rcl}
Q_2 (X) &=&  x_4 x_1^2 +x_5 x_2^2 +x_6 x_3^2 + x_7 x_1 x_2 +x_8 x_2 x_3 +x_9 x_1 x_3\\
Q_3 (X)&=& 8x_4 x_5 x_6 -2x_6x_7^2 -2x_4 x_8^2 -2x_5 x_9^2 +2 x_7 x_8 x_9 \end{array} \] 
 \end{lemma}
 
 The proof is a verification. Note that this parametrization of the generic matrix $X$, instead of the one in the case of $Q_1 $, is merely a reordering of the variables and hence the expression $det^3 (X)$ of the determinant is in the orbit of $det_3 (X)$, and hence the orbit closures of $det(X)$ and $det^3 (X)$ match. The expansion of $det^3 (\lambda_2 (t)\cdot X)$ puts us into the familiar situation where $Q_2$ is the limit, $f_b =Q_3$ with  and $a=1$ and $b=3$. 

\begin{lemma}
The stabilizer of ${\cal H}_2$ of $Q_2$ within $gl(X)$ has dimension 17.  The tangent of exit is the form $Q_3$ and the dimension of its stabilizer with ${\cal K}$, the stabilizer of $det^3(X)$, is $8$.
\end{lemma}
We may write $Q_2 $ as the inner product:
\[ Q_2 (X)=R \cdot C=\left[
\begin{array}{cccccc}
x_1^2 & x_2^2 & x_3^2 & x_1 x_2 & x_2 x_3 & x_1 x_3 \\ \end{array} \right] 
\left[ \begin{array}{c}
x_4 \\
x_5 \\
x_6 \\
x_7 \\
x_8 \\
x_9 \\
\end{array} \right] \]
Let $Y=\{ x_1 ,x_2 , x_3 \}$ and $Z=\{ x_4 ,\ldots ,x_9 \}$ and let $\rho : GL(Y)\rightarrow GL(Sym^2 (Y))$ be the representation of $GL(Y)$. We see that for any $A\in GL(Y)$, we have:
\[ Q_2 (X)= R \cdot \rho (A)^T (\rho(A)^{T})^{-1} \cdot C \]
Thus, we have an embedding $\alpha :GL(Y) \rightarrow {\cal H}_2 \subseteq GL(Y)\times GL(Z)$ of $GL(Y)$ into ${\cal H}_2 $, given by $A\rightarrow A \times (\rho(A)^T)^{-1}$. The image ${\cal R}=\alpha (GL(Y))$ is the reductive part of $H_2$ and has dimension equal to $dim(GL(Y))=9$. We may further factor ${\cal R}$ as follows. Let $D\subset gl(Y)$ be multiples of the identity matrix. We may then write $gl(X)=D\oplus sl(X)$ giving us a decomposition of ${\cal R}={\cal D}\oplus {\cal R}'$, where $\alpha (D)={\cal D}$ is a $1$-dimensional subspace of ${\cal R}$, and ${\cal R}'=\alpha (sl(X))$. 

Let us now come to the nilpotent part which is a subspace of $Hom(Z,Y)$. Let $A\in Hom(Z,Y)$ be a $6\times 3$-matrix and let us write the action of $A$ on $Q_2 $ as:
\[ A \cdot Q_2 (X)=R \cdot C=\left[
\begin{array}{cccccc}
x_1^2 & x_2^2 & x_3^2 & x_1 x_2 & x_2 x_3 & x_1 x_3 \\ \end{array} \right] 
\left( \left[ \begin{array}{c}
x_4 \\
x_5 \\
x_6 \\
x_7 \\
x_8 \\
x_9 \\
\end{array} \right] 
+
\left[\begin{array}{ccc}
a_{11} & a_{12} & a_{13} \\
& & \\
& \vdots & \\
& &\\
& \vdots & \\
a_{61} & a_{62} & a_{63} \\
\end{array} \right] 
\left[ \begin{array}{c}
x_1 \\
x_2 \\
x_3 \\
\end{array} \right] \right)
\]
The requirement that $A$ stabilizes $Q_2$ is tantamount to the condition that:
\[ p(A,Y)=R\cdot A \cdot \left[ \begin{array}{c}
x_1 \\
x_2 \\
x_3 \\
\end{array} \right] =0\]
Since $p(A,Y)\in Sym^3 (Y)$ and $dim(Sym^3 (Y))=10$, the condition $p(A,Y)=0$ translates to $10$ equations. Thus, we define ${\cal A}$ as follows:
\[ {\cal A}=\{A \in Hom(Z,Y)| p(A,Y)=0 \} \]
Then we have ${\cal A}\rightarrow {\cal H}_2$ and that $dim({\cal A})=dim(Hom(Z,Y))-10=8$. Thus ${\cal H}_2 ={\cal R}\oplus {\cal A}$ and $dim({\cal H}_2 )=8+9=17$. 
The stabilizer dimension of $Q_3$ is through a direct computation.
$\Box $
\begin{lemma}
We have $H_b ={\cal R}' \oplus {\cal A}$ and ${\cal K}_0 ={\cal H}_b$. Let 
$\ell_2 =\log (t^{-1/3}\lambda_2 (t))$, then $\ell_2 \in {\cal H}_2$ so that ${\cal H}_2 ={\cal H}_b \oplus \ell_2 $ and $[\ell_2  ,{\cal H}_b ]\subseteq {\cal H}_b$. 
Moreover, $\ell_2 \cdot det^3 (X)=Q_3$, the tangent of exit. ${\cal H}_b$ is an ideal of ${\cal H}_2$ and condition (3) of Theorem ~\ref{theorem:hoffman} holds. Finally, as in Case (ii) Section ~\ref{sec:stabcond}, with $b-a=2$, we have for every element $\mf{a}\in {\cal A}_Y$, there is an element $\mf{s}$ in ${\cal S}_1$ such that $\mf{a}\cdot f_b =\mf{s} \cdot Q_1$. 
\end{lemma}

\begin{remark}
The algebras ${\cal H}_2$ and ${\cal K}_0$ provide an extremely important insight into the possible structure of the Levi decomposition of ${\cal K}_0$. The coupling of the $End(Y)$ and $End(Z)$ components through different representations of the same algebra $gl(3)$, and the representation on $End(Z,Y)$ as one which interleaves between the two, provides an important template for the Lie algebra structure of ${\cal K}_0$. 
\end{remark}

Finally, we cover the case of a limit of $det^3(X)$ which is not one of the components of the orbit closure. Let $Y=\{ x_1 ,\ldots ,x_4 \}$, $Z=\{ x_5 ,\ldots ,x_9 \}$ so that $X=Y\cup Z$. Let $\lambda_3 (t):\C^* \rightarrow GL(X)$ be such that $\lambda_3 (t)\cdot x=x$ for all $x\in Y$, while $\lambda_3 (t)\cdot x= tx$ for all $x \in Z$. We may write:
\[ \lambda_4 (t) \cdot det^3 (X)= t[x_1 (x_5 x_9 -x_6 x_8 )+x_7 (x_2 x_6 -x_3 x_5 )]+t^2 [-x_4 (x_2 x_9 -x_3 x_8)] \]
We may write this as $tQ_4 +t^2 Q'_4$, where $Q_4 =x_1 (x_5 x_9 -x_6 x_8 )+x_7 (x_2 x_6 -x_3 x_5 )$ and $Q'_4 =-x_4 (x_2 x_9 -x_3 x_8)$. We see that $dim({\cal H})=21$ and exceeds $dim({\cal K})+1$. Moreover, ${\cal K}$ has pure elements of degree $-1$ and $1$ with respect to $\ell_4$ and ${\cal K}_{\ell_4 f}$ is of dimension $8$. 

\begin{remark}
The stabilizer ${\cal H}_4$ of $Q_4$ is of dimension $21$. The tangent of exit is $Q'_4$ and its stabilizer within ${\cal K}$ is of dimension $8$. Let $\ell_4 = \log (t^{-1/3} \lambda_4 (t))$, then $\ell_4 \in {\cal H}_4$ and $\ell_4 \cdot det^3 (X) =Q'_4$, the tangent of exit. 
\end{remark}
 We end this section with the following table of filtered dimensions for various 1-PS operating on the form $det_3$, the limit $g$ and the dimensions of the graded components of the stabilizers of the limit $g$ and the tangent of exit $\ell f$.  ${\cal K}_{\ell f}$ is the subalgebra of triple stabilizers. These calculations illustrate Lemma ~\ref{lemma:pure} and Proposition ~\ref{prop:gradedK0}.

\[ \begin{array}{|c|c|c|c|c|c|c|} \hline 

1-PS & forms & \multicolumn{3}{|c|}{dim(({\cal K}_0)_i)}& \ell f & dim({\cal K}_{\ell f})  \\
 & & \multicolumn{3}{|c|}{dim(({\cal H})_i)}&  & \\ \hline
 & (f,g) & 1 & 0 & -1 && \\ 
 \hline \hline 
\ell_1 & det_3 & 0 & 8 & 8 & & - \\
& Q_1   & 0 & 8 & 8 & Q'_1 & 0+4+0 \\  \hline 
\ell_2 & det_3 & 0 & 8 & 8 & & - \\
& Q_2 & 0 & 8 & 8 & Q_3 & 0+8+0 \\  \hline 
\ell_4 & det_3 & 1 & 10 & 5 & & - \\
 & Q_4 & 1 & 13 & 7 & Q'_4 & 1+6+1 \\  \hline 
\end{array} \]

\section{Local Stabilizers and Lie algebra Cohomology obstructions }
\label{sec:cohomo}
We begin this section with two examples. Both concern the action of a 1-PS on a form $f$ with leading term $g$ and direction of approach $f_b$. We continue to use the notation ${\cal K}(t)$, ${\cal K}_0$ from the previous sections. We give an explicit description of the limiting algebra in both these examples. This motivates the question of recovering ${\cal K}(t)$ from the ${\cal K}_0$.

\begin{ex} \label{example:O2}
Let $X=\{ z,y\}$ and let $f=(y^2 +z^2 )^2$. Let $g=y^4$ and note that with $\lambda (t)\cdot y=y$ and $\lambda (t)\cdot z =t\cdot z$, gives us:
\[ f(t)=\lambda(t) \cdot f =y^4 +2t^2 \cdot y^2 z^2 +t^4 z^4 \]
Thus, $\lim_{t\rightarrow 0} \lambda (t) f(t) =g$, and that $a=0$ and $b=2$ in the notation of this section.

Let us order $X$ as $x_1 =z$ and $x_2=y$, and let $e_{11}, e_{12},e_{21},e_{22}$ be basis vectors of $gl(X)=gl_2$. In other words, we have: 
\[ e_{11}=\left[ \begin{array}{cc}
1 & 0 \\ 0 & 0 \\ \end{array} \right] \: \: 
e_{12}=\left[ \begin{array}{cc}
0 & 1 \\ 0 & 0 \\ \end{array} \right] \: \: 
e_{21}=\left[ \begin{array}{cc}
0 & 0 \\ 1 & 0 \\ \end{array} \right] \: \: 
e_{22}=\left[ \begin{array}{cc}
0 & 0 \\ 0 & 1 \\ \end{array} \right] \: \: \]
The stabilizer ${\cal H}$ of $g$ consists of matrices $\C \cdot e_{11}+\C \cdot e_{12}$, as show below:
\[ {\cal H}= \left[ \begin{array}{cc} 
* & * \\
0 & 0 \\ \end{array} \right] \]

We may choose ${\cal S}=\C \cdot e_{21}+\C \cdot e_{22}$ and note that $e_{21}\cdot g= 4y^3 z $ and $e_{22} \cdot g= 4g$. Thus $TO_g $ consists of the vector $\C \cdot y^4  +\C \cdot y^3 z$. We may decompose $V=TV_g$ as:
\[ V= (\C \cdot y^4 +\C \cdot y^3 z )\oplus (\C \cdot y^2 z^2 + \C \cdot yz^3 + \C \cdot z^4 )=TO_g \oplus N. \]
We may solve the equation from Proposition~\ref{prop:main}(ii) as follows. As per the recipe of that equation, we to find $\mf{h}$, $\mf{s}$ which when applied to the left and right of the expressions below within brackets, gives $n=0$.  
\eat{\[ 
\mf{h} f^+ (t)= \mf{s} (g+f^+(t)) + n  \] 
as:}
We compute the terms of \[ 
\mf{h} (2t^2 y^2 z^2 + t^4 z^4)= \mf{s} \cdot (y^4 +2t^2 y^2 z^2 + t^4 z^4) + n,  \] 
in the following table:
\[ \begin{array}{|c|c|c|} \hline 
\mf{h} & \lambda_{\cal S} & \lambda_N \\ \hline \hline
e_{11} & 0 & 4t^2 y^2 z^2 +4t^4 z^4\\ \hline 
e_{12} & t^2 e_{21} & 0 \\ \hline  
\end{array} \]
The stabilizer condition that $\lambda_N =0$ gives us the stabilizer of $f(t)$ as the element $\mf{k}(t)=e_{12}-t^2 e_{21}$ as shown below:
\[ \mf{k}(t)= \left[ \begin{array}{cc}
0 & 1 \\
-t^2 & 0 \\ \end{array} \right] \]
One may check, that ${\cal K}(t)=\C \cdot \mf{k}(t)$, and that ${\cal K}_0 = \C \cdot e_{12}$. 
We see that $f_b =4y^2 z^2$ and verify that $e_{12} \cdot f_b =4y^3 z \in TO_g$, and thus ${\cal K}_0 \subseteq {\cal H}_b$, the stabilizer of $\overline{f_b}$ under the $\star$-action of ${\cal H}$. $\Box $
\end{ex}

\begin{ex} \label{example:O3}
Let $X=\{ z, y_1 ,y_2 \}$ in that order and let $f=(y_1^2 +y_2^2 +z^2)^2$ and $g=(y_1^2 +y_2^2)^2$. Let $\lambda$ be such that $\lambda (t)\cdot y_i =y_i $ while $\lambda (t)\cdot z =tz$. Thus:
\[ f(t)= (y_1^2 +y_2^2 )^2 +2 t^2 z^2 (y_1^2 +y_2^2)+ t^4 z^4 \]
Let us have the following elements:
\[ {\cal K}(t)= \left\{ \left[ \begin{array}{ccc}
0 & a & b \\
-at^2 & 0 & c \\
-bt^2 & -c & 0 \\ \end{array} \right] | a,b,c\in \C \right\} \: \: 
{\cal H}= \left\{ = \left[ \begin{array}{ccc}
d & a & b \\
0 & 0 & c \\
0 & -c & 0 \\ \end{array} \right] |a,b,c,d \in \C \right\}
\]
Whence 
\[ {\cal K}_0 =\{ \left[ \begin{array}{ccc}
0 & a & b \\
0 & 0 & c \\
0 & -c & 0 \\ \end{array} \right] | \:  a,b,c\in \C \}
\]
We note that ${\cal K}_0 \subseteq {\cal H}$ and that ${\cal K}_0$ stabilizes $\overline{f_b} =\overline{(y_1^2 +y_2^2 )z^2}$ for the $\star$-action of ${\cal H}$. In other words, ${\cal K}_0 \cdot f_b \subseteq TO_g$. 

It is instructive to form the multiplication table for ${\cal K}(t)$. Let us put: 
\[ \mf{k}_1 (t)=\left[ \begin{array}{ccc}
0 & 1 & 0 \\
-t^2 & 0 & 0 \\
0 & 0 & 0 \\ \end{array} \right] 
\mf{k}_2 (t) = \left[ \begin{array}{ccc}
0 & 0 & 1 \\
0 & 0 & 0 \\
-t^2 & 0 & 0 \\ \end{array} \right] 
\mf{k}_3 (t)= \left[ \begin{array}{ccc}
0 & 0 & 0 \\
0 & 0 & 1 \\
0 & -1 & 0 \\ \end{array} \right] 
\]
We then have the following structure constants of the Lie bracket:
\[ \left[ \begin{array}{c} 
\mf{k}_1 (t) ,\mf{k}_2 (t)  \\
\mf{k}_1 (t),\mf{k}_3 (t) \\
\mf{k}_2 (t),\mf{k}_3 (t) \end{array} \right] =
\left[ \begin{array}{rrr} 
0 & 0 & -t^2\\
0 & 1 & 0\\
-1 & 0 & 0 \\
\end{array} \right] 
\left[ \begin{array}{c} 
\mf{k}_1 (t) \\
\mf{k}_2 (t) \\
\mf{k}_3 (t)\end{array} \right] 
\]
Here $\mf{k}_1 (t) ,\mf{k}_2 (t)$ on the left denotes 
$[\mf{k}_1 (t) ,\mf{k}_2 (t)]$ and likewise for the other rows on the left.   We note that $\{ \mf{k}_i (0)|i=1,2,3\}$ is a basis for ${\cal K}_0$ and the structure constants are precisely those of ${\cal K}(t)$ specialized to $t=0$.
\end{ex}

In the above example the structure constants $(\alpha_{ij}^k(t))$ of ${\cal K}(t)$ were elements of $\C [t]$ and where putting $t=0$ gave us ${\cal K}_0$. Moreover we see  that the $\C$-algebra ${\cal K} (t_0)$, i.e., ${\cal K}(t)$ instantiated at a $t_0 \in \C$, may be very different from ${\cal K}_0$, and perhaps ``more semisimple'' than both ${\cal H}$ or ${\cal K}_0$. We look at the process of constructing ${\cal K}(t)$ from ${\cal K}_0$. This naturally leads us to a Lie algebra cohomology formulation. We continue with the notation and transversality assumption of Section ~\ref{assume:traneverse}.

For any system $(\alpha_{ij}^k (t))$ of structure constants of ${\cal K}(t)$,  the verification that these do satisfy the Jacobi conditions is in the polynomial ring $\C [t]$. Let us consider the ring $\C[\epsilon_D]\cong \C[t]/(t^D)$ and Lie algebras over $\C[ \epsilon_D ]$. If $(\alpha_{ij}^k (t))$ satisfy the Jacobi identity, then so do $(\alpha_{ij}^k (\epsilon_D))$. Let us call the Lie algebra so generated as ${\cal K}_{D-1}$. Note that ${\cal K}_{D-1}$ is a Lie algebra, both over $\C$ and over $\C [\epsilon_D ]$ and that ${\cal K}_0$ matches with our earlier definition. Note that if $dim_{\C[t]} ({\cal K} (t))=k$, then each ${\cal K}_{D-1}$ is a free $\C [\epsilon_D]$-module of the same dimension. Moreover, for general $a>0$, we also have surjective Lie algebra homomorphisms (over $\C$) ${\cal K}_{a} \rightarrow {\cal K}_{a-1}$. 
Finally, for any ${\cal K}(t)$ over $\C[t]$, there is a large $D$, that the structure constants $(\alpha_{ij}^k (t))$ may be verified in $\C [\epsilon_D ]$. However, the corresponding ${\cal K}_{D-1}$ will not match ${\cal K}(t)$ in structure. With this remark, let us study the structure of ${\cal K}_1$. 

Let $\C [\epsilon]$ denote the ring formed by adding the variable $\epsilon $ with $\epsilon^2=0$.  For any $\C$-vector space or algebra $X$, let $X[\epsilon] $ denote the extension of coefficients, i.e., $X[\epsilon]=\C [\epsilon]\otimes_{\C} X$. 
Note that $X[\epsilon] =X \oplus \epsilon X$ as vector spaces. We also define the ``evaluation at zero'' map $e_0 :X[\epsilon] \rightarrow X$ as follows. For an element $x=x_0 +\epsilon x_1 \in X[\epsilon]$, we set $e_0 (x)=x_0$. This will also be denoted as $x(0)$. 
Let $V[\epsilon] $ and ${\cal G}[\epsilon]$ be obtained from $V$ and ${\cal G}$ as above. Note that ${\cal G}[\epsilon]$ is a Lie algebra and $V[\epsilon]$, a ${\cal G}[\epsilon]$ module. By the same token $V[\epsilon]$ is also an ${\cal H}[\epsilon]$ module under the $\star$-action.

\begin{lemma}
${\cal G}[\epsilon]$ is a Lie algebgra over $\C$ and the map $e_0 :{\cal G}[\epsilon] \rightarrow {\cal G}$ defined as $e_0 (\mf{g}+\epsilon \mf{g}')= \mf{g}$ is a Lie $\C$-algebra homomorphism.
\end{lemma}

Let us consider the element $p=g+\epsilon f_b \in V[\epsilon]$ and construct the stabilizer of the point $p$ in ${\cal G}[\epsilon]$. For a typical element $\mf{g}_0 +\epsilon \mf{g}_1$, we see that the stabilizer condition implies:
\begin{equation} \label{eqn:Hp}
\mf{g}_0 \cdot g=0 \: \: \mbox{  and  } \mf{g}_0 \cdot f_b + \mf{g}_1 \cdot g=0 
\end{equation}
The first condition states that $\mf{g}_0 \in {\cal H}$. The second condition further implies that $\mf{g}_0 \in {\cal H}_b \subseteq {\cal H}$, the stabilizer of $\overline{f_b }$ under the $\star $-action of ${\cal H}$ on $V$. This condition may also be interpreted as follows. Let us identify $TO_g $ with the quotient ${\cal G}/{\cal H}$, given simply by $d_g: \mf{g} \rightarrow \mf{g}\cdot g \in TO_g$. 

\begin{lemma}
\label{lemma:TasHmod}
The map $d_g :{\cal G}/{\cal H} \rightarrow TO_g \subseteq V$ is an isomorphism of ${\cal H}$-modules. 
\end{lemma}

\noindent 
{\bf Proof}: Note that ${\cal H}$ acts on $V$ as maps in $End(V)$ under the Lie bracket. Let $\mf{s}\in {\cal G}$ and let $\{ \mf{s} \}$ denote the coset $\mf{s}+{\cal H}$. The operation of $\mf{h}\in {\cal H}$ on ${\cal G}/{\cal H}$ is given by $\mf{h}\cdot \{ \mf{s}\}=\{ [\mf{h},\mf{s}] \}$. Whence, we then have:
\[ \begin{array}{rcl}
d_g (\{ [\mf{h},\mf{s}]\})&=& [\mf{h},\mf{s}]\cdot g \\
&=& (\mf{h}\cdot \mf{s} -\mf{s}\cdot \mf{h})\cdot g\\
&=& \mf{h}\cdot \mf{s} \cdot g \mbox{  (since $\mf{h}\cdot g=0$)} \\
&=& \mf{h} \cdot d_g (\mf{s}) 
\end{array} \]
This completes the proof. $\Box $

\begin{lemma} \label{lemma:Hp}
The stabilizer ${\cal H}_p \subseteq {\cal G}]\epsilon]$ is given by the set $\mf{h}_0 +\epsilon \mf{g}$ such that (i) $\mf{h}_0 \in {\cal H}_b$, and (ii) 
$\mf{g}\in d_g^{-1} (-\mf{h}_0 f_b )$. \end{lemma}

\noindent 
{\bf Proof}: Let us look at the two parts of Eq. \ref{eqn:Hp}. The first implies that $\mf{h}_0 \in {\cal H}$. The second implies that, in fact, $\mf{h}_0$ stabilize $\overline{f_b}$, and thus $\mf{h}_0 \in {\cal H}_b$. This proves (i). For (ii), note that the second condition of Eq. \ref{eqn:Hp}, gives us that $\mf{h}_0 \cdot f_b +d_g (\mf{g} )=0$. This implies (ii). $\Box $

Let us now define the map $d_b :{\cal H}_b \rightarrow {\cal G}/{\cal H}$ as follows. For any element $\mf{s}\in {\cal G}$, let $\{ \mf{s} \}$ denote the coset $\mf{s}+{\cal H}$ in ${\cal G}/{\cal H}$. For an $\mf{h}\in {\cal H}_b$, let $d_b (\mf{h})$ be that element $\{ \mf{s} \} \in {\cal G}/{\cal H}$ such that $\mf{s} \cdot g=\mf{h}\cdot f_b$, i.e., $d_g (\{\mf{s} \})=\mf{h}\cdot f_b$.  

\begin{defn} 
A derivation (see, for example, \cite{liecohomology1})~ $d$ from ${\cal G}$ to a ${\cal G}$-module W is a ${\mathbb C}$-linear map from ${\cal G}$ to $W$ such that for all $\mf{h}_1, \mf{h}_2 \in {\cal G}$ we have 
\[ d([\mf{h}_1,\mf{h}_2]) = \mf{h}_1 \circ d(\mf{h}_2) - \mf{h}_2 \circ d(\mf{h}_1).\]
\end{defn}

\begin{lemma}
(1) ${\cal G}/{\cal H}$ is an ${\cal H}_b$-module. The map $d_b :{\cal H}_b \rightarrow {\cal G}/{\cal H}$ is a derivation. (2) Moreover, the stabilizer ${\cal H}_p \subseteq {\cal G}]\epsilon]$ is also given by the set $\mf{h} +\epsilon \mf{g}$ such that (i) $\mf{h}_0 \in {\cal H}_b$, and (ii) 
$\mf{g}\in d_b (-\mf{h}_0 )$.
\end{lemma}

\noindent
{\bf Proof}: Since $[{\cal H}_b ,{\cal H}]\subseteq {\cal H}$, the quotient ${\cal G}/{\cal H}$ is an ${\cal H}_b$-module under the adjoint action. 
Let $\mf{h}_1 ,\mf{h}_2 \in {\cal H}_b$.  Suppose $\mf{h}_i f_b =\mf{s}_i g $ for some $\mf{s}_i \in {\cal G}$, whence $d_b (\mf{h}_i )=\{ \mf{s}_i \}$, the coset of $\mf{s}_i$. We must show that $d_b ([ \mf{h}_1 ,\mf{h}_2 ])= \mf{h}_1 \circ d_b(\mf{h}_2) - \mf{h}_2 \circ d_b(\mf{h}_1) $ as elements of
${\cal G}/{\cal H}$. We check this through the isomorphism $d_g $. Applying $d_g$ to the RHS, we get:
\[ \begin{array}{rcl} 
(\mf{h}_1 \circ d_b(\mf{h}_2) - \mf{h}_2 \circ d_b(\mf{h}_1) )\cdot g &=& (\{ [\mf{h}_1 ,\mf{s}_2]\}-\{ [\mf{h}_2 ,\mf{s}_1 ]\}) \cdot g \\
&=& \mf{h}_1  \mf{s}_2  g -\mf{h}_2  \mf{s}_1  g \mbox{   (since $\mf{h}_i \cdot g=0$)}
\end{array} \]
We check this with the LHS
\[ d_g ([\mf{h}_1 , \mf{h}_2 ])=[\mf{h}_1 ,\mf{h}_2 ] \cdot f_b = \mf{h}_1 \mf{h}_2 f_b -\mf{h}_2 \mf{h}_1 f_b = (\mf{h}_1 \mf{s}_2 -\mf{h}_2 \mf{s}_1 ) \cdot g\]
This proves (1). Concerning (2), the only change from lemma \ref{lemma:Hp} is in (ii). It is easily seen that $\mf{g}\in d_g^{-1} (-\mf{h}_0 f_b)$ is equivalent to the condition that $\mf{g}\in d_b (-\mf{h}_0 )$. $\Box $
\begin{remark}
The above claim is equivalent to the statement that ${\cal H}_p$ is a Lie subalgebra of ${\cal G}[\epsilon]$.
\end{remark}

We collect some properties of ${\cal H}_p$. 
\begin{lemma}
The image $e_0 ({\cal H}_p)$ is ${\cal H}_b$. Moreover, the nilpotent subalgebra $\epsilon {\cal H} \subseteq {\cal G}[\epsilon]$ is a subset of ${\cal H}_p$. 
\end{lemma}

\begin{defn} \label{epdeform}
Let ${\cal K}\subseteq \Hb$ be a subalgebra and let $dim_{\C} ({\cal K}) = K$. Let subalgebra $\overline{{\cal K}}\subseteq {\cal H}_p$ be such that (i) $\overline{\cal K}$ is generated as a $\C[\epsilon]$-module by $K$ elements $\mf{k}_1 , \ldots ,\mf{k}_K$, and (ii) $e_0 (\overline{\cal K})={\cal K}$. 
We call $\overline{\cal K}$ an $\epsilon$-deformation of ${\cal K}$. Note that $dim_{\C} (\overline{\cal K})=2K$. 
\end{defn}

The condition that $\overline{\cal K}$ is generated by the $K$ elements as a $\C[\epsilon]$-module implies that we can write the Lie bracket using the structure coefficients $(\alpha_{ij}^k )$, as:
\[ [\mf{k}_i , \mf{k}_j ]= \sum_k \alpha_{ij}^k \mf{k}_k \mbox { for some $(\alpha_{ij}^k ) \in \C[\epsilon]$} \]

\begin{lemma}
If $\overline{\cal K}$ is an $\epsilon$-deformation of ${\cal K}$ then $(\mf{b}_i) =(e_0 (\mf{k}_i )) $ is a $\C$-basis for ${\cal K}$ and its structure coefficients are precisely $(\alpha_{ij}^k (0))$, the evaluations of $(\alpha_{ij}^k )$ at $0$. Thus, we have:
\[ [\mf{b}_i , \mf{b}_j ]= \sum_k \alpha_{ij}^k (0) \mf{k}_k \]
\end{lemma}
\begin{prop} \label{prop:localmain}
The existence of $\overline{\cal K}$ is equivalent to the existence of a derivation $\overline{d_b}: {\cal K} \rightarrow {\cal G}/ {\cal K}$ which extends the derivation $d_b:{\cal K} \rightarrow {\cal G}/{\cal H}$ as shown in the diagram below:
\begin{center}
\begin{tikzcd}
{{\cal G}/{\cal K}} \arrow[r] 
& {{\cal G}/{\cal H}} \\
{\cal K} \arrow[u,"\overline{d_b}"] \arrow[r] & \arrow[u,"d_b"] \Hb
\end{tikzcd}
\end{center}
\end{prop}

\noindent 
{\bf Proof}: Suppose that $\overline{d_b}$ is such an extension. Let us define $\overline{\cal K}$ as follows:
\[ \overline{\cal K} =\{ \mf{h}+\epsilon \mf{s}| \mf{h}\in {\cal K}, \mf{s}\in {\cal G} \mbox{ such that } \mf{s} \in \overline{d_b}(-\mf{h}) \} \oplus \epsilon {\cal K}\]
Clearly, since $\overline{d_b}$ extends $d_b$, we have $\overline{d_b} (-\mf{h}) \subseteq d_b (-\mf{h})$ and hence every element of $\overline{\cal K}$ stabilizes $p$. Thus $\overline{\cal K}\subseteq {\cal H}_p$. Other properties of Defn. \ref{epdeform} are easily verified. Thus, all that remains to show is that $\overline{\cal K}$ is closed under the Lie bracket, i.e., it is a subalgebra of ${\cal H}_p$. 

Let us now verify the Lie bracket condition. For any elements $\mf{k}_1 =\mf{h}_1 +\epsilon \mf{s}_1 $ and $\mf{k}_2 = \mf{h}_2 +\epsilon \mf{s}_2 $ we must check that $[ \mf{k}_1 ,\mf{k}_2 ] \in \overline{\cal K}$. 
But
\[ [ \mf{k}_1 ,\mf{k}_2 ]= [\mf{h}_1 ,\mf{h}_2 ] +\epsilon ([\mf{h}_1 , \mf{s}_2 ]+[\mf{s}_1 ,\mf{h}_2 ] ) 
\] 
Since $[\mf{h}_1 ,\mf{h}_2 ] \in {\cal K}$, there is an element $[\mf{h}_1 ,\mf{h}_2 ]+\epsilon \mf{s} \in \overline{\cal K}$ with the condition that $\mf{s} \in \overline{d_b} (-[\mf{h}_1 ,\mf{h}_2 ])$. We must check that
\[ \mf{s}-([\mf{h}_1 \mf{s}_2 ]+[\mf{s}_1 ,\mf{h}_2 ]) \in {\cal K}. \]
Now, as elements of ${\cal G}/{\cal K}$, we have:
\[ \begin{array} {rcl}
[\mf{h}_1, \mf{s}_2 ]+[\mf{s}_1 ,\mf{h}_2 ] &=&  [\mf{h}_1, \mf{s}_2 ]-[\mf{h}_2 ,\mf{s}_1 ] \\
&=&
\mf{h}_1 \cdot \overline{d_b} (-\mf{h}_2) - \mf{h}_2 \cdot \overline{d_b }(-\mf{h}_1) \\
&=& \overline{d_b }(-[\mf{h}_1 ,\mf{h}_2 ]) \mbox{   since $\overline{d_b}$ is a derivation}  \\
&=& \mf{s} + {\cal K}\\
\end{array} \]
Thus $\mf{s}-[\mf{h}_1 ,\mf{s}_2 ]-[\mf{s}_1 ,\mf{h}_2 ] \in {\cal K}$. 
Whence:
\[ [\mf{h}_1 +\epsilon \mf{s}_1 ,\mf{h}_2 +\epsilon \mf{s}_2 ] = [\mf{h}_1 ,\mf{h}_2 ]+\epsilon \mf{s}+\epsilon \mf{h} \mbox{ with  } \mf{h} \in {\cal K} \]
This proves that $\overline{\cal K}$ in indeed a Lie subalgebra. Note that $\epsilon {\cal K} \subseteq \epsilon \overline{\cal K} \subseteq \overline{\cal K}$ and thus if $\mf{h}+\epsilon \mf{s} \in \overline{\cal K}$ then so is $\mf{h}+\epsilon \mf{s}'$, where $\mf{s}'\in \mf{s}+{\cal K}$.   

Conversely, if $\overline{\cal K}=\{ \mf{h}+\epsilon \mf{s} | \mf{h}\in {\cal K} \}$ is such an extension, then put $\overline{d_b} (\mf{h})=-\mf{s}+{\cal K}$. 
Since $\overline{\cal K}\subseteq {\cal H}_p$, we have $-\mf{s} + {\cal H}= d_b (\mf{h})$ as well. Thus $\overline{d_b}(\mf{h})\subseteq d_b (\mf{h})$ and $\overline{d_b}$ extends $d_b$. 

We must now show that $\overline{d_b }$ is a derivation. Towards this, let $\mf{h}_i +\epsilon \mf{s}_i $ and $[\mf{h}_1 ,\mf{h}_2 ]+\epsilon \mf{s}$ be elements of $\overline{\cal K}$. By the closure of $\overline{\cal K}$ under Lie bracket, we have: 
\[ [\mf{h}_1 ,\mf{s}_2]+[\mf{s}_1 ,\mf{h}_2] \in \mf{s}+{\cal K} \]
Let us verify the derivation condition:
\[ \begin{array}{rcl}
\mf{h}_1 \circ \overline{d_b }(\mf{h}_2 )-\mf{h}_2 \circ \overline{d_b }(\mf{h}_1 )     &=&  -[\mf{h}_1 ,\mf{s}_2 ] +[\mf{h}_2 , \mf{s}_1 ] +{\cal K}  \\
&=& -([\mf{h}_1 ,\mf{s}_2 ]+[\mf{s}_1 ,\mf{h}_2 ])+{\cal K} \\
&=& -\mf{s}+{\cal K} \\
&=& \overline{d_b}([\mf{h}_1 , \mf{h}_2 ])\end{array} \]
This proves the proposition. 
$\Box$

\begin{remark}
We have the exact sequence of ${\cal K}$-modules \cite{liecohomology1}:
\[ 0 \longrightarrow {\cal H}/{\cal K} \longrightarrow {\cal G}/{\cal K} \longrightarrow {\cal G}/{\cal H} \longrightarrow 0 \]
and the corresponding long exact sequence of cohomology modules:
\[ \begin{array}{rcl}
0 &\longrightarrow &H^0 ({\cal K},{\cal H}/{\cal K}) \longrightarrow H^0 ({\cal K},{\cal G}/{\cal K}) \longrightarrow H^0 ({\cal K}, {\cal G}/{\cal H}) \\
&\longrightarrow & 
H^1 ({\cal K},{\cal H}/{\cal K}) \longrightarrow H^1 ({\cal K},{\cal G}/{\cal K}) \longrightarrow H^1 ({\cal K}, {\cal G}/{\cal H}) \\
&\longrightarrow & 
H^2 ({\cal K},{\cal H}/{\cal K}) \longrightarrow \ldots
\end{array}
\]
Since both $d_b$ and $\overline{d_b}$ are derivations, they belong to the spaces $H^1 ({\cal K},{\cal G}/{\cal H})$ and $H^1 ({\cal K},{\cal G}/{\cal K})$ respectively (but they may be $0$). Whence, parts of the long exact sequence which are relevant to the extension of $d_b$ to $\overline{d_b}$ are:
\[ 
\longrightarrow  
H^1 ({\cal K},{\cal H}/{\cal K}) \longrightarrow H^1 ({\cal K},{\cal G}/{\cal K}) \longrightarrow H^1 ({\cal K}, {\cal G}/{\cal H}) \\
\longrightarrow  
H^2 ({\cal K},{\cal H}/{\cal K}) \longrightarrow
\]
\end{remark}
\begin{remark}
Our construction of the $\epsilon$-extension $\overline{\cal K}$ through a module leads to derivations as the direction of infinitesimal extension. This is a variation of the ideas of Nijenhuis and Richardson, \cite{liecohomology2}, where the 2-cocyles are the infinitesimal directions of deformations. All the same, it is likely that ${\cal K}_0$ is not {\em rigid} (again, \cite{liecohomology2}) while ${\cal K}(t_0 )$ are, for generic $t_0 \in \C$. 
\end{remark}

We now come to the main proposition of this section:

\begin{defn}
Let $A(t),f$ and $g$ be such that:
\[ A(t)\cdot f = f(t)=t^a g +t^b f_b + \mbox{ higher terms} \]
We say that $g$ is a regular limit of $f$ via $A(t)$ if (i) in the above expression $f_c \neq 0$ only when $c-a$ is a multiple of $b-a$, and (ii) if ${\cal H}$ is the stabilizer of $g$, then the stabilizer ${\cal K}(1)$ of $f$ is not contained inside ${\cal H}$. 
\end{defn}
\begin{prop}
Let $g$ be a regular limit of $f$ via $A(t)$ and $f(t)$ and $f_b$ be as above. Let ${\cal K}(t)$ be the stabilizing Lie algebra of $f(t)$ and ${\cal H}$ that of $g$. Moreover, let ${\cal K}_0 \subseteq {\cal H}_b$ be the limit of ${\cal K}(t)$, as $t\rightarrow 0$. Let $d_b : {\cal K}_0 \rightarrow {\cal G}/{\cal H}$ be the derivation as above. Then there is a derivation $\overline{d_b}:{\cal K}_0 \rightarrow {\cal G}/{\cal K}_0$ which extends $d_b$. 
\end{prop}

\noindent 
{\bf Proof}: Let $\delta =b-a$. Condition (i) of regularity implies that that after removal of $t^a$, all non-zero terms in the expression for $f(t)$ have powers of $t$ as $t^{k\delta}$ for some $k\geq 0$. 
We may absorb $t^a$ into $A(t)$ and write:
\[ A(t)\cdot f = f(t)=g +t^{\delta} f_b + \mbox{ higher terms in powers of $t^{\delta}$} \]
This, in turn, implies that the matrices $M_N (t)$ and $M_{\cal S} (t)$ have elements in $\C [t^{\delta}]$. Thus the basis for ${\cal K}(t)$ may be written in the form $\{ \mf{k}_i (t^{\delta}) \} _{i=1}^k $. Now condition (ii) that ${\cal K}\not \subseteq {\cal H}$ implies that ${\cal H} \cdot N \not \subseteq N$ and that there is a $\mf{k}\in {\cal K}$ with $\mf{k}=\mf{h}+\mf{s}$, with $\mf{h}\in {\cal H}$ and $\mf{s}\in {\cal S}$. This implies that $\lambda_{\cal S} (\mf{h}\cdot f_b )\neq 0$ and thus $\lambda_{\cal S}(\mf{h} \cdot f^+(t))$ must be of the form $t^{\delta} \mf{s}_1 + \mbox{higher terms}$. Thus some element $\mf{k}_i (t) = \mf{h}_0 +\mf{g}_1 t^{\delta} +\mbox{  higher terms}$ is such that $\mf{g}_1 \neq 0$. This implies that if $(\alpha_{ij}^k (T))$ are the structure constants of ${\cal K}(t)$, then there are some constants $i',j',k'$ such that $\alpha_{i'j'}^{k'} =c_1 t^{\delta} +c_2 t^{2 \delta} + \ldots$, where $c_1 \neq 0$. Whence, there is a non-trivial $\epsilon$-extension ${\cal K}_1$ of ${\cal K}_0$. By Prop. \ref{prop:localmain}, this is equivalent to the assertion that there is a derivation $\overline{d_b}:{\cal K}_0 \rightarrow {\cal G}/{\cal K}_0$ which extends $d_b$.

\begin{remark}
The Lie alegbra $\overline{{\cal K}_0}$ may not yield the ``additional semisimplcity" which may be present in ${\cal K}(t_0)$ (and absent in ${\cal K}_0$), but its structure coefficients do approximate those of ${\cal K}(t)$ better as polynomials in $t$. 
See Example \ref{example:O2}, where ${\cal K}_0$ is the $1$-dimensional algebra generated by $A$ below, while, $\overline{{\cal K}_0}$ is $2$-dimensional (over $\C$) and generated by $A(\epsilon)$.
On the other hand ${\cal K}(t_0)$ is $1$-dimensional and generated by $A(t_0 )$. 
\[ A= \left[ \begin{array}{cc}
0 & 1 \\
0 & 0 \end{array} \right] \mbox{   } 
 A(\epsilon)= \left[ \begin{array}{cc}
0 & 1 \\
-\epsilon & 0 \end{array} \right] \mbox{   }
 A(t_0 )= \left[ \begin{array}{cc}
0 & 1 \\
-t_0^2 & 0 \end{array} \right] \] 
\end{remark}

\begin{remark}
The existence of the limit $\lim_{t\rightarrow 0} f(t)=g$ allows the construction of ${\cal K}(t)$, the limit ${\cal K}_0$ and the tangent of approach $f_b$. We also have the additional $\star $-action of ${\cal H}$ and that ${\cal K}_0 \subseteq {\cal H}_b $, the stabilizer of $\overline{f_b}$ under this $\star$-action. Thus, we have the picture below: 
\[ \begin{array}{ccc}
{\cal K}(t) & \stackrel{lim}{\rightarrow}     & {\cal K}_0  \\
\downarrow &     & \downarrow \\
? & \stackrel{lim}{\rightarrow} & {\cal H}_b \\
\end{array}
\]
In other words, we already know that there is a local ${\cal G}$-action in the vicinity of the orbit $O(g)$ of $g$. Does this action allow us to compute the stabilizer ${\cal G}_{x,n}$ ( possibly over $\C(t)$), of the tangent vector $n\in TV_g$? This may then be applied to the vector $f_b$ to obtain a possible localization of ${\cal K}(t)$. The current recipe of $\overline{{\cal K}_0}$ over $\C [\epsilon]$ is not entirely satisfactory. 
 \end{remark}
 
 \begin{remark}{\bf The codimension $1$ case}
 \label{subsec:codim1}
This case yields two important simplifications. Firstly, the equality of ${\cal H}_b$ and ${\cal K}_0$ simplifies the extension problem as below:
\begin{center}
\begin{tikzcd}
{{\cal G}/\Hb} \arrow[r] 
& {{\cal G}/{\cal H}} \\
 & \arrow[ul,"\overline{d_b}"] \arrow[u,"d_b"] \Hb 
\end{tikzcd}
\end{center}

The second simplification comes from the theorem of Hoffman\cite{hofmann1965lie}, stated below:
\begin{theorem} \label{theorem:hoffman}
Let ${\cal K}\subseteq {\cal H}$ be Lie alegbras over $\R$ such that ${\cal K}$ is of codimension $1$. Then exactly one of the following three cases must be true:
\begin{enumerate}
    \item Let $P$ be the $2$-dimensional parabolic subgroup of upper triangular matrices in $sl_2$. Then there is an ideal ${\cal I}$ such that ${\cal K}/{\cal I}\cong P$ and ${\cal H}/{\cal I} \cong sl_2$. 
    \item Let $D\subseteq P$ the subalgebra of diagonal matrices in $P$. Then there is an ideal ${\cal I}$ such that ${\cal K}/{\cal I}\cong D$ and ${\cal H}/{\cal I} \cong P$
    \item ${\cal K}$ is an ideal of ${\cal H}$. 
    \end{enumerate}
\end{theorem}
Assuming that our Lie algebras are complex extensions of real Lie algebras, we see that ${\cal H}_b\subset {\cal H}$ must fall into one of these. Whence, the $1$-dimensional $\Hb$-module ${\cal H}/\Hb$ has a simple structure in the exact sequence below:
\[ \longrightarrow  
H^1 ({\cal H}_b,{\cal H}/{\cal H}_b) \longrightarrow H^1 ({\cal H}_b,{\cal G}/{\cal H}_b) \longrightarrow H^1 ({\cal H}_b, {\cal G}/{\cal H}) 
\longrightarrow  
H^2 ({\cal H}_b,{\cal H}/{\cal H}_b) \longrightarrow
\]
\end{remark}

\section{The single matrix under conjugation}
\label{sec:conj}
In this section we use the local model to solve a different problem - computing the projective limits of stable and polystable points of the vector space of $n\times n$ matrices under conjugation. This problem   has been  studied extensively in the (affine)  setting of semisimple Lie groups under the adjoint action. The proof of the main theorem in this section, Theorem~\ref{thm:conjugation:main}, illustrates the use of the local parametrization and the map $\theta$, to obtain projective closures. The statement of the main theorem was inspired by calculations performed using the local model of the neighbourhood of a matrix with a single Jordon block resulting in 
Theorem~\ref{theorem:Jn}. 

\begin{remark}
Experts in the area inform us that the statement of Theorem~\ref{thm:conjugation:main} does not surprise them and that it probably follows from earlier work. However we are not aware of any literature where this question has been discussed.
\end{remark}

Before proceeding to describe our results we set up some notation and recall some well known results about closures of orbits in affine space, in type A.  

Let $V$ be the $GL(X)$-module $End(X)$ where $GL(X)$ acts by conjugation. Thus, given $A\in GL(X)$ and $Y\in V$, we have $A\cdot Y=AYA^{-1}$. 
Polystable points in $V$ are the diagonalizable matrices,  and a generic stable point is a diagonalizable matrix with distinct eigenvalues, also called a regular semisimple element.

It is well known that the null-cone ${\cal N}$ for this action is the set of all nilpotent matrices,see \cite{newstead1978introduction}. The closure of nilpotent matrices in affine space is well understood, see \cite{gerstenhaber1961dominance},\cite{hesselink1976singularities}. The main theorem of this section has a description similar to that of closures of nilpotent matrices in affine space, which we review in Remark ~\ref{nilp-closure}.  

If $N\in V$ is a nilpotent matrix, in its $G$-orbit is the Jordan canonical form of $N$. Let us quickly recall this. The $b\times b$-matrix $J_b $ is given by the conditions: $J_b (i,j)=0$ for all $(i,j)$ except for the tuple $(i,i+1)$, for $i=1,\ldots,b-1$, where it is $1$. $J_4 $ is shown below:

\[ J_4 =\left[ \begin{array}{cccc}
0& 1 & 0 & 0 \\
0 & 0 & 1 & 0 \\
0 & 0& 0 & 1 \\
0 & 0 & 0 & 0\end{array} \right] \]
$J_b $ is determined by the properties that (i) it is a $b\times b$ matrix, (ii) $J_b^a \neq 0$ for $a<b$, and (iii) $J_b^b =0$. Given $N$, any $n\times n$ nilpotent matrix, its Jordan canonical form is given as the block diagonal matrix with blocks $J_{\lambda_1},\ldots, J_{\lambda_k }$, with $\lambda_1 \geq \ldots \geq \lambda_k >0$ and $\sum_i \lambda_i =n$, or in other words, a partition $\lambda $ of $n$. This partition $\lambda $ will be called the nilpotent signature of $N$.

Let us assume now that $X$ is an $n$-dimensional complex vector space and that we have chosen a basis of $X$. With respect to this basis every element in $V$ can be identified with an $n \times n$ matrix. 

Let $J_n $ be the nilpotent matrix given by $J_n (i,i+1)=1$ for $i=1,\ldots,n-1$, and zero otherwise. This is the single largest orbit in ${\cal N}$, and ${\cal N}$ is the closure of the orbit of $J_n $. 

In the next section we study the neighbourhood of $J_n$ in a direction normal to the orbit of $J_n$ under conjugation. We use the local model to show that a polystable point whose projective orbit closure contains $J_n$ must be necessarily stable.

In Section~\ref{jab} we consider the case $J_{a,b}$ and describe the semisimple and nilpotent matrices in its neighborhood. 

In Section \ref{ssfinal}, motivated by the above two computations, we state our main theorem which describes the nilpotent matrices which are in the projective closure of polystable points. This is described by the combinatorics of multiplicities of the eigenvalues of the matrix.

To simplify notation we don't use fraktur symbols to denote Lie algebra elements in the next two sections. Instead we use letters $A, B,C \ldots$ to denote matrices.

\subsection{The neighbourhood of $J_n$}

Recall that 
\[ J_n =\left[ \begin{array}{ccccc}
0& 1 & 0 & \ldots & \\
0 & 0 & 1 & 0 & \ldots \\
 & \vdots & & \vdots & \\ 
0 & 0& \ldots & 0 & 1 \\
0 & 0 & \ldots &  0 & 0\end{array} \right]. \]

The action $\rho $ of the Lie algebra $gl(X)$ on $V$ is given by $\rho (A)(Y)=A\cdot Y-Y\cdot A$, where $A\in gl(X)$ and $Y\in V$. For $k=-(n-1),\ldots,-1,0,1,\ldots, n-1$, let $Z_k $ be the $n\times n$-matrix where $Z_k (i,j)=0$ unless $j=i+k$, and then $Z_{i,i+k}=1$. Thus $J_n =Z_1$. We see that $Z_i Z_j =Z_{i+j}$ if $i,j\leq 0$ or $i,j\geq 0$. 

\begin{lemma}
\label{orbitJn}
Let $J_n $ be as above. Then:
\begin{enumerate} 
\item The tangent space of the orbit $O$ of $p=J_n $ under the $GL(X)$ action is
\[ TO_p =\{ B \in End(X) |\sum_{i=-k+1}^n b_{i,i+k}= 0  \mbox{ for $-(n-1)\leq k \leq 0 $} \} \]
\item A complement to $TO_p$ within $TV_p =End(X)$ is the space of matrices ${\cal C}_n$ is
\[ {\cal C}_n (c) =\left[ \begin{array}{lcccc}
-c_{n-1} & 0 & 0 & \ldots & \\
-c_{n-2} & 0 & 0 & 0 & \ldots \\
\vdots  & \vdots & & \vdots & \\ 
-c_1 & 0& \ldots & 0 & 0 \\
-c_0 & 0 & \ldots &  0 & 0\end{array} \right] \]
where $c=[c_0,\ldots,c_{n-1}]^T \in \C^n$ is a column vector. Note that the affine space $J_n +{\cal C}_n$ is the space of all companion matrices:
\[ J_n +{\cal C}_n =\left[ \begin{array}{lcccc}
-c_{n-1} & 1 & 0 & \ldots & \\
-c_{n-2} & 0 & 1 & 0 & \ldots \\
\vdots  & \vdots & & \vdots & \\ 
-c_1 & 0& \ldots & 0 & 1 \\
-c_0 & 0 & \ldots &  0 & 0\end{array} \right] \]
where $c_0,\ldots,c_{n-1} \in \C$. 
\item The stabilizer Lie algebra ${\cal H}$ of $J_n$ is the collection of all matrices $\sum_{i=0}^{n-1} a_i Z_i$, with $a_i \in \C$
\item A complement ${\cal S}$ to ${\cal H}$ within $gl(X)$ is the collection of all matrices $B\in gl(X)$ such that $B_{1,i}=0$ for $i=1,\ldots ,n$ and for every element $v\in TO_p$, there is a unique element $a\in {\cal S}$ such that $\rho(a)(J_n )=v$. 
\end{enumerate}
\end{lemma}

\noindent
{\bf Proof}: The action  of the Lie algebra $gl(X)$ on $p$ is given by $A\rightarrow AJ_n -J_n A=B$. But $B=(b_{ij})=(a_{i,j-1}-a_{i+1,j})$. We thus see that 
\[ \begin{array}{lrcll}
(I)& \sum_{i=1}^{n-k} b_{i,i+k}&=& a_{1,k}-a_{n-k+1,n} & \mbox{for $1\leq k\leq n-1$}\\
(II) & \sum_{i=-k+1}^n b_{i,i+k}&=& 0 & \mbox{for $-(n-1)\leq k \leq 0 $} \\
\end{array} \]
In fact, it is easily shown that for any matrix $B$ such that the $n$ equations constituting $(II)$ hold, we may find an $A\in gl(X)$ such that $B=AJ_n -J_n A$. This proves the first assertion. The second assertion follows from the fact that ${\cal C}_n \cap TO_p =\overline{0}$, and that $dim({\cal C}_n )=n$. Assertion (3) is classical. Assertion (4) follows from (1) and (2). $\Box $ 

\begin{remark}\label{rem:companionminpoly}
The minimal polynomial of $J_n+{\cal C}_n (c)$ is $p(X)=X^n +c_{n-1} X^{n-1}+\ldots + c_0 X^0 $.
\end{remark}

\noindent
{\bf Proof}: This is classical. Let us use $T$ for the matrix ${\cal C}_n (c)$. If $e_n=[0,\ldots,0,1]^T \in \C^n $ is the $n$-th column vector, then it is easily seen that $T^i \cdot e_n =e_{n-i}$ and that 
\[ \begin{array}{rcl}
T^n \cdot e_n &=&  -c_{n-1}e_1 -c_{n-2} e_2 +\ldots -c_0 e_n \\
&=& -c_{n-1}T^{n-1} \cdot e_n -c_{n-2} T^{n-2} \cdot e_n +\ldots -c_0 T^0 \cdot  e_n \\
\end{array} \]
Thus, we have $(T^n +c_{n-1}T^{n-1} + \ldots+ c_0 T^0 )\cdot e_n =0$. This, and the fact that $\{ T^i \cdot e_n $ are linearly independent vectors show that the polynomial $p$ is the minimial polynomial of $T$. $\Box $

\eat{Let us now compute that map $\theta (n)$. Since $\theta(n)$ is zero on $N \subset V$, to understand the higher compositions of $\theta(n)$ we need to study the maps $\lambda_{\cal S} (\cdot n) : {\cal S}\rightarrow {\cal S}$ given by $\theta (n) (s)=\lambda_{\cal S} (\cdot n)$, and the associated map $\lambda_{N} (s\cdot n)$. 
}
We now prove the following lemma.
\begin{lemma}
$\theta^i (n)$ is identically zero for $i \geq 2$.
\end{lemma}

\noindent
{\bf Proof}:
Since $\theta(n)$ is zero on $N \subset V$, it suffices to calculate it on the orbit which we have identified with ${\cal S}$. Given a matrix $B\in {\cal S}$ such that $B(1,i)=0$ for all $i$ and an element $C\in {\cal C}_n $ shown below, we see that $BC-CB$ is precisely $BC$ since $CB=0$. 
\[ C =\left[ \begin{array}{lcccc}
-c_{n-1} & 0 & 0 & \ldots & \\
-c_{n-2} & 0 & 0 & 0 & \ldots \\
\vdots  & \vdots & & \vdots & \\ 
-c_1 & 0& \ldots & 0 & 0 \\
-c_0 & 0 & \ldots &  0 & 0\end{array} \right] \]

We may thus write $B=[0,b_2^T ,\ldots, b_n^T ]^T$, where $b_i $ are row matrices and $C=[c,0,\ldots,0]$, where $c$ is a column matrix, to get
\[ BC-CB= \left[ \begin{array}{cccc}
0 & 0 &\ldots & 0 \\
b_2 \cdot c & 0 & \ldots & 0\\
\vdots & & \vdots & \\
b_n \cdot c & 0 & \ldots & 0 \\
\end{array} \right] \]
Whence $BC-CB \in N$ \eat{and thus, in the notation of the earlier section, $V\stackrel{\lambda_{\cal S}}{\longrightarrow} {\cal S} \stackrel{\cdot n}{\longrightarrow} V \stackrel{\lambda_N }{\longrightarrow} V=0$. {\color{red}  Should be ${\cal S} \stackrel{\cdot n}{\longrightarrow} V \stackrel{\lambda_{\cal S} }{\longrightarrow} V=0$}. 
{\color{blue}  Should be ${\cal S} \stackrel{\cdot n}{\longrightarrow} V \stackrel{\lambda_S }{\longrightarrow} V=0$}.}
Since $\theta^i(n)$, $i \geq 2$, is the composition of projection of $\theta^{i-1}(n)$ onto ${\cal S}$ followed by $\theta$, it follows that $\theta^i (n)$ is zero for every $i > 1$.$\Box$  

To carry out the program of the local model outlined in Theorem~\ref{thm:main} we need to understand how ${\cal H}$ acts on the space $N$, as given in the sequence given below. Note from the description of
${\cal H}$ in Lemma~\ref{orbitJn}, ${\cal R}$, the reductive part of ${\cal H}$ is zero, so in fact 
${\cal H} = {\cal Q}$.

\eat{\[   {\cal H} \stackrel{\cdot n}{\longrightarrow}V\stackrel{\lambda_{\cal S}}{\longrightarrow} {\cal S} \stackrel{\cdot n}{\longrightarrow} V \stackrel{\lambda_N }{\longrightarrow} V=0?\]}

We show 
\[   {\cal H} \stackrel{\cdot n}{\longrightarrow}V\stackrel{\lambda_{\cal S}}{\longrightarrow} {\cal S} \stackrel{\cdot n}{\longrightarrow} 0\]

For $A\in {\cal H}$ and $C\in N$ we have the first arrow:
\[ CA-AC=\left[ \begin{array}{cccc}
c_{n-1} a_0 & c_{n-1} a_1 & \ldots &  c_{n-1}a_{n-1} \\
c_{n-2} a_0 & c_{n-2} a_1 & \ldots &  c_{n-2}a_{n-1} \\
\vdots & & \vdots &  \\
c_{0} a_0 & c_{0} a_1 & \ldots & c_{0}a_{n-1} \\
\end{array} \right] -
\left[ \begin{array}{cccc}
\sum_{i=0}^{n-1} c_{n-1-i} a_i & 0 & \ldots &   0 \\
\sum_{i=0}^{n-2} c_{n-2-i} a_i & 0 & \ldots & 0 \\
\vdots & & \vdots &  \\
 c_0 a_0 & 0 & \ldots & 0\\
\end{array} \right] \]
Since the expression is linear for ${\cal H}$, let us assume that $a_i=1$ for a particular $i$, and $a_j=0$ for $j\neq i$. Whence, we have: 
\[CA_i -A_i C=
\left[ \begin{array}{ccccccccc}
-c_{n-i} & 0 & \ldots & 0 \ldots & & 0 &0 & 0 & 0\\
-c_{n-i-1} & 0 & \ldots & 0 \ldots & & 0 & 0 & 0 & 0\\
\vdots &\vdots &  &\vdots \ldots & &\vdots & 0\\
-c_0 &\vdots & \vdots & & &\vdots & 0 & 0 & 0\\
\vdots & & \vdots & & & c_{n-i}& 0 & 0 & 0\\
0 & 0 & \ldots& 0 \ldots  & & c_{n-i-1}& 0 & 0 & 0\\
\vdots & & \vdots &  & &\vdots & 0 & 0 & 0\\
0 & 0 & \ldots & & & c_0& 0 & 0 & 0\\
\end{array} \right] \]

The sum of the entries in the diagonals $j, -(n-1) \leq j \leq 0$ are all zero, so from Lemma~\ref{orbitJn} this is in the orbit $TO_{J_n}$.

Thus there is an $S\in {\cal S}$ such that $S\cdot J_n -J_n \cdot S=O_i $. Since $S\cdot C-C\cdot S=0$ we have
\[   {\cal Q} \stackrel{\cdot n}{\longrightarrow}V\stackrel{\lambda_{\cal S}}{\longrightarrow} {\cal S} \stackrel{\cdot n}{\longrightarrow} V =0\]

\eat{We illustrate the proof through an example when 
$n=5$ and $j=3$.
\[ O_3=\left[ \begin{array}{ccccc}
-c_2 & 0 & c_4 & 0 & 0 \\
-c_1 & 0 & c_3 & 0 & 0 \\
-c_0 & 0 & c_2 & 0 & 0 \\
0 & 0 & c_1 & 0 & 0 \\
0 & 0 & c_0 & 0 & 0 \\
\end{array} \right] \]
}
\eat{For this example $S_3 $ is as given below, along with the matrix $C_3 $:
\[ S_3=\left[ \begin{array}{ccccc}
 0 & 0 & 0 & 0 & 0 \\
-c_2 & 0 & c_4 & 0 & 0 \\
-c_1& -c_2& c_3 & c_4 & 0 \\
-c_0 & -c_1 & 0 & c_3 & c_4 \\
0 & -c_0 & 0 & 0 & c_3\\
\end{array} \right] 
C_3=\left[ \begin{array}{ccccc}
 c_4 & 0 & 0 & 0 & 0 \\
c_3 & 0 & 0 & 0 & 0 \\
c_2 & 0& 0 & 0 & 0 \\
c_1 & 0 & 0 & 0 & 0 \\
c_0  & 0 & 0 & 0 & 0\\
\end{array} \right] 
\]

}
$\Box $

Our final computation is the dimension of the stabilizer of an arbitrary point sitting normally over $J_n $. 
\begin{lemma}
The dimension of the stabilizer of any point $J_n +{\cal C}_n (c)$ is $n$.
\end{lemma}

\noindent
{\bf Proof}: Let $T={\cal C}_n (c)$. By Theorem \ref{thm:main}, the stabilizer of the point $J_n +T$ is governed by the equation:
\[ \Delta r \cdot T 
+\lambda_N \circ (1+\theta (n))^{-1} (\Delta q \cdot T )=0 \]
where ${\mf r} \in {\cal R}, \mf{q} \in {\cal Q}$, the reductive and nilpotent Lie subalgebras of ${\cal H}$. As mentioned earlier ${\cal R} = 0$.
\eat{Whence ${\cal R}=0$ and the generic element of the stabilizer ${\cal H}$ of $J_n$ is given by $A$. Evaluating $TA-AT$, we get:
\[ TA-AT=\left[ \begin{array}{cccc}
c_{n-1} a_0 & c_{n-1} a_1 & \ldots &  c_{n-1}a_{n-1} \\
c_{n-2} a_0 & c_{n-2} a_1 & \ldots &  c_{n-2}a_{n-1} \\
\vdots & & \vdots &  \\
c_{0} a_0 & c_{0} a_1 & \ldots & c_{0}a_{n-1} \\
\end{array} \right] -
\left[ \begin{array}{cccc}
\sum_{i=0}^{n-1} c_{n-1-i} a_i & 0 & \ldots &   0 \\
\sum_{i=0}^{n-2} c_{n-2-i} a_i & 0 & \ldots & 0 \\
\vdots & & \vdots &  \\
 c_0 a_0 & 0 & \ldots & 0\\
\end{array} \right]
\]}
From the calculations in the preceeding paragraph, the image of ${\cal Q} \cdot n$ is in ${\cal S}$ and so $\lambda_N \circ (1+\theta (n))^{-1} (\Delta q \cdot T )=0$ for all ${\mf q}$. So every elements of ${\cal Q}$ may be supplemented by an element of ${\cal S}$ to given an element of the stabilier of $T + J_n$. $\Box$

\begin{remark}
That the matrix $T+J_n$ is in companion form actually tells us that its minimal polynomial is of degree $n$
\end{remark}

We now come to the main theorem:
\begin{theorem} \label{theorem:Jn}
Let $A \in End(X)$ be polystable under the action of $SL(X)$ by conjugation. If $J_n $ is a projective limit of a matrix $A$ under the adjoint action, then $A$ has $n$ distinct eigenvalues and is stable. Moreover, for any matrix $A$ such that $A$ has distinct eigenvalues $\lambda_1 ,\ldots, \lambda_n $, there is a family ${\cal A}(t)$ such that (i) the eigenvalues of ${\cal A}(t)$ are $t\lambda_1 ,\ldots ,t\lambda_n $ and (ii) $\lim_{t\rightarrow 0} {\cal A}(t)=J_n $.  
\end{theorem}

\noindent
{\bf Proof}: Under the hypothesis, we know that there is a $GL(X)$-conjugate $A'$ such that $\| A'-J_n \|< \epsilon$. By the local model construction, we may further assume that $A'$ is of the form $J_n +T$, with $T\in N$. Since $A$ is polystable, so is $A'$ and so it is diagonalizable. Since $A'=J_n +T$ as above, its minimal polynomial is of degree $n$. Thus, $A'$ and therefore $A$ has $n$ distinct eigenvalues. This proves the first part.

For the second part, let ${\cal C}_n (c)$ be the companion form of $A$. Thus, the vector $c=(c_{n-1} ,\ldots, c_0)$ are the coefficients of the characteristic polynomial of $A$. Define $c(t)=(tc_{n-1}, t^2 c_{n-2}t^2,\ldots, t^n c_0)$ and let ${\cal A}(t)=J_n +{\cal C}(c(t))$, i.e., the companion matrix with column vector $c(t)$. Then ${\cal A}(t)$ has eigenvalues $(t\lambda_1,\ldots,t\lambda_n )$ and $\lim_{t\rightarrow0} {\cal A}(t)=J_n$. 
$\Box $

\begin{ex}
Let us compute explicitly through ${\cal S}$-completion, the stabilizer of the matrix $Z_4 = J_4 +C_4$ sitting "normally" above $J_4$:

\[ Z_4 =\left[ \begin{array}{cccc}
0 & 1 & 0 & 0 \\
0 & 0 & 1 & 0\\
0 & 0 & 0 & 1\\
1 & 0 & 0 & 0 \\ \end{array} \right] 
\: \: \: 
C_4 =\left[ \begin{array}{cccc}
0 & 0 & 0 & 0 \\
0 & 0 & 0 & 0\\
0 & 0 & 0 & 0\\
1 & 0 & 0 & 0 \\ \end{array} \right]
\]
The matrix $Z_4 $ has the minimal polynomial $X^4-1$. The stabilizing Lie algebra of $Z_4 $ is given by the basis elements $Z_4^i$, for $i=0,\ldots,3$. 

The complementary space ${\cal S}$ is the 12-dimensional space of matrices shown below:
\[ {\cal S} =\left[ \begin{array}{cccc}
0 & 0 & 0 & 0 \\
* & * & * & *\\
* & * & * & *\\
* & * & * & * \\ \end{array} \right] \] The elementary matrices $E_{ij}$, for $i=2,3,4$ and $j=1,2,3,4$ form a basis for ${\cal S}$. The tangent space to the orbit of $J_4 $ is given by the basis of 12 matrices $E_{ij}\cdot J_4 -J_4 \cdot E_{ij}$.  The stabilizer ${\cal H}$ of $J_4$ is the space of matrices generated by $J_4^i$, $i=0,\dots,3$. The stabilizer space of $Z_4$ is again $4$-dimensional and given by the elements $h+s$ where $h\in {\cal H}$ and $s\in {\cal S}$ such that $[s,J_4]+[h,C_4]=0$. Let us select the element $h=J_4^3$ to obtain $[h ,C_4 ]$ as shown below:
\[ h =\left[ \begin{array}{cccc}
0 & 0 & 0 & 1 \\
0 & 0 & 0 & 0\\
0 & 0 & 0 & 0\\
0 & 0 & 0 & 0 \\ \end{array} \right] 
\: \: \: 
[h ,C_4 ]= \left[ \begin{array}{cccc}
1 & 0 & 0 & 0 \\
0 & 0 & 0 & 0\\
0 & 0 & 0 & 0\\
0 & 0 & 0 & -1 \\ \end{array} \right]
\]
One may check that $[s ,J_4]$ below and verify that $[s,J_4 ]=-[h ,C_4]$:
\[ s =E_{21}+E_{32}+E_{43}=\left[ \begin{array}{cccc}
0 & 0 & 0 & 0 \\
1 & 0 & 0 & 0\\
0 & 1 & 0 & 0\\
0 & 0 & 1 & 0 \\ \end{array} \right] 
\: \: \: 
[s ,J_4 ]= \left[ \begin{array}{cccc}
-1 & 0 & 0 & 0 \\
0 & 0 & 0 & 0\\
0 & 0 & 0 & 0\\
0 & 0 & 0 & 1 \\ \end{array} \right]
\]
Thus $s+h$ shown below, stabilizes $Z_4$.
\[ s +h =\left[ \begin{array}{cccc}
0 & 0 & 0 & 1 \\
1 & 0 & 0 & 0\\
0 & 1 & 0 & 0\\
0 & 0 & 1 & 0 \\ \end{array} \right] \]
\end{ex}

\subsection{The $J_{a,b}$ case}\label{jab}

We now consider the partition $(a,b)$ where $a+b=n$ and $a\geq b$. The matrix $J_{a,b}$ may be written in the block-diagonal form as below:
\[ J_{a,b} =\left[ \begin{array}{cc}
J_a & 0 \\
0 & J_b \end{array} \right] \]

\begin{lemma} \label{lemma:abstab}
The stabilizer ${\cal H}$ of $J_{a,b}$ is of the form:
\[ \left[ \begin{array}{cc}
Z_a & Y_{ab} \\
Y_{ba} & Z_b \end{array} \right] \]
Where $Z_m$ is $m\times m$ such that (i) $Z_m(i,j)=0$ for all $i>j$, and (ii) $Z_m(i+1,j+1)=Z_m(i,j)$. The forms for $Y_{ab}$ and $Y_{ba}$ are given below:
\[ Y_{ab}= \left[ \begin{array}{c}
Z'_b \\
0 \end{array} \right] \: \: \: \: 
Y_{ba}=\left [ \begin{array}{cc}
0 & Z''_b \end{array} \right] \]
The dimension of ${\cal H}$ is $a+3b$, $a$ coming from $Z_a$ and $b$ each from the
three matrices $Z_b, Z'_b,Z''_b$.
\end{lemma} 

\begin{lemma} \label{abnormal}
A normal space to the orbit of $J_{a,b}$ is ${\cal C}$ below:
\[ {\cal C}=\left[ \begin{array}{cc}
{\cal C}_a & W_{ab} \\
W_{ba} & {\cal C}_b \end{array} \right] \]
where ${\cal C}_a$ and ${\cal C}_b$ are as shown below:
\[ {\cal C}_a (c) =\left[ \begin{array}{lcccc}
-c_{a-1} & 0 & 0 & \ldots & \\
-c_{a-2} & 0 & 0 & 0 & \ldots \\
\vdots  & \vdots & & \vdots & \\ 
-c_1 & 0& \ldots & 0 & 0 \\
-c_0 & 0 & \ldots &  0 & 0\end{array} \right] \: \: \:  
{\cal C}_b (d) =\left[ \begin{array}{lcccc}
-d_{b-1} & 0 & 0 & \ldots & \\
-d_{b-2} & 0 & 0 & 0 & \ldots \\
\vdots  & \vdots & & \vdots & \\ 
-d_1 & 0& \ldots & 0 & 0 \\
-d_0 & 0 & \ldots &  0 & 0\end{array} \right] \]
and $W_{ab}$ and $W_{ba}$ are of the following form:
\[ W_{ab}= \left[ \begin{array}{ccc}
&&\\
&0_{(a-1)\times b} &\\
&& \\ \hline
\alpha_1 & \ldots & \alpha_b \end{array} \right] \: \:\: 
W_{ba} = \left[ \begin{array}{c|c}
\beta_1 &  \\
\vdots & 0_{b\times{(a-1)}}  \\
\beta_b & \end{array} \right] \] 
where $c_i ,d_j,\alpha_r ,\beta_s \in \C$ and $O_{r\times s}$ are zero matrices. The dimension of ${\cal C}$ is also $a+3b$. 
\end{lemma}

\noindent 
The proofs of these lemmas are straightforward. We now come to two main lemmas.

\begin{lemma} \label{lemma:abminimial}
For any $C\in {\cal C}$, let $p_C$ be the minimal polynomial of $T=J_{a,b}+C$. Then the degree of $p_C$ is at least $a$. 

\end{lemma}

\noindent 
{\bf Proof}: We may write $T$ in the jumbo-size matrix below:
\[ \left [ 
\begin{array}{lcccc|lcccc}
-c_{a-1} & 1 & 0 & \ldots & & & \\
-c_{a-2} & 0 & 1 & 0 & \ldots  & 
& &&\\
\vdots  & \vdots & & \vdots & & & \multicolumn{3}{c}{0_{(a-1)\times b}} &\\ 
-c_1 & 0& \ldots & 0 & 1  & & && \\
-c_0 & 0 & \ldots &  0 & 0 & \alpha_1 & \alpha_2 &\ldots & &\alpha_b \\ \hline 
\beta_1 &&&&&-d_{b-1} & 1 & 0 & \ldots & \\
\beta_2 &&&&&-d_{b-2} & 0 & 1 & 0 & \ldots \\
\vdots &\multicolumn{3}{c}{0_{b\times(a-1)}} &&\vdots  & \vdots & & \vdots & \\ 
&&&&&-d_1 & 0& \ldots & 0 & 1 \\
\beta_{b} &&&&&-d_0 & 0 & \ldots &  0 & 0 
\end{array}
\right] \]
Let $e_1 ,\ldots ,e_a ,e_{a+1}, \ldots ,e_{a+b}$ be a basis of $\C^n$ (as column vectors). \eat{Let us also write $e_{a_j}$ as $f_j$, for convenience.} From the matrix structure of $T$, it is clear that for $i=2,\ldots,a$ we have $T(e_i )=e_{i-1}$. Thus if $Q=T^m +\gamma_{m{-1}} T^{m-1} +\ldots +\gamma_0 T^0$ is any operator, where $m<a$, then $Q(e_a )=e_{a_m}+\gamma_{m-1} e_{a-m-1}+\ldots + \gamma_0 e_a$. When $Q(e_a)\neq 0$. Thus, $T$ cannot have a minimal polynomial of degree less than $a$. $\Box$ 

\begin{lemma}
\label{lemma:abnull}
Let $T$ be as above and $\lambda \in C$, be any number. Let $ker(T-\lambda I)$ be the space of all $v$ such that $T(v)=\lambda v$, then $dim(ker(T-\lambda I))\leq 2$. Moreover, any $v$ in this space is determined by $v_1$ and $v_{a+1}$. 
\end{lemma}

\noindent 
{\bf Proof}:
The condition that $(T-\lambda I)v=0$ gives us the following conditions:
\[ \begin{array}{rl}
-c_{a-k+1}v_1 -\lambda v_{k-1} +v_k =0 & \mbox{for $k=2,\ldots,a$} \\
\beta_{i-1} v_1 -d_{b-i+1} v_{a+1}-\lambda v_{a+i-1}+v_{a+i}=0 & \mbox{for $i=2,\ldots,b$} \\
\end{array} \]
This proves the claim. $\Box$

For $T$ as above, we write \[ T_a =\left[ \begin{array}{lcccc}
-c_{a-1} & 1 & 0 & \ldots & \\
-c_{a-2} & 0 & 1 & 0 & \ldots \\
\vdots  & \vdots & & \vdots & \\ 
-c_1 & 0& \ldots & 0 & 1\\
-c_0 & 0 & \ldots &  0 & 0\end{array} \right] \: \: \:  
{T}_b =\left[ \begin{array}{lcccc}
-d_{b-1} & 1 & 0 & \ldots & \\
-d_{b-2} & 0 & 1 & 0 & \ldots \\
\vdots  & \vdots & & \vdots & \\ 
-d_1 & 0& \ldots & 0 & 1\\
-d_0 & 0 & \ldots &  0 & 0\end{array} \right] \]

Note that, by Remark~\ref{rem:companionminpoly}, the minimal polynomial of $T_a$ (respectively $T_b$) 
is $T^a + c_{a-1}T^{a-1} + \ldots + c_0 T^0$ (respectively $T^b + d_{b-1}T^{a-1} + \ldots + d_0T^0$).

\begin{lemma}\label{lemma:abminpolya}
Suppose that the minimal polynomial $p$ of $T$ has degree $a$. Then $p$ is the minimal polynomial of $T_a$ and
the minimal polynomial of $T_b$ divides $p$.
\end{lemma}

\noindent
{\bf Proof:} Let $p(T) = T^a + \gamma_{a{-1}} T^{a-1} +\ldots +\gamma_0 T^0$ be the minimal polynomial of $T$ of degree $a$.
Using the notation as in the proof of Lemma~\ref{lemma:abminimial}, $p(T) (e_a) = 0$ implies that
$\gamma_i = c_i$ for $i=0,\ldots, a{-1}$ and $\beta_j = 0$ for $j=1, \ldots, b$.  Thus, 
$p(T) = T^a + c_{a{-1}} T^{a-1} + \ldots + c_0T^0$ is also the
minimal polynomial of $T_a$.

We now write $T$ in block upper
triangular form
\[
T = \left[
\begin{array}{c|c}
T_a & * \\ \hline
0_{b\times a} & T_b
\end{array}
\right]
\]
It is easy to see that for any polynomial $Q$, $Q(T)$ is of the form
\[
Q(T) = \left[
\begin{array}{c|c}
Q(T_a) & * \\ \hline
0_{b\times a} & Q(T_b)
\end{array}
\right]
\]
As $p(T)=0$, this shows that $p(T_b)=0$. Therefore, the minimal polynomial of $T_b$ divides $p$.
$\Box$

\begin{prop} \label{prop:abnil}
The nilpotent matrices in the space $J_{a,b}+{\cal C}$ have signature $(a',b')$ with $a'\geq a$. Moreover, $J_{a,b}$ is in the projective orbit closure of $J_{a',b'}$. 
\end{prop}

\noindent 
{\bf Proof:} Let $J$ be a nilpotent matrix of the above form and let $a'$ be the largest component of its signature. The conditions on the minimum degree of the minimal polynomial implies that $a'\geq a$. The dimension of the kernel of $J$,  with $\lambda =0$ in lemma \ref{lemma:abnull}, is only $2$. That forces $J$ to have only two components.
Thus, $J=J_{a',b'}$ with $a'\geq a$. That every such $J_{a',b'}$ is obtained is verified by putting all $c$'s, $d$'s, $\beta$'s and $\alpha$'s to zero except that $\alpha_i =t$, where $t\neq 0$. One may verify that this matrix $J_i (t)$ is nilpotent with signature $a'=a+b-i+1$ and $b'=i-1$. We also see that $\lim_{t\rightarrow 0}J_i (t)=J_{a,b}$, and thus $J_{a,b}$ is in the projective orbit closure of $J_{a',b'}$.  $\Box$

\begin{prop} \label{prop:abdiag}
The diagonalizable matrices in the space $J_{a,b}+{\cal C}$ have at least $a$ distinct eigenvalues. Moreover,
each eigenvalue has multiplicity at most $2$. Moreover, for any such matrix $A$, the matrix $J_{a,b}$ is in its projective closure. 
\end{prop}

\noindent
{\bf Proof:}  For a diagonalizable matrix $T$, the number of distinct eigenvalues is equal to the degree of
the minimal polynomial. Further, for an eigenvalue $\lambda$ of $T$, the (algebraic) multiplicity of $\lambda$ is
equal to its geometric multiplicity, namely $dim(ker(T-\lambda I))$. The claim now follows from Lemmas~\ref{lemma:abminimial} and \ref{lemma:abnull}.

Coming to the second assertion, let $\{ \lambda_1 ,\ldots ,\lambda_{b'}, \lambda_{b'+1},\ldots, \lambda_{a'}\}$ be an ordering of the eigenvalues of $T$ such that $\lambda_1 ,\ldots,\lambda_{b'}$ have multiplicity $2$ and thus $a'+b'=n$. By extending the argument of Theorem \ref{theorem:Jn}, $J_{a',b'}$ is in the projective orbit closure of $T$. By prop. \ref{prop:abnil}, we have $a'\geq a$, and thus $J_{a,b}$ is in its projective orbit closure of $J_{a',b'}$, and therefore of $T$.  
$\Box$

\subsection{Orbit closure of polystable points}\label{ssfinal}
We end the section with a complete characterization of the orbit closure of polystable points in the vicinity of a single nilpotent matrix. The closure is completely determined by the combinatorics of the spectrum of the polystable point under consideration. We give an algebraic proof which is motivated by the constructions of the previous section. 

Before we begin, let us recall some definitions. A partition $\alpha$ of $n$ is a sequence $(\alpha_1 \geq \cdots \geq \alpha_r >0)$ such that $\sum_{i=1}^r \alpha_i =n$.

\begin{defn}
Let $\alpha $ and $\beta $ be partitions of $n$. The transpose $\gamma$ of $\alpha $ is given by the numbers $\gamma_1 \geq \ldots \geq \gamma_s >0$ where $\gamma_i =| \{ j | \alpha_j \geq i\}|$, and is denoted by $\alpha^T$. It is easy to check that $\gamma $ too is a partition of $n$. We say that $\alpha$ dominates $\beta$ iff for all $i>0$, $\sum_{j=1}^i \alpha_j \geq \sum_{j=1}^i \beta_j$. This is denoted by $\alpha \unrhd \beta$ or, equivalently, $\beta  \unlhd \alpha$. It is easy to see that
$\alpha \unrhd \beta$ iff $\beta^T \unrhd \alpha^T$.
\end{defn}

The order described above on partitions of $n$ is a partial order called the dominance order.

\begin{remark} \label{nilp-closure} Let $X$ be a nilpotent matrix with signature $\lambda$. The orbit closure of $X$ under conjugation contains all nilpotent matrices with partition signature $\lambda'$ with
$\lambda \unrhd \lambda'$. This is a celebrated result proved by Gerstenhaber, \cite{gerstenhaber1961dominance}. 
\end{remark}

If an $n\times n$-diagonalizable matrix has eigenvalues $\mu_1 ,\mu_2, \ldots , \mu_s $ with multiplicities $\lambda := \lambda_1 \geq \lambda_2 \geq \ldots \geq \lambda_s$,  we call $\lambda$, the spectrum partition of that matrix.

\begin{prop}
Let ${\cal X}_k^r$ denote the set of all $n\times n$ matrices $x$ for which there exist  $\lambda_1 , \ldots , \lambda_k$  such that  $\mathrm{rank}((x-\lambda_1 I)(x-\lambda_2 I)\ldots (x-\lambda_k I))\leq r$ and
$\mathrm{det}(x - \lambda_i I)=0$ for all $i=1, \ldots, k$. Then ${\cal X}_k^r$ is a projective variety. Moreover, ${\cal X}_k^r$ is invariant under the action of $GL_n$ on $x$ by conjugation.   
\end{prop}

\noindent
{\bf Proof}: Let $X$ be an $n\times n$ matrix of variables and let $Y=\{y_1 ,\ldots, y_k\}$ be another set of indeterminates.

Let the affine variety ${\cal V} \subseteq {\mathbb C}^{n^2+k}$ be
defined by the homogeneous ideal  $J$ in ${\mathbb C}[X,Y]$
generated by $\mathrm{det}(X - y_iI) = 0$ for $i=1,\ldots, k$. We define
the projection morphism $\phi: {\cal V} \to {\mathbb C}^{n^2}$ as follows:
$\phi(A,b_1, \ldots, b_k) = A$. 
The induced map of the co-ordinate rings
is simply the inclusion of ${\mathbb C}[X] \hookrightarrow {\mathbb C}[{\cal V}]={\mathbb C}[X, Y]/J$. As $\mathrm{det}(X-y_i I)=0$, each
$y_i \in {\mathbb C}[{\cal V}]$ is integral
over ${\mathbb C}[X]$. This shows that $\phi$ is
a finite map and hence maps closed sets to
closed sets.

Consider the matrix $Z=(X-y_1 I)(X-y_2 I)\ldots (X-y_k I))$. Let $P,Q$
be subsets of $[n]$ such that $|P|=|Q|=r+1$ and let $f_{P,Q}$ be the homogeneous form $\mathrm{det}(Z_{P,Q})$, i.e., the determinant of the $(P,Q)^{\mathrm{th}}$ minor of $Z$. Consider the ideal $I \subseteq {\mathbb C}[{\cal Z}]$ generated by all such forms. This is a homogeneous ideal in the variables of $X\cup Y$. Its variety, say ${\cal U} \subseteq {\cal V}$ is indeed of points  $(A, \lambda_1 ,\ldots , \lambda_k ) \in {\cal V}$ such that $\mathrm{rank}((A-\lambda_1 I)(A-\lambda_2 I)\ldots (A-\lambda_k I))\leq r$. 
The set ${\cal X}_k^r$ is precisely the
image of the closed set ${\cal U}$ under
the finite map $\phi$. As a result,
${\cal X}_k^r$ is also a closed set. The 
fact that ${\cal X}_k^r$ is closed under
homothety and invariant under the conjugation
action of $GL_n$ is obvious.
$\Box$ 

\begin{prop}\label{prop:diag}
Let $x$ be an $n\times n$ diagonalizable matrix with spectrum partition $\lambda = (\lambda_1 \geq \lambda_2 \geq \cdots \lambda_k \geq \cdots \geq \lambda_s)$.  The $x\in {\cal X}_k^r$ iff 
$\lambda_{1}+\ldots +\lambda_{k}\geq n-r$.  
\end{prop}

\noindent 
The proof is obvious.

\begin{prop}\label{prop:oldnil}
Let $x=J_{\theta_1} \oplus \ldots \oplus J_{\theta_s}$ be a nilpotent matrix with partition signature $(\theta_1 \geq \ldots \geq \theta_s)$.  Then $x\in {\cal X}_k^r$ iff $rank(x^k)\leq r$. This is given by the condition $\sum_{i: \theta_i > k} (\theta_i -k)\leq r$. 
\end{prop}

It will be useful to restate this condition in terms of the transpose partition.
\begin{prop}\label{prop:nil}
Let $x=J_\theta$ be a nilpotent matrix with partition signature $\theta=(\theta_1 \geq \ldots \geq \theta_s)$ whose transpose
partition is $\theta^T=(\theta'_1 \geq \cdots \geq \theta'_t)$. 
Then $x\in {\cal X}_k^r$ iff 
$\theta'_{1}+\ldots +\theta'_{k}\geq n-r$.  
\end{prop}

We state our main theorem.
\begin{theorem}
\label{thm:conjugation:main}
Let $x$ be a diagonalizable matrix in $End(X)$ with spectrum partition $\lambda$. The projective orbit closure of $x$ under the conjugation action by $GL(X)$ in ${\mathbb P}(End(X))$  contains nilpotent matrices whose partition signature is $\theta$
iff $\theta \unlhd \lambda^T$.
\end{theorem}
{\bf Proof}
Let $y$ be a nilpotent matrix with partition signature $\theta$ such
that $\theta \not\!\!\unlhd ~\lambda^T$. This implies that
$\lambda \not\!\!\unlhd ~\theta^T=(\theta'_1 \geq \cdots \geq \theta'_s)$. As a result, there is an index $k$ such that
\[ \sum_{j=1}^k \lambda_j > \sum_{j=1}^k \theta'_j \]
Now we choose $r$ such that $\sum_{j=1}^k \lambda_j = n-r$ and
apply Propositions~\ref{prop:diag} and \ref{prop:nil} to conclude
that $x \in {\cal X}_k^r$ and $y \not\in {\cal X}_k^r$. This allows
us to conclude that $y$ is not in the projective orbit closure of $x$.

In view of Remark~\ref{nilp-closure}, it only remains to show
that a nilpotent matrix $J$ with partition signature 
$\lambda^T= (\lambda'_1 \geq \cdots \geq \lambda'_s)$
belongs to the projective orbit closure of $x$. Towards this, we
can easily adapt the argument in the proof of Theorem \ref{theorem:Jn}.
by writing $x$ as a block diagonal matrix $x_1 \oplus \cdots \oplus x_s$
where diagonal matrix $x_i$ is of size $\lambda'_i$ and
has distinct eigenvalues, namely those eigenvalues of $x$ which have
multiplicity at least $i$. We skip the straightforward details.
$\Box$

Now we extend the above theorem to all matrices. Towards this, fix
a matrix $x$ in a Jordan canonical form with $s$ distinct eigenvalues $\mu_1, \mu_2, \ldots, \mu_s$. For eigenvalue $\mu_i$, we associate a partition
$\lambda_i = (\lambda_{i1} \geq  \lambda_{i2} \geq \ldots )$ which
records the sizes of Jordan blocks with diagonal entries $\mu_i$,
in non-increasing order. We now define the transpose block-spectrum partition
of $x$ to be the partition
$\chi = (\chi_1 \geq \chi_2 \geq \ldots )$ where 
$\chi_j = \lambda_{1j} + \lambda_{2j} + \ldots \lambda_{sj}$. Thus, $\chi$ 
is simply the ``sum'' of the partitions $\lambda_1, \ldots, \lambda_s$.
Note that when $x$ is a diagonal matrix, its transpose block-spectrum
partition is nothing but the transpose of its spectrum partition.

\begin{theorem}
Let $x$ be a matrix in $End(X)$ with transpose block-spectrum partition $\chi$. The projective orbit closure of $x$ under the conjugation action by $GL(X)$ in ${\mathbb P}(End(X))$  contains nilpotent matrices whose partition signature is $\theta$
iff $\theta \unlhd \chi$.
\end{theorem}
{\bf Proof} Let $y$ be a nilpotent matrix with partition signature $\theta=(\theta_1 \geq \theta_2 \geq \ldots)$
such that $\theta \not\!\!\unlhd ~\chi$. This means, there is an index $\ell$
with
$$\sum_{j=1}^{\ell}  \theta_j > \sum_{j=1}^{\ell} \chi_j$$
Henceforth we fix the least index $\ell$
with this property. This further implies that $\theta_{\ell} > \chi_\ell$.

Now we set $k=\chi_{\ell+1}$, $r = \sum_{i=1}^{\ell} (\chi_{i} - \chi_{\ell+1})$ and
claim that $x \in {\cal X}_k^r$ and $y \not\in {\cal X}_k^r$.
We first focus on showing that $y \not\in {\cal X}_k^r$. By Proposition~\ref{prop:oldnil}, it suffices to show that
\[ \sum_{i: \theta_i > k} (\theta_i - k) > r \]
As $k = \chi_{\ell+1}$ and $\theta_\ell > \chi_\ell \geq \chi_{\ell+1}$, we have
\[ 
\sum_{i: \theta_i > k} (\theta_i - k) \geq \sum_{i=1}^{\ell} (\theta_i - \chi_{\ell+1})
> \sum_{i=1}^{\ell} (\chi_i - \chi_{\ell+1}) = r
\]

Now we turn our focus on showing that $x \in {\cal X}_k^r$. We make use of the notation developed for defining $\chi$ just before the proof of the theorem. In this notation, we show that
\[
\mathrm{rank} \left( (x - \mu_1 I)^{\lambda_{1,{\ell+1}}}
(x - \mu_2 I)^{\lambda_{2,{\ell+1}}} \ldots 
(x - \mu_s I)^{\lambda_{s,{\ell+1}}}
\right) = r
\]
Note that, this would imply that $x \in {\cal X}_k^r$ as  $k = \chi_{\ell+1}= \lambda_{1,{\ell+1}} + \ldots + \lambda_{s, {\ell+1}}$.
As the generalized eigenspaces corresponding to different eigenvalues are independent and $(x - \mu_i I)$ is nilpotent on the generalized
eigenspace corresponding to eigenvalue $\mu_i$ and invertible on
other generalized eigenspaces, we can use Proposition~\ref{prop:oldnil} to conclude that
\[
\mathrm{rank} \left( (x - \mu_i I)^{\lambda_{i,{\ell+1}}} \right)=
(\lambda_{i,1} - \lambda_{i,\ell+1}) + \ldots + (\lambda_{i,\ell}- \lambda_{i,\ell+1})  =
\sum_{j=1}^{\ell} (\lambda_{i,j} - \lambda_{i,\ell+1}  )
\] and hence, 
\[ \begin{array}{ccl}
\mathrm{rank} \left( (x - \mu_1 I)^{\lambda_{1,{\ell+1}}}
(x - \mu_2 I)^{\lambda_{2,{\ell+1}}} \ldots 
(x - \mu_s I)^{\lambda_{s,{\ell+1}}} \right) & = & 
\sum_{i=1}^{s} \sum_{j=1}^{\ell} (\lambda_{i,j} - \lambda_{i,\ell+1}  ) \\
& = & 
\sum_{j=1}^{\ell} \sum_{i=1}^{s} (\lambda_{i,j} - \lambda_{i,\ell+1}) \\& = & 
\sum_{j=1}^{\ell} \left( ( \sum_{i=1}^{s} \lambda_{i,j}) - (\sum_{i=1}^{s} \lambda_{i,\ell+1}) \right) \\
& = & 
\sum_{j=1}^{\ell} (\chi_j - \chi_{\ell+1}) \\
& = & r
\end{array}
\]
We have thus proved that $x \in {\cal X}_k^r$ and $y \not\in {\cal X}_k^r$. This shows that $y$ is not in the projective orbit closure of $x$. 

In order to complete the proof, we finally show
that a nilpotent matrix of partition signature $\chi= (\chi_1, \ldots, \chi_s)$ belongs to the projective orbit closure of $x$.
We first write $x$ as a block diagonal matrix $x_1 \oplus \cdots \oplus x_s$ where matrix $x_j$ is of size $\chi_j$ and 
$x_j = y_{1,j} \oplus y_{2,j} \cdots \oplus y_{s,j}$ where
$y_{i,j}$ is a Jordan block of size $\lambda_{i,j}$ with
diagonal entry $\mu_i$. The crucial property of each $x_j$ is
that different Jordan blocks in it correspond to different eigenvalues
and hence the minimal polynomial of $x_j$ equals its characteristic
polynomial.
As a result, each $x_j$ can be conjugated to obtain $x'_j$ which
is a companion matrix. We now
consider $x'=x'_1 \oplus \cdots \oplus x'_s$ to which 
$x$ can be block-conjugated.

Now we use the 1-PS family $A(t)= \mathrm{diagonal}(t,t^2,\ldots, t^n)$
and observe that the leading term of $A(t).x'$ is the matrix
$t^{-1}J_\chi$ where $J_\chi=
J_{\chi_1} \oplus \ldots \oplus J_{\chi_s}$ is a nilpotent matrix whose partition
signature is $\chi$. This shows
that $J_\chi$ belongs to the projective orbit closure of $x$ and completes the proof.
$\Box$.

We now conclude this section
with an example.
\begin{ex}
Consider the following matrices
in Jordan canonical forms
\[
x_1 = \left[ \begin{array}{ccc}
1 & 0 & 0 \\
0 & 1 & 0 \\
0 & 0 & -1\\
\end{array} \right] \: \: \:
x_2 = \left[ \begin{array}{ccc}
1 & 1 & 0 \\
0 & 1 & 0 \\
0 & 0 & -1\\
\end{array} \right] \: \: \:
y_1 = \left[ \begin{array}{ccc}
0 & 1 & 0 \\
0 & 0 & 1 \\
0 & 0 & 0 \\
\end{array} \right] \: \: \:
y_2 = \left[ \begin{array}{ccc}
0 & 1 & 0 \\
0 & 0 & 0 \\
0 & 0 & 0 \\
\end{array} \right]
\]
Clearly, $x_1$ is diagonal and
its spectrum partition as well as
transpose spectrum partition is $(2,1)$. 
The transpose block-spectrum
partition of $x_2$ is $(3)$.
Note that $x_1$ is in the affine
orbit closure of $x_2$.
Further, $y_1$ and $y_2$ are nilpotent matrices with partition
signatures $(3)$ and $(2,1)$ respectively. 

The above results allow us to conclude
that the  projective orbit closure
of $x_1$ contains $y_2$ but not $y_1$, and the projective orbit
closure of $x_2$ contains both
$y_1$ and $y_2$.
\end{ex}

\section{Local Differential Geometry of Orbits}
\label{sec:diffgeom}
So far we have explored algebraic techniques which relate properties of the point $y$, the tangent vector $\ell \cdot y \in T_y O(y)$ and the ultimate limit $x$ for paths which arise from 1-PS.
A natural question is whether these paths satisfy additional differential geometric properties, or if there are other special paths which take us to specific limit points $x$ and if the tangent vectors of exit can be determined using the local geometry at $y$.  We formulate this as the following {\bf path problem}.

\begin{defn}
Given a point $y \in V$ and an $x\in \overline{O(y)}$, the {\bf path problem} is to determine a path $\pi$ from $y$ to $x$ satisfying specific properties derived from those of $y$ and $x$. 
\end{defn}

This question is intriguing even when the orbit of $y$ is affine and we wish to reach representative points of codimension $1$ components, i.e., the boundary of the orbit $O(y)\subseteq \PP V$. Specifically, is there an optimization problem formulation to arrive at the codimension $1$ components and the tangents of exit from $y$ leading to these components. 

An important step in this direction would be to develop the local manifold geometry of the orbit. We have already seen in Remark \ref{rem:SFTtheta}, that the function $\Phi : {\cal S}\times N \to {\cal S}$ of the local model is intimately connected to the 2-form $\Pi$. However, one technical issue is that our interest lies in the closure in projective space. In Lemma~\ref{lemma:Tilde} we extend Remark \ref{rem:SFTtheta} to compute the Riemannian curvature tensor on a chart of the projective space containing $y$. This is developed further in Section \ref{sec:cyclic}, where we combine the optimization approach and the computation of the curvatures for a specific problem.  

Let us begin the discussion with the familiar 1-PS situation with a $GL(X)$-module $V$, and $x,y\in V$ with the action of a 1-PS $\lambda (t)$ as shown below:
\[ \lambda (t) \cdot y = x_a \cdot t^a +  x_b t^b + \ldots + x_d t^d \]
Here $x=x_a$ is the projective limit of $y$ via the 1-PS $\lambda (t)$. We will call $a$ as the degree of the limit $x$ under $\lambda $. If $\ell =\log_t(\lambda (t))$ then $\ell \cdot y$ is the tangent of exit of the curve $\lambda (t) \cdot y$ as it exits $y$ and approaches its limit $x=x_a$. We see how the two may be related by a smooth optimization problem. 

\subsection{Degree and the structure of the orbit closure}

In this subsection, $O(y)$ will denote the $SL(X)$-orbit of $y$, $\overline{O(y)}$ its closure in $V$, and $\PP \overline{O(y)} \subset \PP V$ will denote the projective closure of the orbit $O(y)$. The structure of $\PP \overline{O(y)}$ consists of closures $\PP \overline{O(x)}$ of various lower dimensional orbits $O(x)$ of limit points $x$ of the orbit $O(y)$. 
If $O(y)$ is closed in $V$ and the stabilizer $K$ of $y$ is reductive, then in fact, the set $\PP \overline{O(y)}- \PP O(y)$ is the union of a finite number of co-dimension 1 projective varieties $\PP \overline{O(x_i )}$. This offers the first decomposition of the limit points of $O(y)$. 

And yet, while the boundary of $\PP O(y)$ may be of small co-dimension, the limit $x$ of a typical 1-PS $\lambda (t)$ acting on $y$ may be embedded deep within this boundary as the following examples show. 

\begin{ex}
Let $A$ be an $n\times n$-matrix such that $A(1,n)\neq 0$. Let $\lambda (t)$ be the 1-PS consisting of diagonal matrices $diag([1,t,\ldots, t^{n-1}])$. 
Then the projective limit of $\lambda(t) A \lambda (t)^{-1}$ is the matrix $E_{1n}$ with degree $-(n-1)$. On the other hand, we know that if $A$ has distinct eigenvalues then the orbit $\PP \overline{O(J_n )}$ of $J_n$ is of co-dimension 1 within $\PP \overline{O(A))}$. The point $J_n$ is approached by a special 1-PS which is conjugate to $\lambda (t)$ and where the degree of $J_n$ is $-1$. 
\end{ex}

\begin{ex} 
Consider the variables $X=\{ x_1 ,\ldots ,x_9 \}$ and the form $det_3 (X) \in Sym^3(X)$. The dimension of the stabilizer $K\subseteq GL(X)$ of $det_3$ is $16$. There is another basis ${\cal X}=\{X_1 ,\ldots , X_9\}$ of $\C X$ and a 1-PS $\lambda_2 (t)$ with the following properties ( see \ref{subsection:det3} for details). 
Let $Y=\{X_1,X_2 ,X_3 \}$ and $Z=\{ X_4 ,\ldots ,X_9 \}$. Let $\lambda_2 (t)$ be such that $\lambda_2 (t)(y)=y$ for $y\in Y$, while $\lambda_2 (t)(z)=tz$ for all $t\in Z$. Then $\lambda_2 \cdot det_3 ({\cal X})= tQ_2 ({\cal X}) + t^3 Q_3 ({\cal X})$, where $Q_2 ({\cal X})=X_1^2 X_4 +X_2^2 X_5 +X_3^2 X_6 +X_1 X_2 X_7 + X_2 X_3 X_8 +X_1 X_3 X_9 $. 
Moreover, the degree of $Q_2$ is $1$ and its orbit is one of the co-dimension 1 forms in the projective closure of the orbit of $det_3$, and the dimension of its stabilizer is 17.

On the other hand, for $\lambda'_2 (t)$ defined as $\lambda'_2 (t)x_i =x_i$ for $i=1,2,3$ and $\lambda'_2 (t)(x_j )=tx_j$, for $j=4,\ldots, 9$, we have the limit $Q'_2 =x_1 x_2 x_3$ of degree $0$. The stabilizer of $Q'_2$ is of dimension $19$. Indeed, for most $\mu (t)$ conjugate to $\lambda_2 $, we will have $\mu (t) (det_3 (X))= t^0 P_0 +t^1 P_1 +t^2 P_2 +t^3 P_3$, i.e., with the degree of the leading term $P_0$ as zero, which is lower than that of $Q_2$. Moreover, the dimension of the orbit of this generic limit is lower than that of $Q_2$.  \end{ex}

Thus, an important problem is to locate directions, which in the algebraic case, are 1-PS, to specific limit points. A foundation for such analysis was provided by Kempf, \cite{kempf1978instability}, in the case when $\lambda (t) \subseteq SL(X)$ and the degree $a$ at which the limit is found is non-negative. The case of $a>0$ corresponds to $y$ being in the nullcone with the limit $x=0$. The case $a=0$ corresponds to a limit $x\neq 0$, with the orbit of $O(x)$ lying in the affine closure of the $SL(X)$-orbit of $y$. 

A key construction in Kempf's analysis is the computation of an optimal 1-PS $\lambda (t)$ within a fixed torus $T$. The uniqueness of this $\lambda $ within one torus allows its comparsion with optimal points in other tori, and the identification of a canonical parabolic subgroup $P(\lambda )\subseteq GL(X)$ to be associated with $y$ and $x$. The first computation has a generalization to projective limits, which we outline below.

Let us fix a maximal torus $T\subseteq GL(X)$. For simplicity, let us assume that it is the space of invertible diagonal matrices $D_n$. Let ${\cal D}_n$ be the Lie algebra of $D_n$. Let 
\[ O_{n-2} =\{ p \in \R^n \mbox{ such that } | \sum_i p_i^2 =1 \mbox{ and } \sum_i p_i =0\}\]
For any $p\in \R^n$, let us denote by $\overline{p}$ the $n\times n$ diagonal matrix with $p$ as its diagonal. 
We regard the action of $\overline{p}$ as an element of $gl(x)$, the Lie algebra of $GL(X)$. Let $t^p$ denote the diagonal matrix with diagonal $[t^{p_1}, \ldots ,t^{p_n} ]$. Note that for $p\in O_{n-2}$, we have $\overline{p}\in {\cal D}_n \subseteq sl(X)$ and $ t^p \in D_n \subseteq SL(X)$. 

Let $\Xi \subseteq {\cal D}_n^*$ be the weight-space of the module $V$ under $T$. For any $v\in V$ let $v=\sum_{\chi \in \Xi} v_{\chi}$ be the decomposition of $v$ by this weight-space. 
Let $\Xi (v)$ be those $\chi \in \Xi$ for which $v_{\chi} \neq 0$. Thus $\Xi (v)$ is the support of $v$. To simplify notation we omit the subscripts in ${\cal D}_n, D_n$.
\begin{lemma}
For any element $\ell \in {\cal D}$, we have:
\[ \begin{array}{rcll}
\ell \cdot v &=& \sum_{\chi} \langle \ell , \chi \rangle v_{\chi} & \mbox{ via the Lie algebra action} \\
t^{\ell} \cdot v &=& \sum_{\chi} t^{\langle \ell , \chi \rangle} v_{\chi} & \mbox{ via the group action}\\
\end{array} \]
\end{lemma}

\noindent
We now state Kempf's basic optimization result. 
\begin{prop}
(Kempf) We say that $v\in V$ is the the null-cone (with respect to ${\cal D}$), if there is an $\ell \in {\cal D}$ such that for all $\chi\in \Xi (v)$ we have $\langle \ell ,\chi \rangle >0$. Let $\mu (\ell ,v)= \min \{ \langle \ell ,\chi \rangle | \chi \in \Xi (v) \}$. 
If $v$ is unstable then there is a unique $\ell_0 \in sl(X)$ such that $\mu(\ell_0 ,v)>0$ and $\mu (\ell_0 ,v) > \mu (\ell ,v)$ for all $\ell \neq \ell_0$. 
\end{prop} 

\noindent
We have the following alternate optimization formulation.
\begin{prop}
\label{prop:diff-kempf}
Let $v \in V$ be arbitrary. Let ${\cal L}_{\alpha} (v)=\{ \ell \in {\cal D} | \mu (\ell ,v)>\alpha \}$. Suppose that ${\cal L}_{\alpha} (v)$ is non-empty. Consider the function 
\[ f(t, \ell) = \sum_{\chi \in \Xi} \| v_{\chi} \|^2 t^{-\langle \ell , \chi \rangle} \]
For fixed $v\in V$ and $t \in \R$ the above function is a smooth function on the compact set $O_{n-2}$. Let $\ell_f (t)$ minimize $f$ on $O_{n-2}$. 
Then there is a constant $t_0 >1$ depending on $v$, such that for all $t>t_0$ we have $\ell_f (t)\in {\cal L}_{\alpha} (v)$.  
\end{prop}

\noindent 
{\bf Proof}: Let $\ell_K \in {\cal L}_{\alpha} (v)$ be any element and let $\mu_K =\mu (\ell_K ,v)$. Note that $\mu_K >\alpha $. Whence, for $t>1$, we have $f(t ,\ell_K ) < \| v \|^2 t^{-\mu_K}$. 

Let $\Xi_- (t) \subseteq \Xi(v)$ be those $\chi \in \Xi (v)$ such that $\langle \ell_f (t) ,\chi \rangle \leq \alpha$. Similarly, let $\Xi_{+} (t)$ be those $\chi \in \Xi (v)$ such that $\langle \ell_f (t) ,\chi\rangle >\alpha$. The claim is that $\Xi_- (t) $ is empty for $t$ sufficiently large. 

Suppose not and there is a $\Xi_- \subseteq \Xi (v)$ such that $\Xi_- =\Xi_- (t)$ for a divergent sequence $(t_i )_{i=1}^{\infty}$ with $t_i >1$ for all $i$. Then there is a $\nu \in \Xi_- \subseteq \Xi (v)$ such that $\langle \ell_f (t_i ) , \nu \rangle \leq \alpha$ for all $i$. Whence $f(t_i ,\ell_i )\geq \| v_{\nu} \|^2 t^{-\alpha} >0$ for all $t_i$. But $f(t, \ell_K ) <\| v \|^2 t^{-\mu_K}$ for all $t$. 
Thus we have:
\[ \| v_{\nu} \|^2 t^{-\alpha}\leq f(t_i ,\ell_i ) \leq f(t_i , \ell_K ) <\| v \|^2 t_i^{-\mu_K} \]
This gives us $t_i^{\mu_k -\alpha} < \| v\|^2 /\| v_{\nu}\|^2$, with $\mu_K>\alpha $. 

This is a contradiction since we can choose a $t_i$ large enough for this to not hold. This proves the proposition. $\Box $ 

\begin{remark}
For any $\ell \in {\cal L}_{\alpha} (v)$, let $\alpha_0 = \mu (\ell ,v)$ and $\Xi_{\ell} =\{ \chi \in \Xi (v)| \langle \ell ,\chi \rangle =\alpha_0 \}$. Let $v_{\ell} =\sum_{\chi \in \Xi_{\ell}} v_{\chi}$. Then $v_{\ell}$ is a projective limit of $v$. We do see that larger values of $\alpha$ do lead to more fundamental limit points of $v$. 
\end{remark}

Thus, if $x$ is the minimum degree limit for the point $y$, and is approachable by a 1-PS within $T$, then the optimization function $f(t,\ell)$ on ${\cal D} \subseteq sl(X)$ does indeed offer an algorithm to discover this. 
However, in the general case, this must be done across all tori within $sl(X)={\cal G}$ in a unified manner. An example where this has been achieved is in designing efficient algorithms to decide  intersections of orbit closures in the space $(Mat_n)^k$, i.e., $k$-tuples of matrices with $GL_n \times GL_n$ acting diagonally on the the left and right. ~\cite{allen2018operator} formulate this problem as a geodesic convex optimization problem on the space of positive definite
matrices with their natural hyperbolic metric.

Thus, what is needed are global optimization functions $\Tilde{f}(\mf{g},y)$, where $\mf{g} \in {\cal G}$. 
or more specifically on ${\cal G}/{\cal K}$ where ${\cal K}$ is the stabilizer of $y$. Since ${\cal G}/{\cal Y}$ is the tangent space of the $GL(X)$-orbit of $y$, with a ${\cal K}$-structure, the differential geometric structure of this space needs closer examination. 
The Riemannian approach, and its associated tensors and the geometry of geodesic paths offers an alternate analysis of this algebraic object. We begin with curvature of the orbit $O(y)$ in the next subsection. 

\subsection{The Curvature Tensors}
\label{sec:diff}
We now look at the differential geometry of $y$ as a point in $\PP V$. 

Remark \ref{rem:SFTtheta} 
already connects the algebraic gadgets $\theta $ and $\Phi$ with the local 2-form $\Pi$ of the orbit $O(y)\subset V$. Thus, it is worthwhile to investigate if there is a differential geometric equivalent of the triple stabilizer condition. The hope is that the local geometry may guide us to identify potential tangent vectors and standard paths (e.g., along 1-PS or geodesics) which lead us to various limit points. 

However, one difficulty in pursuing this further is that while the Lie algebra acts on $V$, the topology and the computation of the limit happens in $\PP V$. Whence, we must develop a local chart ${\cal C}$ for the point $y$ as a member of $\PP V$ and transfer the action of ${\cal G}$ to this chart. 

The approach we adopt here is to restrict this computation to the real case. Whence, we assume that $V$ is equipped with a positive definite metric $\langle \cdot , \cdot \rangle$. Moreover, let us assume that a basis for $V$ is chosen such that $\langle v,v' \rangle =v^T v'$, for any $v,v' \in V$ (treated as column vectors).

The chart ${\cal C}$ is constructed from the ball $B_{\epsilon} (y)$ as follows:
\[ {\cal C}= \{ v\in B_{\epsilon} (y) \mbox{  such that  } \| v\| =\|y\| \} \]
We define the projection $\pi_{\cal C} : B_{\epsilon} (y) \rightarrow {\cal C}$ as 
\[ \pi_{\cal C} (v)= \frac{\| y \|}{\| v\|} v \]

\begin{defn}
We say that the representation $\rho :GL(X)\rightarrow GL(V)$ is of (A) degree zero if $\R^* I_X $ is in the kernel of this map. Otherwise, (B) we say that $\rho$ has a non-zero degree. In this case, and if $V$ is irreducible, then there is an integer $d \neq 0$ such $\rho (tI_X)=t^d I_V$. 

We say that a point $v$ is simple if in the case (A) above, the orbit $O(y)$ does not contain any other multiple of $y$, or in the case (B) above, the condition $\mf{g} \in {\cal G}$ and $\mf{g}\cdot v \in \R v$ implies that $\mf{g} \in \R I_X$. In this section, in case (B), we restrict ${\cal G}$ to $sl(X)$ and the orbit $O(v)$ to be the $SL(X)$-orbit of $v$. 
\end{defn}
Let us define $\Tilde{O}(y)$ as the image of the  projection of $O(y)\cap B_{\epsilon} (y)$ onto ${\cal C}$. We have the following simple lemma:
\begin{lemma}
Suppose that $y$ is a simple point. Then we may choose $\epsilon>0$ such that all points $v\in B_{\epsilon} (y)$ are also simple. Then the projection $\pi_{\cal C}: O(y) \rightarrow \Tilde{O}(y)$ is a diffeomorphism at all points. 
\end{lemma}

\noindent 
{\bf Proof}: The condition of simplicity is equivalent to $v \not \in {\cal G}\cdot v$ for (A) ${\cal G}=gl(X)$ or (B) ${\cal G}=sl(X)$. This is clearly an open condition. The second condition is immediate. $\Box $

The above lemma allows us to define $\rho_{\cal C}: {\cal G} \rightarrow Vec({\cal C})$ given by $\rho_{\cal C} (\mf{g})(v)= \pi_{\cal C}^* (\rho (\mf{g})(v))$. Our next task is to compute the second fundamental form $\Pi_{\cal C}$ on ${\cal C}$. As before, we will use the vector fields $\rho_{\cal C}(\mf{s})$ for $\mf{s}\in {\cal S}$, to compute this.

\begin{lemma} \label{lemma:Tilde}
Let $y\in \Tilde{O}(y) \subset {\cal C}\subset V$ be a simple point with stabilizer ${\cal K}$ and complement ${\cal S}$ of dimension $K$. 
Let $X=\rho(\mf{s})(v)= \nabla S v$, be the linear vector field, in the notation of Remark  ~\ref{rem:SFTtheta}. For $i=1,\ldots ,K$, let $\mf{s}_i$ be a basis of ${\cal S}$, $X_i =\rho (\mf{s}_i) =\nabla S_i v$ and  $\Tilde{X_i}=\pi^*_{\cal C}(X_i )=\rho_{\cal C} (\mf{s}_i )$. Then:
\begin{enumerate} 
\item The normal space $\Tilde{N}$ to $\Tilde{O}(y)$ within ${\cal C}$ is $\pi^*_{\cal C}(y) (N)$.
\item Let $\pi^*_{\cal C}(y)$ be the projection at point $y$ and $\Tilde{S}_i =\pi^*_{\cal C} (y) S_i$, then:
\[ (D_{\Tilde{X}_j} \Tilde{X}_i )(y)= - \frac{y^T S_i y}{y^T y} \Tilde{S}_j y + \Tilde{S}_i \Tilde{S}_j y \]
\item Then the second fundamental form $\Pi_{\cal C}$ of $\Tilde{O}(y)$ as a submanifold of ${\cal C}$ is given as follows: 
\[ \Pi_{\cal C} (\Tilde{X}_j ,\Tilde{X}_i)(y) =\lambda_{\Tilde{N}} (D_{\Tilde{X}_j} \Tilde{X}_i )(y) =\lambda_{\Tilde{N}} \Tilde{S}_i \Tilde{S}_j y \]
\item Finally, if $TO(y)=T\Tilde{O}(y)$, i.e., the orbit and its projection osculate at $y$, and $\langle D_{X_i} y, y \rangle =0$, i.e., $y^T S_i y=0$ for all $i$, then:
\[ \Pi_{\cal C} (\Tilde{X}_j ,\Tilde{X}_i)(y)= \pi_{\cal C}^* (y)(\Pi (X_j , X_i )(y)) \]
\end{enumerate}
\end{lemma}

\noindent 
{\bf Proof}: The first part is easily proved by taking up cases (A) and (B) separately. 

Next, the field $\Tilde{X}_i =\rho_{\cal C} (\mf{s}_i )$ is given by:
\[ \Tilde{X}_i = \pi_{\cal C}^* (v) (\nabla S_i v) = \nabla (I-\frac{1}{v^T v}vv^T )S_i v \]
Since $\pi^*_{\cal C}(y)=(I-\frac{1}{y^T y}yy^T ) $, we have:
\[ \Tilde{X}_j (y)= \nabla \pi^*_{\cal C}(y) S_i y = \nabla \Tilde{S_i } y\]
We then have:
\[ \begin{array}{rcl}
(D_{\Tilde{X}_j} \Tilde{X}_i )(y)&=& \pi^*_{\cal C} (y) \circ [[ {\Tilde{X}_j}(y) (I-\frac{1}{v^T v} vv^T)(y)] S_i y + (I-\frac{1}{y^T y} yy^T) [{\Tilde{X}_j }(y) (S_i v)(y)] ]\\
&=& \pi^*_{\cal C}(y) \circ [[ 2 (\Tilde{S}_j y^T y) yy^T S_i y -\frac{1}{y^T y}[ \Tilde{X}_j (y) (vv^T)(y)] S_i y + \pi^*_{\cal C}(y) S_i \pi^*_{\cal C}(y) S_j y]]\\
&=& \pi^*_{\cal C}(y) [- \frac{1}{y^T y} \Tilde{S}_j yy^T S_i y -\frac{1}{yT y} yy^T \Tilde{S}_j^T  S_i y + \Tilde{S}_i \Tilde{S}_j y]\\
&=& \pi^*_{\cal C} (y) [- \frac{y^T S_i y}{y^T y} \Tilde{S}_j y -\frac{y^T\Tilde{S}_j S_i y}{yT y} y + \Tilde{S}_i \Tilde{S}_j y]\\
&=& - \frac{y^T S_i y}{y^T y} \Tilde{S}_j y + \Tilde{S}_i \Tilde{S}_j y\\
\end{array} 
\]
This proves (2). Now, since $\Tilde{N}=\pi^*_{\cal C}(y)(N)$, we have $\Tilde{S}_j y \in T_y \Tilde{O(y)}$, we have $\lambda_{\Tilde{N}} (\Tilde{S}_j y) =0$. Whence we have:
\[ \Pi_{\cal C} (\Tilde{X}_j ,\Tilde{X}_i)(y) =\lambda_{\Tilde{N}} (D_{\Tilde{X}_j} \Tilde{X}_i )(y) =\lambda_{\Tilde{N}} \Tilde{S}_i \Tilde{S}_j y \]
This proves (3). In the case when $O(y)$ and $\Tilde{O}(y)$ osculate, i.e., $TO_y =T\Tilde{O}_y$, i.e., the tangent planes match, we have $\Tilde{S}_i =\pi^*_{\cal C} (y) S_i =S_i$, whence
\[ \Pi_{\cal C} (\Tilde{X}_j ,\Tilde{X}_i)(y)= \lambda_{\Tilde{N}} S_i S_j y = \pi^*_{\cal C}(y) \lambda_N S_i S_j y = \pi_{\cal C}^* (\Pi (X_j , X_i )(y)) \]
This proves (4). $\Box $

We now recall the fundamental result on Riemannian curvature tensor for smooth manifolds.

\begin{prop}
Let $M$ be a submanifold of a Riemannian manifold $P$. Let $X_i ,X_j $ and $X_k$ be vector fields on $M$. Then the Riemannian curvature tensor of $M$ is given by:
\[ R(X_i ,X_j )X_k= \overline{D}_{X_i} \overline{D}_{X_j} X_k -\overline{D}_{X_j} \overline{D}_{X_i} X_k -\overline{D}_{[X_i ,X_j ]} X_k \]
Moreover, if $X_l $ is another vector field on $M$ then 
The quantity $r_{ijkl}=\langle R(X_i ,X_j )(X_k)(x),X_l (x)) \rangle$ equals the quantity $\langle \Pi(X_i ,X_l)(x), \Pi(X_j ,X_k)(x) \rangle -\langle \Pi (X_j ,X_l)(x),\Pi(X_i ,X_k)(x)\rangle$.
\end{prop}

The proof of this is standard.

\begin{ex}
Let us consider $\R^n$ with the standard basis $e_1 ,\ldots ,e_n$ and under the action of $SO(n)$ and Lie algebra $so(n)$. Consider  the vector $x=(r,0,\ldots ,0)=r\cdot e_1$, whose orbit is $S=rS^{n-1}$, the sphere of radius $r$ in $\R^n$. Note that in this case, $O(x)=\Tilde{O}(x)$. 
The tangent space $T_x S$ is the span of $e_2 ,\ldots ,e_n$. We have $\R^n =\R e_1 \oplus T_x S$ and thus the normal space is only $1$-dimensional. Let $\mf{s}_i \in so(n)$ be the element with matrix $S^i$ with the only non-zero entries as $S^i (1,i)=-1/r, S^i (i,1)=1/r$. 
Thus $X_i =\nabla S^i \overline{v}$ is the corresponding vector field such that $X_i (x)=e_i$.   

For $i=2,\ldots,n$ and $j=2,\ldots,n$, we may build that data $c^i_j$, i.e., the ``normal'' part of $j$-th column of the matrix $S_i$. Since the normal part is only $1$-dimensional, i.e., the first row, we have:
\[ \Pi (X_j ,X_i )=c^i_j =\left\{ \begin{array}{lr}
0 & \mbox{if $i\neq j$} \\
-1/r & \mbox{if $i=j$}
\end{array} \right. \]
The Ricci curvature $R^c_{jk}$ is the sum below:
\[
R^c_{jk} =\sum_i r_{ijki}=\sum_i c^i_i c^k_j - c^i_j c^k_i \]
This gives us:
\[ R^c_{jk} =\left\{ \begin{array}{lr}
0 & \mbox{if $j\neq k$} \\
(n-1)/r^2 & \mbox{if $j=k$}
\end{array} \right. \]
\end{ex}

\begin{ex}
Let us consider the adjoint action of $GL_n$ on $Matn_n$, i.e., $\rho (A)(X)=AXA^{-1}$, where $A\in GL_n$ and $X\in Matn_n$, and the corresponding Lie algebra action of $gl_n$, given by $\rho (a)(X)=aX-Xa$, where $a\in gl_n$. 

Let $x=diag(\overline{\lambda})$, the diagonal matrix with distinct eigenvalues $\overline{\lambda} =(\lambda_1 ,\ldots, \lambda_n )$. Then ${\cal H}={\cal D}$, the Lie algebra of all diagonal matrices and ${\cal S}$ is the linear span of $\{ E_{ij} | i\neq j\}$, the non-diagonal elementary matrices. Note that $[E_{ij},x]=(\lambda_j -\lambda_i )x$ and thus ${\cal S}$ is both the complement to ${\cal D}$ within $gl_n$ and to complement in $Matn_n$ to $T_x O$ the tangent space of the orbit. 
Note that in the standard inner product on matrices, we have $Tr([E_{ij},x]x^T)=0$ and thus the point $x$ satisfies the osculation conditions of Lemma \ref{lemma:Tilde} (4). 

Let $e_{ij}=E_{ij}/(\lambda_j -\lambda_i)$ and note that $\rho (e_{ij})(x)= e_{ij}$. We now compute $\Pi (e_{rs},e_{pq})=c^{pq}_{rs}(x)=\lambda_N D_{e_{rs}} \rho(e_{pq})(x)$. This is given by $\lambda_N$ of the column vector corresponding to $e_{rs}$ in the matrix $\rho (e_{pq})$. Since $N$ consists of diagonal matrices, this equals the diagonal entries in the action of $e_{pq}$ on $Mat_n$, in other words, the diagonal matrix in $[e_{pq},e_{rs}]$. Thus, the only $(r,s)$ for which this happens is $(q,p)$ and thus we see that:
\[ c^{pq}_{rs} =\left\{ \begin{array}{lr}
0 & \mbox{if $(rs)\neq (qp)$} \\
d_{pq}& \mbox{if $(rs)=(qp)$}
\end{array} \right. \]
where $d_{pq}$ is the diagonal matrix below:
\[ \left\{ \begin{array}{ll}
d_{pq}(p,p)&= -(\lambda_q -\lambda_p )^2 \\
d_{pq}(q,q)& =(\lambda_q -\lambda_p )^2 \\
d_{pq}(i,j)&=0 \mbox{   for all other tuples $(i,j)$} \\
\end{array} \right. \]
Finally, note that $\langle d_{pq}, x \rangle \neq 0$ unles $\lambda_p \neq \lambda_q$, $\Pi_{\cal C} (c_{rs}^{pq})$ will need a final projection which is skipped here.
\end{ex}

\begin{ex}
Let us consider the case when $n=2m$ in the example above and
\[ x=\left[  \begin{array}{cc}
\lambda I_m & 0 \\
0 & \mu I_m \\
\end{array} \right] \mbox{    with $\lambda ,\mu \neq 0$}\]
Clearly, the stabilizer ${\cal H}$ of $x$ is the Lie algebra $gl_m \times gl_m $ and ${\cal S}$ is the space below:
\[ {\cal S}= \left\{ =\left[  \begin{array}{cc}
0 & X \\
Y & 0 \\
\end{array} \right] | X,Y \in gl_m \right\}  \]
Note that ${\cal S}=T_x O$ and that $N$ is the space of all block-diagonal matrices:
\[ N= \left\{ =\left[  \begin{array}{cc}
Z & 0 \\
0 & W \\
\end{array} \right] | Z,W \in gl_m \right\}  \]
We use the ${\cal H}$-invariant inner product $\langle A,B \rangle =Tr(AB)$ on the space $gl_n$. We decompose ${\cal S}={\cal S}^+ \oplus {\cal S}^-$, as given below:
\[ {\cal S}^+= \left\{ =\left[  \begin{array}{cc}
0 & X \\
0 & 0 \\
\end{array} \right] | X \in gl_m \right\} \mbox{ and    } {\cal S}^-= \left\{ =\left[  \begin{array}{cc}
0 & 0 \\
Y & 0 \\
\end{array} \right] | Y \in gl_m \right\}  \]
For matrices $X,Y \in {\cal S}^{\pm}$, it is easy to see that $\Pi (X,Y)(x)=\lambda_N D_X Y (x)$ is precisely the block diagonal component of $[X,Y]$. Whence we have:
\[ \Pi (X,Y)(x)=\left\{ 
\begin{array}{rl}
\left[ \begin{array}{ll}
XY & 0 \\
0 & -YX
\end{array} \right] & \mbox{when $X\in {\cal S}^+, Y\in {\cal S}^-$} \\ \noalign{\vskip 3mm}
\left[ \begin{array}{cc}
-YX & 0 \\
0 & XY
\end{array} \right] & \mbox{when $X\in {\cal S}^-, Y\in {\cal S}^+$} \\ \noalign{\vskip 3mm}
0 & \mbox{otherwise} 
\end{array} \right. \]

\end{ex}

\subsection{The Cyclic Shift Matrix}
\label{sec:cyclic}
We now illustrate an optimization formulation of the path problem and its relationship with local curvatures. 

Let $X$ be the space $Mat_n$ as before acted upon by $gl_n$ my conjugation. Let us examine the case of the approach of a special matrix to its projective limit point expressed in a basis in which the required 1-PS is easily written.

Let $\Z_n =\{ 0,\ldots ,n-1\}$ be the set of integers under modulo $n$ addition. This will also be an index set for defining the matrix $\mf{c}_n$ (or simply $\mf{c})$ as the cyclic shift matrix:
\[ \mf{c} (i,j)=\left\{ \begin{array}{rl}
1 & \mbox{   if $j-i=1$} \\
0 & \mbox{   otherwise} \end{array} \right. \]
Thus $\mf{c}$ is the matrix below:
\[ \mf{c} = \left [ \begin{array}{ccccc}
0 & 1 & \ldots & & 0 \\
0 & 0 & 1 & \ldots & 0 \\
\vdots & & \vdots & & \\ 
0 & \ldots & 0 & 0 & 1 \\
1 & \ldots & 0 & 0 & 0 \\
\end{array} \right] \]
Let $\lambda (t)$ be the $1$-PS below:
\[ \lambda (t) =\left [ \begin{array}{ccccc}
t & 0 & \ldots & & 0 \\
0 & t^2 & 0 & \ldots & 0 \\
\vdots & & \vdots & & \\ 
0 & \ldots & 0 & 0 & t^n \\
\end{array} \right] \]
We the see that $\log (\lambda (t))$ is the matrix $\ell $ below:
\[ \ell = \left [ \begin{array}{ccccc}
1 & 0 & \ldots & & 0 \\
0 & 2 & 0 & \ldots & 0 \\
\vdots & & \vdots & & \\ 
0 & \ldots & 0 & 0 & n \\
\end{array} \right] \]
The action of a generic element $\mf{a}=(a_{ij})$ of $gl_n$ on $\mf{c} $ is given by:
\[ t_n = \mf{a} \cdot \mf{c} = a\cdot \mf{c} - \mf{c} \cdot a =
\left [ \begin{array}{ccccc}
a_{0,n-1}-a_{10} & a_{00}-a_{11} & \ldots & & a_{0,n-2}-a_{1,n-1} \\
\vdots & & \vdots & & \\ 
a_{n-1,n-1}-a_{0,0} & a_{n-1,0}-a_{0,1}&\ldots  &  & a_{n-1,n-1}-a_{0,n-1} \\
\end{array} \right] \] Thus the matrix $t_n$ is given by:
\[ t_n (i,j)= a_{i,j-1}-a_{i+1,j} \]
We use the standard inner product on both $Mat_n$ as well as $gl_n$, viz., $\langle x,y \rangle =Tr(x\overline{y}^T)$, for any matrices $x,y$. Since $(\mf{c}^i)^T =\mf{c}^{n-i}$, we have the following straightforward lemma:
\begin{lemma}
The tangent space $T_{n} X$ to $\mf{c}$ is perpendicular to $\mf{c}$. In other words $Tr(t_n \mf{c}^T)=0$ for all $\mf{a} \in gl_n$. 
\end{lemma}
The stabilizer ${\cal H}$ of $\mf{c}$ is the space $\oplus_{i=0}^{n-1}  \mf{c}^i$. Let us compute the complement ${\cal S}$ to ${\cal H}$ under the above inner product. For $k=0,\ldots, n-1$, let $D_k$ be the indices 
\[ D_k =\{ (i,i+k)|i=1,\ldots, n\} \]
Thus $D_k$ is $k$-th shifted diagonal, i.e., the support $c_n^k$. 
Thus \[ {\cal S}= \{ \mf{s}=(s_{ij}) \mbox{ such that } \sum_{(i,j)\in D_k} s_{ij} =0 \mbox{ for all $k$} \} \] 
Let us look at $\gamma (\mf{s})= \| \mf{s}\cdot \mf{c} \| /\| \mf{s} \|$, the {\em compression} achieved by $\mf{s}$. We show that $\ell$ above is one of the elements of ${\cal S}$ which minimize $\gamma $. For convenience, let us denote by ${\cal S}_k$ as the elements $(s=(s_{ij})\in {\cal S}$ where $j-i=k$, i.e., matrices in ${\cal S}$ with support in $D_k$.  Thus, ${\cal S}=\oplus_{k=0}^{n-1} {\cal S}_i$. 

Let us assume that $\gamma $ achieves a minima at $\mf{x} \in {\cal S}$ and that $\mf{t}=\mf{x}\cdot \mf{c}$. We may decompose these elements as above, by their support and express $\mf{x}$ and $\mf{t}$ as a sum $\mf{x}=\sum_{i=0}^{n-1} \mf{x}_i$ and $\mf{t}=\sum_{i=0}^{n-1} \mf{t}_i$. We observe that:
\[ \gamma (\mf{x})=\frac{\sqrt{\|\mf{t}_0 \|^2 + \ldots +\| \mf{t}_{n-1}\|^2}}{\sqrt{\|\mf{x}_0 \|^2 + \ldots +\| \mf{x}_{n-1}\|^2}}  \]
Since $\mf{x}_i \cdot \mf{c} =\mf{t}_{i+1}$, the minima is also achieved where $\mf{x}=\mf{x}_i$ for some $i$. We choose $i=0$; the other cases are similar. 
Thus $\mf{x}\in {\cal S}_0$ is a diagonal matrix $(x_{ii})$ such that $\sum_i x_{ii}=0$. We see then that:
\[ \gamma (\mf{x}) = \frac{\sum_{i=1}^n (x_{ii}-x_{i+1,i+1})^2}{\sum_{i=1}^n x_{ii}^2} \]
The condition that $\sum_i x_{ii}=0$ requires that some $x_{ii}\neq x_{i+1,i+1}$ for some $i$. Suppose then that $x_{nn}-x_{11}=1$.  It is then easily seen that the optimal distribution of the intermediate $x_{ii}$ is obtained by evenly spacing them, i.e., $x_{i+1,i+1}-x_{ii}=1/(n-1)$, for $ i=1,\ldots,n-2$.

Recall $\ell$ as before and let $\overline{\ell}$ be obtained by subtracting constant $\frac{(n+1)}{2}$ from each entry of $\ell$. This makes $\overline{\ell}$ an element of ${\cal S}$. We have the following lemma, which is easily proved.
\begin{lemma} \label{lemma:Lbasis}
Let $\ell_{ij}=\mf{c}^i \overline{\ell} \mf{c}^{-j}$. 
\begin{itemize}
    \item Then $\ell_{ij}\in {\cal S}_{i-j}$ and the set $L_k = \{ \ell_{ij} | i-j =k \}$ is a basis for ${\cal S}_i$. $L=\cup_k L_k$ is a basis for ${\cal S}$. 
\item The elements $\mf{x}\in {\cal S}_k$, for some $k$, optimizing $\gamma $ above are scalar multiples of some element of $L$.
\end{itemize}
\end{lemma}

The above lemma is not a completely happy situation. For $\mf{r}\in {\cal S}_0$ with $\ell_{ii}=\mf{c}^i \overline{\ell} \mf{c}^{-i}$, we see that the leading term of $\ell_{ii} (\mf{c})= [\ell_{ii} , \mf{c}]=\ell_{ii} \mf{c}-\mf{c} \ell_{ii}$ is $\mf{c}^i J_n \mf{c}^{-i}$, which is a cyclic renumbering of the nilpotent flag of $J_n$ and thus $\ell_{ii}$ is not only an optimal tangent of exit, but it also leads to the desired highest dimension limit point. This is the good part.

However, for $i\neq j$, $\ell_{ij}=\mf{c}^i \cdot \overline{\ell} \cdot \mf{c}^{-j} \in {\cal S}_{j-i}$ does  minimize $\gamma$, but we have the following:
\[ LT([\ell_{ij} ,\mf{c}])=\ell_{ij} \mf{c}-\mf{c} \ell_{ij} =(\mf{c}^i J_n \mf{c}^{-1})c^{i-j} \]
This limit is not conjugate to $J_n$, and need not even be a nilpotent matrix. Thus the direction $\ell_{ij} \in {\cal S}_{i-j}$ is a false start from the point $\mf{c}$. The resolution of this situation lies in the analysis of the curvature form at the point $\mf{c}$. Note that thus form $\Pi : {\cal S} \times {\cal S} \rightarrow N$, where we have the basis $\cup_k L_i$ for ${\cal S}$ and $\{ \mf{c}^i | i=0,\ldots ,n-1 \}$ for $N$. 

Our next task will be to find an expression for $\langle \Pi (\ell_{i,j}, \ell_{k,l}), \mf{c}^r \rangle$ for general $i,j,k,l,r$ and to evaluate these for a suitable subset. For this, we will use the lemma:
\begin{lemma} 
With the above notation, we have:
\[ \langle \Pi (\ell_{i,j}, \ell_{k,l}), \mf{c}^r \rangle
=\langle \lambda_N (\rho (\ell_{k,l})(\rho(\ell_{i,j}), (\mf{c}^r))\rangle =\frac{1}{n} Tr([\ell_{kl} , [\ell_{ij}, \mf{c}^1]] \cdot (\mf{c}^r)^T) \]
\end{lemma}

\noindent 
{\bf Proof}: The first equality follows from Lemma~\ref{lemma:Tilde}. For the second equality, note that $\rho(\ell_{k,l})(\rho(\ell_{i,j}) (\mf{c}^r))=[\ell_{kl} , [\ell_{ij}, \mf{c}^1]]$. 
Next, since $\{ \mf{c}^i | i=0,\ldots ,n-1 \}$ is an orthogonal basis in the norm $\langle X,Y\rangle =Tr(XY^T)$, we have that the component along $\mf{c}^r$ is precisely $\langle \cdot , \mf{c}^r \rangle /\langle \mf{c}^r ,\mf{c}^r \rangle$. The result follows since $\langle \mf{c}^r ,\mf{c}^r \rangle =Tr(\mf{c}^r (\mf{c}^r)^T )=n$. $\Box $

\begin{lemma}
Let $L_i =\ell \cdot \mf{c}^i=\ell_{0,-i}$ and let:
\[ \Pi (L_i ,L_j )=\sum_{k=0}^{n-1} p_{ij}^k \mf{c}^k \]
Then $p_{ij}^k$ are given by:
\[ p_{ij}^k = \left\{ \begin{array}{rl}
(n-1)-(i+j) & \mbox{ when $k\neq 0$ and  $k=i+j+1$ in $\Z_n$} \\
0 & \mbox{otherwise}
\end{array} \right. \]
\end{lemma}

\noindent 
{\bf Proof}: Let is use the index set $\Z_n =\{ 0,\ldots , n-1\}$. Let us first note that for any $i,j$, there is at most one $k$ for which $p_{ijk}\neq 0$. This is simply because $\Pi (L_i ,L_j )= [L_j ,[L_i ,\mf{c}^1]]$ has support only in $D_{i+j+1}$. Thus, $k=i+j+1$. 
Let us denote the inner term $[\ell \mf{c}^i ,\mf{c}^1]$ as $\Delta_i$ which is given by:
\[ \Delta_i (r,s)= \left\{ \begin{array}{rl}
-1 & \mbox{if $r\in \{ 0,\ldots, n-2\}$ and $s=r+i+1$ in $\Z_n$} \\
n-1& \mbox{$r=n-1$ and $s=r+i+1$ in $\Z_n$} \\
0 & \mbox{otherwise} \\
\end{array} \right.
\]
$L_j$ itself is given by:
\[ L_j (r,s)= \left\{ \begin{array}{rl}
r & \mbox{if $r\in \{0,\ldots, n-1\}$ and $s=r+j$ in $\Z_n$} \\
0 & \mbox{otherwise} \\
\end{array} \right.
\]
Then the matrix of interest $L_j \Delta_i -\Delta_i L_j$. We see that:
\[ \begin{array}{rcl}
Tr(L_j \Delta_i (\mf{c}^{k})^T)&=& \sum_{r=0}^{n-1} (L_j)_{r,r+j} (\Delta_i)_{r+j,r+k}\\
&=& \sum_{r=0, r+j \neq n-1}^{n-1} (-r)+ r\cdot (n-1)\delta_{r+j,n-1} \\
&=& \sum_{r=0}^{n-1} (-r)+ r\cdot (n)\delta_{r+j,n-1} \\
&=& -n(n-1)/2+n(n-1-j) \\
&=& n(n-1)/2-j\cdot n \\
\end{array} 
\]
In the same manner, we have:
\[ \begin{array}{rcl}
Tr(\Delta_i L_j (\mf{c}^{k})^T)&=& \sum_{r=0}^{n-1} (\Delta_i)_{r,r+i+1} (L_j)_{r+i+1,r+k}\\
&=& \sum_{r=0}^{n-2} (\Delta_i)_{r,r+i+1} (L_j)_{r+i+1,r+k}+(n-1)(L_j)_{n+i,n-1+k}\\
&=& -\sum_{r=0}^{n-2-i} r -\sum_{n-2-i+1}^{n-2} (r-n)+(n-1)i\\
&=& -\sum_{r=0}^{n-2} (r+i+1) +\sum_{r=n-2-i+1}^{n-2} (r+i+1-n)+(n-1)i\\
&=& -(i+1+n+i-1)(n-1)/2 +in+i(n-1) \\
&=& -(n-1)n/2 +in \\ 
\end{array} \]
This gives us:
\[ p_{ij}^k = \frac{1}{n} Tr((L_j \Delta_i -\Delta_i L_j)(\mf{c}^k)^T)=(n-1)-(i+j)\]
where $k=i+j+1$. This proves the lemma. $\Box$

\begin{ex}
Let $\Pi^k$ be the matrices $(p^k_{ij})$, then for $n=5$ we have the following.
\[ P^1 =\left[ \begin{array}{rrrrr}
4 & 0 & 0 & 0 & 0 \\
0 & 0 & 0 & 0 & -1 \\
0 & 0 & 0 & -1 & 0\\
0 & 0 & -1 & 0 & 0 \\
0 & -1 & 0 & 0 & 0 \\
\end{array} \right]
P^2 =\left[ \begin{array}{rrrrr}
0 & 3 & 0 & 0 & 0 \\
3 & 0 & 0 & 0 & 0 \\
0 & 0 & 0 & 0 & -2\\
0 & 0 & 0 & -2 & 0 \\
0 & 0 & -2 & 0 & 0 \\
\end{array} \right] \]
\[ P^3 =\left[ \begin{array}{rrrrr}
0 & 0 & 2 & 0 & 0 \\
0 & 2 & 0 & 0 & 0 \\
2 & 0 & 0 & 0 & 0\\
0 & 0 & 0 & 0 & -3 \\
0 & 0 & 0 & -3 & 0 \\
\end{array} \right]
P^4 =\left[ \begin{array}{rrrrr}
0 & 0 & 0 & 1 & 0 \\
0 & 0 & 1 & 0 & 0 \\
0 & 1 & 0 & 0 & 0\\
1 & 0 & 0 & 0 & 0 \\
0 & 0 & 0 & 0 & -4 \\
\end{array} \right]
\]
\end{ex}

Note that $\Pi (L_i ,L_j)$ will have no component along $\mf{c}^0$, the identity matrix, simply because the matrix $[L_i [L_j,\mf{c}^1]]$ will have zero trace. Next, the normal $\mf{c}^1$ is indeed special since, the in any local model of the projective orbit $\mathbb{P}O(\mf{c})$, the normal $\mf{c}$ to $O(\mf{c})$ will collapse. We then have the following important lemma:
\begin{lemma}
For the element $\mf{c} \in Mat_n$, we have $\Pi (L_0 ,L_0)=0$ for the submanifold $\mathbb{P}O(\mf{c}) \subseteq \mathbb{P}V$. Moreover, if $n$ is odd then of all $L_i$, it is unique in having this property. 
\end{lemma}

This analysis does show that the limiting 1-PS is indeed special in the local differential geometry of the starting point $y$. However, the exact connection between this and the algebraic conditions which lead us to the limit $x$, and its properties, is not clear.

\section{Conclusion}
Our motivating question has been to understand the conditions under which a point $x\in V$ is in the projective orbit closure of $y\in V$, when both $x$ and $y$ have distinctive stabilizers, resp. $H$ and $K$, with Lie algebras ${\cal H}$ and ${\cal K}$. Our primary objective has been to provide some initial Lie algebraic tools to analyse this important problem. 

Towards this, we develop a local model for $x$, i.e., an open neighborhood of $x$ and a Lie algebra action of ${\cal G}$ with certain useful properties. 
We show that this open neighborhood has a convenient product structure $M\times N$, where $M$ corresponds to a neighborhood of the orbit $O(x)$ identified with ${\cal S}$, a complement to ${\cal H}$ in ${\cal G}$, and $N$, a complementary space to $T_x O(x)\subseteq V$, the tangent space of the orbit. 
We show that for points $y=x+n \in M\times N$, the Lie algebra ${\cal K}$ of the stabilizer $K$ of $y$ is determined by a subspace of ${\cal H}$. 
Indeed, ${\cal K}$ is the ${\cal S}$-completion of this subspace of ${\cal H}$, where ${\cal S}$ above is a complement to ${\cal H}$ in ${\cal G}$. 
Moreover, we show that the computation of ${\cal K}$ is made effective by a map $\theta $, which measures the ``shear'' due to ${\cal Q}$, the unipotent radical of ${\cal H}$, and its close relative $\Phi$, which is intimately connected with the curvature of the orbit at $x$. This alignment between ${\cal H}$ and ${\cal K}$ is an  extension of Luna's results on reductive stabilizers, but at the Lie algebra level. 

Our main application is on forms, i.e., on the space $V=Sym^d (X)$, and points $x,y\in V$, which have distinctive stabilizers, and where $x$ is obtained from $y$ by taking limits through 1-parameter subgroups $\lambda (t)$ of $GL(X)$. We show that that there is a leading term subalgebra ${\cal K}_0$ of ${\cal K}$, the stabilizer of $y$, which appears as a subalgebra of ${\cal H}$, the stabilizer of $x$. We also show that a generic $\lambda $ leads to solvable ${\cal K}_0$, while non-generic $\lambda$ connect the form $x$, the tangent of approach, the tangent of exit and the form $y$ through triple-stabilizers. 

The connection between ${\cal K}_0$ and ${\cal K}$ is explored further by considering ${\cal K}$ as an extension of ${\cal K}_0$. Under a transversality assumption, there is an extension ${\cal K}(\epsilon)$ of ${\cal K}_0$ whose structure constants agree with that of ${\cal K}$ up to the first order. This leads us to derivations of Lie algebras and certain Lie cohomology conditions. These tools are applied to the co-dimension 1 subvarieties within orbit closures. 

The local model is next applied to the celebrated problem of computing projective closures of orbits of a single matrix under conjugation. Here we reprove the classical result of \cite{gerstenhaber1961dominance, hesselink1979desingularizations} on the orbit closure of a nilpotent matrix, as well as prove, perhaps for the first time,  a result on the projective closures of a general semisimple matrix. 

Finally, we consider the {\em path problem}, of the choice of paths which lead out of $y$ and their limit points $x$.  We show how Kempf's algebraic construction for a torus does have an unconstrained optimization formulation, and thus points to local gradient descent schemes at the local tangent space of $y$. We then connect this with the map $\Phi $ above and compute curvatures at $y$ and show, through an example, that these may indicate good choices of paths out of $y$. 

 Coming back to GCT, our motivation for the problems considered in this paper, identifying linkages between forms and their stabilizer groups has been an important theme in the GCT approach. 
The expectation is that such linkages would lead to obstructions, or witnesses of why a form $g$ with stabilizer $H$, cannot be in the orbit closure of another form, $f$, with stabilizer $K$, where both $H$ and $K$ are distinctive. Indeed, in \cite{mulmuley2001geometric} and \cite{mulmuley2008geometric} the authors proposed representation theoretic obstructions which arose from a generalization of the Peter-Weyl theorem to this situation. These were $G$-modules in coordinate rings of the orbit closure of $g$ and $f$, and the obstruction arose from the  multiplicities of these modules in the respective rings. Whether these multiplicities obstructions will be enough to separate $g$ from $f$ has been intensely analysed. 

Perhaps the local model allows for another source of obstructions connecting ${\cal H}$ and ${\cal K}$, the corresponding Lie algebra stabilizers. If indeed $g$ is obtained as a limit of $f$, then the local model gives us two intermediate structures which link ${\cal H}$ and ${\cal K}$. The first is a ``tangent of approach'', a form $f_b$, along which $f$ approaches $g$, and the tangent of exit $\ell f$. The second is the Lie algebra ${\cal K}_0$, which is obtained from ${\cal K}$ through a limiting process. ${\cal K}_0$ is also a sublagebra of ${\cal H}$, which stabilizes the vector $\overline{f_b}$ in the ``quotient'' representation $V/TO_g$ of ${\cal H}$.  This also requires an alignment in certain Lie algebra cohomology modules, which may be the place for potential obstructions. We believe that these obstructions may be especially computable when the orbit of $g$ is expected to be within the orbit of $f$ as a sub-variety of codimension one.  

\bibliography{references}
\appendix\section{Algebraic complexity and the geometric complexity theory - a brief introduction}
\label{appendix-gct}
Let $Y$ denote a vector of indeterminates suitably arranged, (e.g., as an $m\times m$-matrix) and let $Sym^d (Y)$ be the space of forms of degree $d$.  
A natural model to compute forms are {\em straight-line programs} (also known as {\em arithmetic circuits}). Such a program consists  of an ordered sequence of arithmetic operations involving $+,-,*$ (called lines), starting with variables $y_{i}\in Y$ and constants from ${\mathbb C}$. A form $f \in Sym^d (Y)$ is said to have an efficient computation if there is straight-line program computing $f$ in which the number of lines is bounded by a polynomial in $|Y|$ and $d$.  
The determinant form $det_m (Y) \in Sym^m (Y)$ (where $Y$ is an $m\times m$-matrix) is known be be efficiently computable. On the other hand, for the permanent form, $perm_m (Y)$, no such efficient computation is known.  The analysis of this divergence has been a central problem of the theory.

In a remarkable construction, Valiant\cite{Valiant79} showed that for a form $f(Y) \in Sym^d (Y)$, if there is a program\footnote{In this construction the program corresponds to an {\it arithmetic formula}.} of $L$ lines then there is an $O(L)\times O(L)$-matrix with entries as linear forms $A_{ij} (Y')$ in the variables $Y'=Y\cup \{z\}$ (where $z$ is merely a degree homogenizing variable) such that $det_L (A(Y'))=z^{(L-d)} f(Y)$.

\eat{An operation of great importance in computational complexity is that of reduction, i.e., the use of an efficient algorithm for one form, say $f\in Sym^m(X)$, for the computation of another form, say $g\in Sym^n (Y)$. 
The simplest mechanism to implement a reduction is to substitute for $X$, a linear combination $A(Y)$, i.e., $h(Y)=f(A(Y))\in Sym^m (Y)$. It is clear that such forms $h$, after suitable renaming of variables, are in the projective orbit closure of $f$.

Next, an important result in the theory\cite{Valiant79} is that the determinant is {\em universal}.
More precisely, given any form $g(Y)\in Sym^m (Y)$, there is an enlargement $X\supseteq Y$, an $n \geq m$, and an $A:X\rightarrow Y$ 
such that $g'(X)=det_n (A(Y))$, where $g'(X)$ is a suitable degree homogenized version of $g(X)$. This reduction of $g$ to $det_n$ (or for that matter, to any form $f$), upto homogenization, may be denoted as $g\preceq det_n$ (or resp. $g\preceq f$). 

One of the  celebrated questions of the theory is to show that $perm_m$ that $perm_m$ not reducible to $det_n$, i.e., $perm_m \not \preceq det_n$, when 
$n$ is subexponential in $m$. }

We may thus define the determinantal complexity of a function $f \in Sym^d (Y)$ as the smallest integer $n$ such that $f$ can be written as the determinant of an $n \times n$ matrix whose entries are linear forms in the variables $Y'=Y\cup \{z\}$. Valiant's construction shows  that the determinant is {\em universal}, and that every form $f$ has a finite determinantal complexity.

Let us now connect determinantal complexity with the projective orbit closure of the form $det_n (X)\in Sym^n (X)$. Suppose that the determinantal complexity of $f \in Sym^d (Y)$ is $n$ (naturally, with $n> |Y|$). 
Then (i) define $X=\{ x_{ij}, 1 \leq i, j\leq n\}$ as a set of new $n^2$ variables, (ii) identify $y_{ij}$ with $x_{ij}$ for $1 \leq i,j \leq m$ and $z$ with $x_{nn}$, and finally (iii) define the form $f_{d,n}(X) := x_{nn}^{n-d} f(Y)  \in Sym^n(X)$. The condition that $f_{d,n}(X) = z^{n-d}det(A(Y'))$ tell us that there is a (possibly singular) matrix $A\in End(\C \cdot X)$ such that $det (A\cdot X)=f_{d,n}(X)$.  Whence this homogenized version of $f$ is in the $GL(X)$ projective orbit closure of $det_n$. {\em Thus, computability of forms in various models of computation is intimately connected with the algebraic closures of orbits of universal functions associated with these models.} This theme was first proposed as the Geometric Complexity Theoretic (GCT) approach in the study of computational complexity\cite{mulmuley2001geometric}. 

 A central question in computational complexity has been to determine the determinantal complexity of $perm_m$, and thus to determine those $n$ for which $perm_{m,n} \in \overline{O(det_n )}$. 

It is well known that invariant properties of these functions play an important part in the design of algorithms. For example, the invariance of the determinant to elementary row operations, has paved the way for its efficient evaluation. 
The group of symmetries of the permanent too is large, and also well known. 
In fact, both $perm_m (Y)$ and $det_m (Y)$ as forms in $Sym^m (Y)$ are {\em determined} by their stabilizers within $GL(Y)$ - and yet this divergence in their computational complexity. 

One approach to understanding this divergence is by a careful  examination of their symmetries. More concretely, it is to ask if the question of  $[perm_{m,n} ] \not \in \overline{O([det_n ])}$ can be cast in terms of the group theory,  of linkages, or their absence, between the stabilizers of $perm_{m,n}$ and $det_n $. This broad strategy of finding ``obstructions'' was also proposed as a method in the GCT approach through a sequence of papers, which we review briefly later.

\subsection{Obstructions for membership in orbit closures}

Let us now review a key approach of Geometric Complexity Theory, that of the hunt for obstructions for $[x]$ to be in the orbit closure of $[y]$. This was formulated in \cite{mulmuley2001geometric}, in the context of the permanent and determinant forms, both of which have large and distinctive stabilizers. We outline this briefly. 
Recall that $X$ denotes linear dual of the complex $n^2$-dimensional vector space of $n\times n$-matrices over ${\mathbb C}$.

Let $\Delta [det_n ]$ be the projective orbit closure of $[det_n ]$ in ${\mathbb P}(Sym^n (X))$. If $[perm_{m,n}]$ is in $\Delta[det_n]$ then there would be a $GL(X)$ equivariant embedding of the projective orbit closure of $[perm_{m,n}]$, , $Z_{m,n}$, in $\Delta[det_n]$. This would give rise to a surjective $GL(X)$-morphism from the coordinate ring ${\mathbb C}[\Delta[det_n]]$ of $\Delta[det_n]$, to the coordinate ring 
${\mathbb C}[Z_{m,n}]$ of $Z_{m,n}$. 

The distinctive nature of the stabilizers $K$ of $[det_n ]$ and $H$ of $[perm_{m,n}]$ leads to a complete classification of the $GL(X)$-modules which may appear in the above coordinate rings. Whence, the surjection implies that every irreducible $GL(X)$-representation which occurs in ${\mathbb C}[Z_{m,n}]$ must occur in ${\mathbb C}[\Delta[det_n]]$. An irreducible  representation $\lambda$ of $GL(X)$ which occurs in ${\mathbb C}[Z_{m,n}]$ and does not occur in ${\mathbb C}[\Delta[det_n]]$ is called an occurrence obstruction. Mulmuley and Sohoni conjectured that for every non negative integer $c$ and for infinitely many $m$,  there exist irreducible representations $\lambda$ of $GL(X)$ (with $n=m^c$) which occur in ${\mathbb C}[Z_{m,m^c}]$ but not in ${\mathbb C}[\Delta[det_{m^c}]]$. This conjecture was shown to be false for irreducible representations of $GL(X)$ of shape $n \times d$ in \cite{ikenmeyer2017rectangular}. In \cite{burgisser2019no} the authors showed that the conjecture is false 
for all $\lambda$, when $n \geq m^{25}$. Note that these results do not rule out the GCT approach. If one were to show that the multiplicity of some  irreducible $GL(X)$-module $\lambda$ is more in 
${\mathbb C}[Z_{m,m^c}]$ than in ${\mathbb C}[\Delta[det_{m^c}]]$ then an obstruction has been found. If one were to show this for infinitely many $m$ and every $c$ that would separate the  complexity class $VP$ from $VNP$, the flagship problem of GCT.

Mignon and Ressayre \cite{mignon2004quadratic} showed that the determinantal complexity of $perm_m$ is $\Omega(m^2)$. They consider the zero locus of the determinant  hypersurface $det_n = 0$ in ${\mathbb C}^{n^2}$ and the zero locus of the permanent hypersurface $perm_m = 0$ in ${\mathbb C}^{m^2}$. They show that the Hessian of the determinant at any point on $det_n=0$ has rank at most $2n+1$ and this rank is $m^2$ at a generic point on $perm_m = 0$. Their proof follows from this beautiful calculation. Although their proof does not explicitly use obstructions, Growchow \cite{grochow2015unifying} showed, using a construction of 
$\cite{manivel2013hypersurfaces}$, that their proof can also be seen through the lens of representation theoretic obstructions.

\section{Stabilizer conditions in special cases}
\label{sec:stabcond}
Let's analyse the stabilization conditions of Theorem \ref{theorem:mainform} for an element $\mf{h}\in {\cal K}_0$ in a special case. We have:
\[ \mf{s} \cdot g = \mf{h} \cdot f_b \]

We may split $\mf{s}$ further as $\mf{s}=\mf{s}_{-1}+\mf{s}_0 +\mf{s}_1 $, and $\mf{h}$ as $\mf{h}=\mf{h}_{-1}+\mf{h}_0 +\mf{h}_1$ to give us:
\[ \mf{s}_{-1} \cdot g + \mf{s}_0 \cdot g + \mf{s}_1 \cdot g = \mf{h}_{-1} \cdot f_b + \mf{h}_0 \cdot f_b + \mf{h}_1 f_b \]

Now note that $\mf{s}_0 ,\mf{h}_0$ are elements of $Hom(Y,Y)\oplus Hom(Z,Z)$. Typically, the dimension of $Hom(Y,Y)$ may be much smaller that the dimension of ${\cal K}_0$. Thus, there is a subspace ${\cal K}' \subseteq {\cal K}_0$ of dimension $k-|Y|^2$ where we have that for $\mf{h}\in {\cal K}'$, we have $\mf{h}=\mf{h}_{-1}+\mf{h}_Z +\mf{h}_1$, where $\mf{h}_Z \in Hom(Z,Z)$, and that the stabilizer condition at degree $b$ and $b-1$ is satisfied as follows:
\[ \mf{s}_0 \cdot g + \mf{s}_1 \cdot g= \mf{h}_{-1} \cdot f_b + \mf{h}_Z f_b\]

Noting that $g\in V_a$ while $f_b \in V_b$ gives us the following possibilities, (i) $b=a+1$, (ii) $b=a+2$, and (iii) $b>a+2$. Let us take each of them one by one. 
\begin{enumerate}
    \item {\bf Case $b=a+1$}. We have $\mf{h}_{-1} f_b=\mf{s}_0 \cdot g$ and $\mf{s}_1 g=\mf{h}_Z f_b$. 
    \item {\bf Case $b=a+2$.} We have $\mf{s}_1 g =\mf{h}_{-1} f_b$ and $\mf{h}_Z f_b =0$. 
    \item {\bf Case $b>a+2$.} We have $\mf{s}_1 g =0, \mf{h}_{-1} f_b=0$ and $\mf{h}_Z f_b =0$. 
\end{enumerate}
Each condition gives us certain stringent conditions on the relationship between $g$ and $f_b$. As an example, let us take the case $b=a+1$ and an $\mf{h}\in {\cal K}'$. Let us assume that $|Z|=r$ and $|Y|=s$. 
We may find a basis for $\C \cdot Z$ such that $\mf{h}_Z$ is in Jordan canonical form. Assume for the moment that $\mf{h}_Z$ is diagonal with eigenvalues $\mu_1 ,\ldots ,\mu_r$ where the first $r'$ are non-zero. We may express $f_b $ in this new basis as $f_b =\sum_{\alpha} z_{\alpha} f_b^{\alpha}$, where $z_{\alpha}$ are distinct monomials of degree $a+1$ and $f^{\alpha}_b \in Sym^{d-a-1} (Y)$. By the same token, we may express $g$ as $\sum_{\beta} z_{\beta } g^{\beta}$, where $z^{\beta} \in Sym^a (Z)$ and $g^{\beta} \in Sym^{d-a} (Y)$. The condition $\mf{s}_1 g = \mf{r}_Z b$ gives :
\[ \sum_{\beta} \sum_{i=1}^s L_i (z) \cdot \partial g/\partial y_i = \sum_{\alpha} \mu^{\alpha} z_{\alpha} f^{\alpha}_b \]
This tells us that those ${\alpha}$ such that $\mu^{\alpha}\neq 0$, we have that $f_b^{\alpha}$ are $\C$-linear combinations of $\partial g/\partial y_i$, in other words, the ``{\em minors}'' of $g$. The first condition, {\em viz.}, $\mf{q}f_b =0$ gives us that there are linear relations between the $f^{\alpha}_b$, i.e., the minors of $g$, with coefficients in $Sym^1 (Y)$. These are akin to the Laplace's conditions on the expansion of the determinant.  

\bibliographystyle{alpha}

\end{document}